\documentclass{article}
\usepackage[T1]{fontenc}
\usepackage{lmodern} 
\usepackage[a4paper]{geometry} 
\usepackage{fullpage} 

\usepackage[french]{babel}
\usepackage{amsmath,amsfonts,amssymb,stmaryrd}
\usepackage{enumerate}
\usepackage{amscd}
\usepackage{amsthm}
\usepackage[all]{xy}
\usepackage{dsfont}

\newtheorem{déf}{Définition}[subsection] 
\newtheorem{propr}[déf]{Propriété}
\newtheorem{propo}[déf]{Proposition}
\newtheorem{theo}[déf]{Théorème}

\newtheorem{lem}[déf]{Lemme}
\newtheorem*{theoNonNum}{Théorème} 

\theoremstyle{definition}
\newtheorem{rmq}[déf]{Remarque}

\theoremstyle{definition}

\theoremstyle{definition}
\newtheorem{nota}[déf]{Notation}

\newenvironment{démo}
{\noindent{\textit{Démonstration}.}}{\hfill$\square$}

\begin{document}

\title{Correspondance de Jacquet-Langlands et distinction: cas des représentations  
cuspidales de niveau $0$}
\author{Charlène Coniglio-Guilloton}

\maketitle

\textbf{Résumé:} Soit $\mathbb{K} / \mathbb{F}$ une extension quadratique modérément ramifiée 
de corps locaux non archimédiens. 
Soit ${\rm GL}_m (\mathcal{D})$ une forme intérieure de ${\rm GL}_n (\mathbb{F})$ et 
${\rm GL}_{\mu} (\Delta) = ({\rm M}_m (\mathcal{D}) \otimes_{\mathbb{F}} \mathbb{K})^{\times}$. 
Alors ${\rm GL}_{\mu} (\Delta)$ est une forme intérieure de ${\rm GL}_{n} (\mathbb{K})$ et 
les quotients ${\rm GL}_{\mu} (\Delta) / {\rm GL}_m (\mathcal{D})$ et 
${\rm GL}_{n} (\mathbb{K}) / {\rm GL}_{n} (\mathbb{F})$ sont des espaces symétriques. 
En utilisant la paramétrisation de Silberger et Zink, nous déterminons des critères de 
${\rm GL}_m (\mathcal{D})$-distinction pour les représentations cuspidales de niveau $0$ de 
${\rm GL}_{\mu} (\Delta)$ qui sont l'image d'une cuspidale de niveau $0$ par 
Jacquet-Langlands, puis nous prouvons qu'une représentation cuspidale de niveau $0$ de 
${\rm GL}_n (\mathbb{K})$ est ${\rm GL}_n (\mathbb{F})$-distinguée si et seulement si son image par 
la correspondance de Jacquet-Langlands est ${\rm GL}_m (\mathcal{D})$-distinguée.\\

\textbf{Abstract:} Let $\mathbb{K} / \mathbb{F}$ be a tamely ramified 
quadratic extension of non-archimedean locally compact fields. 
Let ${\rm GL}_m (\mathcal{D})$ be an inner form of ${\rm GL}_n (\mathbb{F})$ and 
${\rm GL}_{\mu} (\Delta) = ({\rm M}_m (\mathcal{D}) \otimes_{\mathbb{F}} \mathbb{K})^{\times}$. 
Then ${\rm GL}_{\mu} (\Delta)$ is an inner form of ${\rm GL}_{n} (\mathbb{K})$ and 
the quotients ${\rm GL}_{\mu} (\Delta) / {\rm GL}_m (\mathcal{D})$ and 
${\rm GL}_{n} (\mathbb{K}) / {\rm GL}_{n} (\mathbb{F})$ are symmetric spaces. 
Using the parametrization of Silberger and Zink, we determine conditions of ${\rm GL}_m (\mathcal{D})$-distinction 
for level zero cuspidal representations of ${\rm GL}_{\mu} (\Delta)$ which are the image of a 
level zero cuspidal representation of ${\rm GL}_n (\mathbb{K})$ by the Jacquet-Langlands correspondence. 
We also show that a level zero cuspidal representation of 
${\rm GL}_n (\mathbb{K})$ is ${\rm GL}_n (\mathbb{F})$-distinguished if and only if its image by 
the Jacquet-Langlands correspondence is ${\rm GL}_m (\mathcal{D})$-distinguished.

\section*{Introduction.}

Nous noterons $\mathbb{K} / \mathbb{F}$ une extension quadratique séparable modérément ramifiée 
de corps locaux non archimédiens. 
On fixe un entier naturel $n$ tel que $n \geq 2$ et on note $G_{\mathbb{F}} = {\rm GL}_n (\mathbb{F})$ 
et $G_{\mathbb{K}} = {\rm GL}_n (\mathbb{K})$. On fixe un diviseur $d$ de $n$ 
($dm = n$) et $\mathcal{D}$ une $\mathbb{F}$-algèbre à division centrale d'indice $d$ 
(i.e de dimension $d^2$ sur son centre $\mathbb{F}$). Enfin, on note 
$H_{\mathbb{F}} = {\rm GL}_m (\mathcal{D})$ et 
$H_{\mathbb{K}} = ({\rm M}_m (\mathcal{D}) \otimes_{\mathbb{F}} \mathbb{K})^{\times}$. 
Il existe un diviseur $\mu$ de $n$ et une $\mathbb{K}$-algèbre à division centrale $\Delta$ 
d'indice $\delta = n / \mu$ tels que $H_{\mathbb{K}} = {\rm GL}_{\mu} (\Delta)$ 
(on remarque que $G_{\mathbb{K}} = H_{\mathbb{K}}$ et $G_{\mathbb{F}} = H_{\mathbb{F}}$ lorsque $d = 1$). 
On a le résultat suivant :

\begin{theoNonNum}
(\cite{JacquetLanglands}, \cite{Rogawski}, \cite{DeligneKazhdanVigneras}, \cite{Badulescu})
Il existe une unique bijection, appelée correspondance de Jacquet-Langlands:
$$
JL : \mathcal{R}^2 (G_{\mathbb{K}}) \rightarrow \mathcal{R}^2 (H_{\mathbb{K}})
$$
telle que pour tous $(g, \widetilde{g}) \in H_{\mathbb{K}} \times G_{\mathbb{K}}$ 
elliptiques réguliers de même polynôme minimal, et toute représentation 
$\pi \in \mathcal{R}^2 (G_{\mathbb{K}})$, on a:
$$
\Theta_{\pi} (\widetilde{g}) 
= (-1)^{\mu \times (\delta-1)} \Theta_{JL (\pi)} (g)
$$
où $\mathcal{R}^2 (G_{\mathbb{K}})$ (resp. $\mathcal{R}^2 (H_{\mathbb{K}})$) sont les 
 classes d'isomorphisme des représentations lisses irréductibles membres de la série discrète 
de $G_{\mathbb{K}}$ (resp. $H_{\mathbb{K}}$) et $\Theta$ désigne le caractère d'Harish Chandra.
\end{theoNonNum}

Lorsque $n = 2$ et $d=2$ alors $G_{\mathbb{K}} = H_{\mathbb{K}} = {\rm GL}_2 (\mathbb{K})$, 
$G_{\mathbb{F}} = {\rm GL}_2 (\mathbb{F})$, $H_{\mathbb{F}} = \mathcal{D}^{\times}$ 
et la correspondance de Jacquet-Langlands est l'application identité. 
Dans ce cas, les travaux de J. Hakim et de D. Prasad nous montrent qu'une série discrète (en niveau quelconque) de 
${\rm GL}_2 (\mathbb{K})$ 
est ${\rm GL}_2 (\mathbb{F})$-distinguée si et seulement si elle est 
$\mathcal{D}^{\times}$-distinguée (on pourra se référer dans \cite{Hakim} au 
Théorème 9.1 page 21 pour les représentations cuspidales de caractère 
central trivial ainsi qu'au Théorème 7.1 page 16 pour la représentation de Steinberg et ses tordues 
ou on pourra retrouver ce résultat dans \cite{Prasad2} Théorème C).\\
Plus généralement, les travaux de J. Hakim et F. Murnaghan dans \cite{HakimMurnaghan1} 
(Théorème 11.1 page 1887) nous donnent des critères de 
${\rm GL}_n (\mathbb{F})$-distinction pour les cuspidales modérées de 
${\rm GL}_n (\mathbb{K})$ (là aussi en niveau quelconque). \\
Dans le travail qui suit, nous généralisons les résultats évoqués précédemment dans le cas des cuspidales de niveau $0$. 
Avec les articles \cite{SilbergerZink1} et \cite{SilbergerZink2}, 
A. Silberger et E. W. Zink montrent que la correspondance 
de Jacquet-Langlands se restreint en une bijection:
$$
JL: \mathcal{R}_0^2 (G_{\mathbb{K}}) \rightarrow \mathcal{R}_0^2 (H_{\mathbb{K}})
$$
où $\mathcal{R}_0^2 (G_{\mathbb{K}})$ (resp. $\mathcal{R}_0^2 (H_{\mathbb{K}})$) sont les classes d'isomorphisme des 
représentations lisses irréductibles membres de la série discrète de niveau $0$.
De plus, les articles \cite{SilbergerZink1} et \cite{SilbergerZink2} nous donnent une paramétrisation de ces séries discrètes 
de niveau $0$ via des paires admissibles modérées et montrent 
que si $\pi \in \mathcal{R}_0^2 (G_{\mathbb{K}})$ est cuspidale de niveau $0$ alors son image par 
la correspondance de Jacquet-Langlands est aussi une représentation cuspidale de niveau $0$. 
Enfin, \cite{SilbergerZink1} et \cite{SilbergerZink2} nous donnent une construction explicite (comme 
induite compacte) de ces représentations à partir de la paire admissible modérée qui leur est associée.\\
Dans un premier temps, nous déterminons des conditions nécessaires et suffisantes de 
${\rm GL}_m (\mathcal{D})$-distinction pour les représentations cuspidales de niveau $0$ de 
${\rm GL}_{\mu} (\Delta)$ qui sont l'image d'une cuspidale de niveau $0$ de 
${\rm GL}_n (\mathbb{K})$ par la correspondance de Jacquet-Langlands. 
Puis, dans un deuxième temps, nous montrons 
qu'une représentation cuspidale 
$\pi \in \mathcal{R}_0^2 (G_{\mathbb{K}})$ est $G_{\mathbb{F}}$-distinguée si et seulement si 
son image par la correspondance de Jacquet-Langlands, $JL (\pi)$, est $H_{\mathbb{F}}$-distinguée, 
généralisant ainsi le résultat d'Hakim (en niveau $0$).\\

Expliquons notre démarche.  
Rappelons qu'une représentation complexe $(\pi, V)$ d'un groupe $G$ est distinguée par un sous-groupe 
$H$ s'il existe une forme linéaire non nulle $\varphi : V \rightarrow \mathbb{C}$ qui est $H$-invariante, ce 
qui revient à dire que l'espace d'entrelacement ${\rm Hom}_{H} (\pi, \mathbb{C})$ est non nul. 
Nous allons utiliser le fait qu'une cuspidale (de niveau $0$) de $G_{\mathbb{K}}$ (resp. $H_{\mathbb{K}}$) est une induite compacte. 
La formule de restriction de Mackey ainsi que la réciprocité de Frobenius 
pour l'induction compacte montrent qu'une représentation 
$\pi = {\rm c-Ind}_{\mathcal{K}}^G \pi_0$ de $G$ est distinguée par un sous-groupe $H$ si et seulement s'il existe 
$s$ dans le double quotient $H \backslash G / \mathcal{K}$ tel que l'espace  
${\rm Hom}_{s^{-1} H s \cap \mathcal{K}} ( \pi_0, \mathbb{C})$ soit non nul. 
Dans notre cas, nous allons restreindre le nombre de doubles classes à étudier en utilisant un résultat de 
\cite{HakimMurnaghan2} (Proposition 5.20 page 94), que nous appliquons dans notre contexte en \ref{ResultatHakimMurnaghan}. 
D'après ce résultat, les seuls éléments $s$ de 
$H \backslash G / \mathcal{K}$ qui contribuent à la distinction vérifient en particulier l'égalité  
$\sigma (s \mathcal{K} s^{-1}) = s \mathcal{K} s^{-1}$ où $\sigma$ est l'élément non trivial du groupe de Galois 
${\rm Gal} (\mathbb{K} / \mathbb{F})$. Dans notre étude, $\mathcal{K}$  
est un sous-groupe ouvert compact modulo le centre maximal 
(par exemple $\mathcal{K} = \mathbb{K}^{\times} {\rm GL}_{\mu} (\mathcal{O}_{\Delta})$), 
fixateur d'un sommet de l'immeuble de Bruhat-Tits de $G$ ($G = G_{\mathbb{K}}$ ou $G = H_{\mathbb{K}}$). 
Enfin, si $X_{\mathbb{K}}$ est l'immeuble de Bruhat-Tits de $G$ ($G = G_{\mathbb{K}}$ ou $G = H_{\mathbb{K}}$), 
$X_{\mathbb{F}}$ celui de $H$ ($H = G_{\mathbb{F}}$ ou $H = H_{\mathbb{F}}$) et 
$j : X_{\mathbb{F}} \rightarrow X_{\mathbb{K}}$ l'injection naturelle entre ces immeubles, les travaux de 
G.Prasad et J.K.Yu dans \cite{PrasadYu}, que nous rappelons en \ref{TheoremeYuPrasad}, 
nous montrent que les sommets de $X_{\mathbb{K}}$ stables par l'action de $\langle \sigma \rangle$ sont exactement 
les sommets de $X_{\mathbb{K}}$ qui sont dans l'image de $X_{\mathbb{F}}$ par $j$.\\
Après avoir appliqué ces deux résultats, il nous reste, la plupart du temps, au plus une double classe à considérer. 
Lorsque nous recherchons des critères de $H_{\mathbb{F}}$-distinction par exemple, 
nous sommes  ramenés à étudier un espace d'entrelacement de la forme 
${\rm Hom}_{s^{-1} H s \cap \mathcal{K}} (\pi_0, \mathbb{C})$,  
où $\pi_0$ est l'inflation à $\mathcal{K} = \mathbb{K}^{\times} {\rm GL}_{\mu} (\mathcal{O}_{\Delta})$ 
d'une cuspidale $\overline{\gamma}$ de ${\rm GL}_{\mu} (k_{\Delta})$, $k_{\Delta}$ étant le corps résiduel de $\Delta$. 
Pour cela, nous utiliserons à plusieurs reprises les résultats des travaux de Lusztig dans \cite{Lusztig},  
que nous rappelons en \ref{TheoremeLusztig}, et qui nous permettent  
de calculer la dimension de l'espace d'entrelacement 
${\rm Hom}_{S} (\overline{\gamma}, \mathbb{C})$, où $({\rm GL}_{\mu} (k_{\Delta}), S, \theta)$ est un espace symétrique 
(ce qui signifie en particulier que l'application 
$\theta : {\rm GL}_{\mu} (k_{\Delta}) \rightarrow {\rm GL}_{\mu} (k_{\Delta})$ 
est une involution et que $S$ est contenu dans l'ensemble des points fixes de 
${\rm GL}_{\mu} (k_{\Delta})$ par $\theta$).\\

Dans le travail qui suit, nous commençons dans la partie 
\ref{PartieInjectionImmeubles} par quelques rappels sur 
les injections d'immeubles, injections que nous rendons explicites en \ref{InjectionImmeubleNRImpair}, 
\ref{InjectionImmeubleTRImpair}, \ref{InjectionImmeubleNRPair} et \ref{InjectionImmeubleTRPair} 
(nous distinguons les cas où l'extension $\mathbb{K} / \mathbb{F}$ est non ramifiée, totalement ramifiée, 
et où $d$ est pair ou impair).\\ 
On en déduit dans les parties \ref{PartieCNSDistinctionImpair} et \ref{PartieCNSDistinctionPair}
des conditions nécessaires et suffisantes de ${\rm GL}_m (\mathcal{D})$-distinction 
pour les cuspidales de niveau $0$ de ${\rm GL}_{\mu} (\Delta)$ qui sont 
l'image d'une cuspidale par Jacquet-Langlands, conditions que nous pouvons lire 
sur la paire admissible modérée associée à la représentation. Explicitons ces conditions. 
Pour cela, fixons $\overline{\mathbb{K}}$ une clôture algébrique de $\mathbb{K}$ et pour tout entier naturel non nul $l$, 
notons $\mathbb{K}_l$ l'extension non ramifiée de degré $l$ de $\mathbb{K}$ contenue dans $\overline{\mathbb{K}}$ 
et $k_{\mathbb{K},l}$ son corps résiduel. On fixe $\varpi_{\mathbb{K}}$ une uniformisante de $\mathbb{K}$. 
On note $k_{\mathcal{D}}$ (resp. $k_{\Delta}$) le corps résiduel de $\mathcal{D}$ (resp. $\Delta$). 
On montre le théorème suivant en \ref{CritereDistinctionNRdImpair}, 
\ref{CritereDistinctionTRdImpair}, \ref{CritereDistinctionNRdPair} et \ref{CritereDistinctionTRPair} :

\begin{theoNonNum}
Soit $\pi \in \mathcal{R}_0^2 ({\rm GL}_{\mu} (\Delta))$ une cuspidale de niveau $0$, 
image d'une cuspidale de niveau $0$ de ${\rm GL}_n (\mathbb{K})$ par la correspondance de Jacquet-Langlands, et 
$(\mathbb{K}_n / \mathbb{K}_{\delta}, \chi)$ la paire admissible modérée associée à $\pi$ 
(en particulier $\chi$ est un caractère modéré de $\mathbb{K}_n^{\times}$). 
On note $\overline{\chi}$ 
la restriction de $\chi$ à $\mathcal{O}_{\mathbb{K}_n}^{\times}$ vue comme caractère de $k_{\mathbb{K},n}^{\times}$.
\begin{itemize}
\item[$\ast$] Si $\mathbb{K} / \mathbb{F}$ est non ramifiée, et $n$ est pair, la représentation $\pi$ n'est pas 
${\rm GL}_m (\mathcal{D})$-distinguée.
\item[$\ast$] Si $\mathbb{K} / \mathbb{F}$ est totalement ramifiée, et $n$ est impair, la représentation $\pi$ n'est pas  
${\rm GL}_m (\mathcal{D})$-distinguée.
\item[$\ast$] Si $\mathbb{K} / \mathbb{F}$ est non ramifiée, et $n$ est impair. 
Soit $\tau$ un générateur du groupe de Galois ${\rm Gal} (k_{\mathbb{K},n} / k_{\mathcal{D}})$. 
Alors, la représentation $\pi$ est 
${\rm GL}_m (\mathcal{D})$-distinguée si et seulement si $\chi$ est trivial sur $\mathbb{F}^{\times}$ et s'il existe 
$\alpha$ dans le groupe de Galois ${\rm Gal} (k_{\mathbb{K},n} / k_{\Delta})$ tel que 
$\overline{\chi}^{-1} \circ \alpha = \overline{\chi} \circ \tau$.
\item[$\ast$] Si $\mathbb{K} / \mathbb{F}$ est totalement ramifiée, et $n$ est pair. 
Soit $l_0$ le corps résiduel de $\mathbb{K}_{n/2}$. 
On fixe $\varpi_{\mathbb{K}}$ et $\varpi_{\mathbb{F}}$ telle que $\varpi_{\mathbb{K}}^2 = \varpi_{\mathbb{F}}$. 
Soit $\eta$ dans $k_{\mathbb{K},n}^{\times} \backslash l_0^{\times}$ tel que $\eta^2 \in l_0^{\times}$. 
Alors, la représentation $\pi$ est 
${\rm GL}_m (\mathcal{D})$-distinguée si et seulement si $\chi$ est trivial sur $\mathbb{F}^{\times}$, 
$\overline{\chi}$ est trivial sur $l_0^{\times}$ et 
$\chi (\varpi_{\mathbb{K}}) \overline{\chi} (\eta) = -1$.
\end{itemize}

\end{theoNonNum}

On pourra remarquer que pour  
$\mu = m =n$ et $d = 1$, on retrouve les critères de ${\rm GL}_n (\mathbb{F})$-distinction des cuspidales (de niveau $0$) 
de ${\rm GL}_n (\mathbb{K})$ de \cite{HakimMurnaghan1}.\\ 

Dans les parties \ref{PartieCNSDistinctionImpair} et \ref{PartieCNSDistinctionPair}, 
on montre également le résultat de multiplicité $1$ suivant :  

\begin{theoNonNum}
Si $\pi \in \mathcal{R}_0^2 ({\rm GL}_{\mu} (\Delta))$ est cuspidale de niveau $0$ et est 
l'image d'une cuspidale de niveau $0$ de ${\rm GL}_n (\mathbb{K})$ par la correspondance de Jacquet-Langlands, 
on a:
$$
{\rm dim} ({\rm Hom}_{{\rm GL}_m (\mathcal{D})} (\pi, \mathbb{C})) \leq 1
$$
\end{theoNonNum}

\noindent 
\textbf{Remarque:} Notons que de tels résultats de multiplicité $1$ sont connus pour ${\rm GL}_n (\mathbb{K})$, 
on pourra se référer à \cite{Flicker} pour vérifier que 
$({\rm GL}_n (\mathbb{K}), {\rm GL}_n (\mathbb{F}))$ est une paire de Gelfand. 
La preuve de ce résultat peut s'étendre à 
$({\rm GL}_{\mu} (\Delta), {\rm GL}_m (\mathcal{D}))$. Néanmoins, dans notre travail, 
nous n'avons pas eu besoin d'utiliser ces propriétés.\\

Enfin, dans la partie \ref{PartieConclusion}, nous démontrons le théorème suivant :

\begin{theoNonNum}
Si $\rho \in \mathcal{R}_0^2 ({\rm GL}_n (\mathbb{K}))$ est une cuspidale (de niveau $0$), alors 
$\rho$ est ${\rm GL}_n (\mathbb{F})$-distinguée si et seulement si 
$JL(\rho)$ est ${\rm GL}_m (\mathcal{D})$-distinguée. 
\end{theoNonNum}

\textbf{Remerciements :} Je tiens tout particulièrement à remercier Paul Broussous et Nadir Matringe pour 
de nombreuses discussions très enrichissantes et instructives ainsi que pour 
leurs précieux conseils lors de l'élaboration et la rédaction de ce travail. 
 

\section{Notations et résultats préliminaires.}

\begin{nota}
  On note $\mathbb{K} / \mathbb{F}$ une extension 
quadratique séparable modérément ramifiée de corps locaux non archimédiens, 
$\sigma$ le générateur de ${\rm Gal} (\mathbb{K} / \mathbb{F})$ 
et on fixe un entier naturel $n$ tel que $n \geq 2$.
Soit $H_{\mathbb{F}}$ une forme intérieure de $G_{\mathbb{F}} = {\rm GL}_n (\mathbb{F})$, 
alors il existe un diviseur $d$ de $n$, 
il existe une $\mathbb{F}$-algèbre à division centrale $\mathcal{D}$ d'indice $d$, tels que 
$H_{\mathbb{F}} = A^{\times}$ avec $A = {\rm M}_m (\mathcal{D})$ une $\mathbb{F}$-algèbre centrale simple 
de degré réduit~$n$. On a donc $n = m \times d$.
Soit ${\rm inv}_{\mathbb{F}} (\mathcal{D}) = \frac{r}{d} \in \mathbb{Q} / \mathbb{Z}$, 
où $(r,d) = 1$, l'invariant de Hasse de $\mathcal{D}$.
Soit $A_{\mathbb{K}} = A \otimes_{\mathbb{F}} \mathbb{K}$, alors $A_{\mathbb{K}}$ est 
une $\mathbb{K}$-algèbre centrale simple. Il existe donc  
un entier naturel $\mu$ non nul et une $\mathbb{K}$-algèbre à division centrale $\Delta$, tels que 
$A_{\mathbb{K}} \simeq {\rm M}_{\mu} (\Delta)$.
Notons $\delta = \sqrt{[\Delta : \mathbb{K}]}$. Alors $n= \mu \times \delta$ et
$H_{\mathbb{K}} = A_{\mathbb{K}}^{\times} \simeq {\rm GL}_{\mu} (\Delta)$ 
est une forme intérieure de $G_{\mathbb{K}} = {\rm GL}_n (\mathbb{K})$.
De plus, on peut voir $H_{\mathbb{F}}$ comme sous-groupe de $H_{\mathbb{K}}$ via l'injection suivante:
$$
H_{\mathbb{F}} = A^{\times} \rightarrow H_{\mathbb{K}} = (A \otimes_{\mathbb{F}} \mathbb{K})^{\times}, \, 
g \mapsto g \otimes 1
$$
On remarque que, si $s/\delta$ est l'invariant de Hasse de $\Delta$ (avec $s$ et $\delta$ premiers entre eux), 
on a :
\begin{eqnarray*}
\frac{s}{\delta} = {\rm inv}_{\mathbb{K}} (\Delta)
& = & {\rm inv}_{\mathbb{K}} ({\rm M}_{\mu} (\Delta)) = {\rm inv}_{\mathbb{K}} (A \otimes_{\mathbb{F}} \mathbb{K}) 
    = [ \mathbb{K} : \mathbb{F} ] \times {\rm inv}_{\mathbb{F}} (A)\\
& = &  2 \times \frac{r}{d} = \frac{2r / (d,2)}{d / (d,2)}  
\end{eqnarray*}
Par suite $\delta = \frac{d}{(d,2)}$. 
\end{nota}

\begin{nota}
Si $\mathbb{L}$ est un corps local non archimédien et $\overline{\mathbb{L}}$ une clôture algébrique de 
$\mathbb{L}$, on note $k_{\mathbb{L}}$ le corps résiduel de $\mathbb{L}$. Pour tout entier naturel non nul $r$, 
on note $\mathbb{L}_r$ l'extension non ramifiée de degré $r$ de 
$\mathbb{L}$ contenue dans $\overline{\mathbb{L}}$ et $k_{\mathbb{L},r}$ son corps résiduel 
(extension de degré $r$ de $k_{\mathbb{L}}$).
On notera $\varpi_{\mathbb{L}}$ une uniformisante de $\mathbb{L}$, $\mathcal{O}_{\mathbb{L}}$ 
l'anneau des entiers de $\mathbb{L}$, $\mathcal{P}_{\mathbb{L}}$ son idéal de valuation et 
$v_{\mathbb{L}}$ une valuation normalisée. On utilisera le même type de notations si 
$\mathcal{B}$ est une $\mathbb{L}$-algèbre à division centrale \mbox{($\varpi_{\mathcal{B}}$, $v_{\mathcal{B}}$, $\cdots$)} 
et on notera ${\rm Nrd}_{\mathcal{B}}$ la norme réduite de $\mathcal{B}$.
Dans toute la suite, on fixe $\overline{\mathbb{K}}$ une clôture algébrique de $\mathbb{K}$.
On notera $k_{\mathcal{D}} = k_{\mathbb{F}, d}$ (resp. $k_{\Delta} = k_{\mathbb{K}, \delta}$) 
le corps résiduel de $\mathcal{D}$ (resp. $\Delta$). 
\end{nota}

\begin{nota}
On dira qu'une représentation complexe lisse et irréductible de $G_{\mathbb{K}}$ (resp. $H_{\mathbb{K}}$) est 
membre de la série discrète si elle est essentiellement de carré intégrable modulo le centre. 
On dira qu'une telle représentation est de niveau $0$ si elle possède un vecteur fixe non 
nul sous l'action du sous-groupe $I_n + \varpi_{\mathbb{K}} {\rm M}_n (\mathcal{O}_{\mathbb{K}})$ 
(resp. $I_{\mu} + \varpi_{\Delta} {\rm M}_{\mu} (\mathcal{O}_{\Delta})$). 
On notera $\mathcal{R}_0^2 (G_{\mathbb{K}})$ (resp. $\mathcal{R}_0^2 (H_{\mathbb{K}})$) les classes d'isomorphisme des 
représentations lisses irréductibles membres de la série discrète de niveau $0$ de $G_{\mathbb{K}}$ 
(resp. $H_{\mathbb{K}}$). 
\end{nota}

D'après les articles \cite{SilbergerZink1} et \cite{SilbergerZink2}, 
on sait que l'on a une bijection naturelle, appelée la correspondance 
de Jacquet-Langlands:
$$
JL: \mathcal{R}_0^2 (G_{\mathbb{K}}) \rightarrow \mathcal{R}_0^2 (H_{\mathbb{K}})
$$

\begin{déf}
\begin{itemize}         
\item[i)] Si $\chi$ est un caractère modéré de $\mathbb{K}_n^{\times}$ 
(i.e trivial sur $1+ \mathcal{P}_{\mathbb{K}_n}$), on notera 
$\overline{\chi}$ sa restriction à $\mathcal{O}_{\mathbb{K}_n}^{\times}$ vue 
comme caractère de $k_{\mathbb{K},n}^{\times}$.
On dira qu'un tel caractère est $k_{\mathbb{K},l}$-régulier si son orbite sous l'action du groupe de Galois 
${\rm Gal} (k_{\mathbb{K},n} / k_{\mathbb{K},l})$ est de longueur maximale $n/l$.
On dira de même qu'un caractère de $\mathbb{K}_n^{\times}$ est $\mathbb{K}_l$-réguliers 
si son orbite sous l'action du groupe de Galois ${\rm Gal} (\mathbb{K}_n / \mathbb{K}_l)$ est de longueur maximale.
\item[ii)] Si $\chi_1$ et $\chi_2$ sont deux caractères modérés de $\mathbb{K}_n^{\times}$, on notera:
$$
\overline{\chi_1} \simeq_{k_{\mathbb{K},l}} \overline{\chi_2}  
$$
lorsque $\overline{\chi_1}$ et $\overline{\chi_2}$ sont dans la même 
${\rm Gal} (k_{\mathbb{K},n} / k_{\mathbb{K},l})$-orbite.
\item[iii)] On appelle paire admissible modérée sur $\mathbb{K}_l$ la donnée d'une extension non ramifiée 
$\mathbb{K}_n$ de $\mathbb{K}_l$ et d'un caractère modéré $\chi$ de $\mathbb{K}_n^{\times}$ qui soit 
$\mathbb{K}_l$-régulier.
\end{itemize}
\end{déf}

\begin{propo}\label{JacquetLanglandsCuspidales}
Soit $JL : \mathcal{R}_0^{2} ({\rm GL}_n (\mathbb{K})) \rightarrow  \mathcal{R}_0^{2} ({\rm GL}_{\mu} (\Delta))$ 
la correspondance de Jacquet-Langlands.
Si $\rho \in \mathcal{R}_0^{2} ({\rm GL}_n (\mathbb{K}))$ est cuspidale, 
alors le paramètre de Silberger et Zink associé est une paire admissible modérée sur $\mathbb{K}$, 
c'est-à-dire 
un caractère modéré $\mathbb{K}$-régulier de $\mathbb{K}_n^{\times}$ que l'on notera $\chi$. 
Alors le type étendu maximal (construit par Silberger et Zink) contenu dans $\rho$ est aussi paramétré par $\chi$.\\
De plus, soit $\pi_{\chi} = \pi = JL (\rho)$. Alors $\pi$ 
est également cuspidale (dans ce cas $\overline{\chi}$ est également un caractère 
$k_{\Delta}$-régulier de $k_{\mathbb{K},n}^{\times}$) et le type étendu maximal contenu dans $\pi$ 
est lui aussi paramétré par le caractère $\chi$.
\end{propo}

\begin{démo}
Pour la preuve, on pourra se référer à \cite{SilbergerZink2} (exemples 0.10 page 184) dont on reprend ici les notations. 
Le paramètre associé au type étendu maximal contenu dans $\rho$ est $w^{n - (f,n)} \chi = \chi$ (où 
$f = n$ est la longueur de l'orbite de $\chi$ sous l'action du groupe de 
Galois ${\rm Gal} (\mathbb{K}_n / \mathbb{K})$ et $w$ est un caractère quadratique).\\
Pour la représentation $\pi$, on a $e = n/f = 1$ et $e^{'} = (e, \mu) = 1$ donc $\pi$ est bien cuspidale. 
De plus, le paramètre associé au type étendu maximal contenu dans $\pi$ est $\widetilde{w}^{\mu - (f,\mu)} \chi = \chi$ 
où $\widetilde{w}$ est un caractère quadratique (ici $f = n$ est divisible par $\mu$).

\end{démo}

\begin{nota}\label{NotationsRepresentationCuspidaleParametreSilbZink}
On fixe $\pi$ dans $\mathcal{R}_0^2 ({\rm GL}_{\mu} (\Delta))$, une représentation cuspidale de niveau $0$ 
de ${\rm GL}_{\mu} (\Delta)$ qui est l'image d'une cuspidale de niveau $0$ de 
${\rm GL}_n (\mathbb{K})$ par la correspondance de Jacquet-Langlands.
Alors le paramètre de Silberger et Zink associé à $\pi$ est une paire admissible modérée 
$(\mathbb{K}_n / \mathbb{K}_{\delta}, \chi)$. 
En particulier, le caractère $\overline{\chi}$ est un caractère $k_{\Delta}$-régulier de 
$k_{\mathbb{K},n}^{\times}$.  
Soit $\overline{\gamma}_0$ la représentation cuspidale de 
${\rm GL}_{\mu} (k_{\Delta})$ obtenue par paramétrisation de Green à partir du caractère $\overline{\chi}$ 
(on rappelle ce qu'est la paramétrisation de Green en \ref{ParametrisationGreen}).
Soit $\gamma_0$ la représentation $\overline{\gamma}_0$ vue comme représentation de 
${\rm GL}_{\mu} (\mathcal{O}_{\Delta})$.
Soit $\pi_0$ la représentation de $\mathbb{K}^{\times} {\rm GL}_{\mu} (\mathcal{O}_{\Delta})$ telle que 
pour tout $x$ dans ${\rm GL}_{\mu} (\mathcal{O}_{\Delta})$ :
$$
\pi_0 (x) = \gamma_0 (x)
$$
et pour tout $x$ dans $\mathbb{K}^{\times}$ :
$$
\pi_0 (x) = \chi (x) {\rm Id}
$$
Alors, d'après les articles \cite{SilbergerZink1} et \cite{SilbergerZink2} :
$$
\pi \simeq 
{\rm c-Ind}_{\mathbb{K}^{\times} {\rm GL}_{\mu} (\mathcal{O}_{\Delta})}^{{\rm GL}_{\mu} (\Delta)} (\pi_0)
$$
\end{nota}

\begin{nota}
Si $\rho \in \mathcal{R}_0^2 ({\rm GL}_n (\mathbb{K}))$ est une cuspidale de niveau $0$, alors 
le paramètre de Silberger et Zink associé à $\rho$ est une paire admissible modérée sur $\mathbb{K}$, 
$(\mathbb{K}_n / \mathbb{K}, \chi_0)$. 
Notons $\overline{\lambda}_0$ la représentation cuspidale de 
${\rm GL}_{n} (k_{\mathbb{K}})$ obtenue par paramétrisation de Green à partir du caractère $\overline{\chi_0}$ 
et $\lambda_0$ la représentation $\overline{\lambda}_0$ vue comme représentation de 
${\rm GL}_{n} (\mathcal{O}_{\mathbb{K}})$.
Soit $\rho_0$ la représentation de $\mathbb{K}^{\times} {\rm GL}_{n} (\mathcal{O}_{\mathbb{K}})$ telle que 
pour tout $x$ dans ${\rm GL}_{n} (\mathcal{O}_{\mathbb{K}})$, on a 
$\rho_0 (x) = \lambda_0 (x)$, 
et pour tout $x$ dans $\mathbb{K}^{\times}$, 
$\rho_0 (x) = \chi_0 (x) {\rm Id}$. 
Alors, d'après les articles \cite{SilbergerZink1} et \cite{SilbergerZink2} :
$$
\rho \simeq 
{\rm c-Ind}_{\mathbb{K}^{\times} {\rm GL}_{n} (\mathcal{O}_{\mathbb{K}})}^{{\rm GL}_{n} (\mathbb{K})} (\rho_0)
$$
Et $\pi = JL(\rho)$ a pour paramètre de Silberger et Zink la paire admissible modérée sur $\mathbb{K}_{\delta}$, 
qui est aussi une paire admissible modérée sur $\mathbb{K}$, 
$(\mathbb{K}_n / \mathbb{K}, \chi_0)$. 
\end{nota}

\section{Injection des immeubles de Bruhat-Tits.}\label{PartieInjectionImmeubles}

Pour les chaînes et les ordres de chaînes, on se réfère à \cite{BushnellFröhlich} pages 212 à 224. 
On introduit pour toute la suite les notations suivantes :

\begin{nota}\label{NotationsImmeuble}
On notera $X_{\mathbb{F}}$ (resp. $X_{\mathbb{K}}$) l'immeuble de Bruhat-Tits de 
${\rm GL}_m (\mathcal{D})$ (resp. ${\rm GL}_{\mu} (\Delta)$).
On identifiera les sommets de $X_{\mathbb{F}}$ (resp. $X_{\mathbb{K}}$) aux 
$\mathcal{O}_{\mathcal{D}}$-chaînes de réseaux (resp. $\mathcal{O}_{\Delta}$-chaînes de réseaux)
de période $1$ d'un $\mathcal{D}$-espace vectoriel (à droite) de dimension $m$ 
(resp. d'un $\Delta$-espace de dimension $\mu$).
Si $L$ est un $\mathcal{O}_{\mathcal{D}}$-réseau d'un $\mathcal{D}$-espace vectoriel de dimension $m$, 
on notera $[L] = \{ L \varpi_{\mathcal{D}}^k: k \in \mathbb{Z} \}$ le sommet associé.\\
Si $\mathcal{L} = (L_k)_{k \in \mathbb{Z}}$ est une chaîne d'$\mathcal{O}_{\mathcal{D}}$-réseaux 
d'un $\mathcal{D}$-espace vectoriel $V$ de dimension $m$, on notera $\mathcal{A} = \mathcal{A} (\mathcal{L})$ 
l'ordre de chaîne associé,  
$\mathcal{A} (\mathcal{L}) = 
\{ a \in {\rm End}_{\mathcal{D}} (V) : \forall k \in \mathbb{Z}, 
a L_k \subseteq L_k \}$,  
et $\mathcal{K} (\mathcal{L}) = \{ g \in ({\rm End}_{\mathcal{D}} (V))^{\times} : 
g \mathcal{A} (\mathcal{L}) g^{-1} = \mathcal{A} ({\mathcal{L}})\}$ son normalisateur dans 
$({\rm End}_{\mathcal{D}} (V))^{\times}$.
On utilisera des notations similaires pour les chaînes d'$\mathcal{O}_{\Delta}$-réseaux.
\end{nota}

On pourra se référer à \cite{Tits} (2.6 page 47) pour les résultats suivants :

\begin{theo} (\cite{Tits})
Il existe une injection naturelle $j : X_{\mathbb{F}} \rightarrow X_{\mathbb{K}}$
vérifiant les trois propriétés suivantes:
\begin{itemize}
\item[a)] L'application $j$ est ${\rm GL}_m (\mathcal{D})$-équivariante, c'est-à-dire que :
$$ 
j(g.x) = g. j(x)
$$
pour tout $g$ dans ${\rm GL}_m (\mathcal{D})$ et tout $x$ dans $X_{\mathbb{F}}$.

\item[b)] L'image de $j$ est incluse dans $ X_{\mathbb{K}}^{{\rm Gal} (\mathbb{K} / \mathbb{F})}$, où 
$ X_{\mathbb{K}}^{{\rm Gal} (\mathbb{K} / \mathbb{F})}$ désigne 
l'ensemble des éléments de $X_{\mathbb{K}}$ qui sont fixes sous l'action du groupe de Galois 
${\rm Gal} (\mathbb{K} / \mathbb{F})$.

\item[c)] L'application $j$ est affine, c'est-à-dire que pour tout appartement $\mathcal{A}$ de 
$X_{\mathbb{F}}$ (vu comme espace affine), il existe un appartement $\mathcal{B}$ de $X_{\mathbb{K}}$ tel que 
$j(\mathcal{A}) \subseteq \mathcal{B}$ et 
$j_{\vert \mathcal{A}} : \mathcal{A} \rightarrow \mathcal{B}$ est une application affine.
\end{itemize}
\end{theo}

Nous allons utiliser le théorème suivant démontré dans \cite{PrasadYu} (Théorème 1.9 page 555) :

\begin{theo}\label{TheoremeYuPrasad} (\cite{PrasadYu})
Soit $G$ un groupe réductif défini sur un corps local non archimédien $\mathbb{L}$ de 
caractéristique résiduelle $p$. Soit $F \subseteq {\rm Aut}_{\mathbb{L}} (G)$ un 
groupe fini d'ordre non divisible par $p$. On note $G^{F}$ les points fixes de $G$ sous l'action de $F$ 
et $H = (G^{F})^{\circ}$ la composante de l'unité de $G^{F}$. 
Soient $X_G$ et $X_H$ les immeubles de Bruhat-Tits de $G$ et de $H$ respectivement. 
Alors, via l'injection naturelle des immeubles $j : X_H \hookrightarrow X_G$, on peut identifier 
$X_H$ à $X_G^{F}$, l'ensemble des points fixes de $X_G$ sous l'action de $F$.

\end{theo}

On en déduit le résultat suivant :

\begin{propo}
L'application $j$ est unique et 
${\rm Im} (j) = X_{\mathbb{K}}^{{\rm Gal} (\mathbb{K} / \mathbb{F})}$.
\end{propo}

\begin{démo}
L'égalité ${\rm Im} (j) = X_{\mathbb{K}}^{{\rm Gal} (\mathbb{K} / \mathbb{F})}$  
découle directement de \ref{TheoremeYuPrasad} car l'extension $\mathbb{K}/ \mathbb{F}$ est 
modérément ramifiée, donc de caractéristique résiduelle distincte de $2$. 
Pour l'unicité, supposons que l'on ait une autre application 
$l : X_{\mathbb{F}} \rightarrow X_{\mathbb{K}}$ vérifiant les trois conditions a), b) et c).
Comme $l(X_{\mathbb{F}}) \subseteq X_{\mathbb{K}}^{{\rm Gal} (\mathbb{K} / \mathbb{F})} 
= j (X_{\mathbb{F}})$, on peut considérer l'application :
$$
k = j^{-1} \circ l: X_{\mathbb{F}} \rightarrow X_{\mathbb{F}}
$$
Alors, par composition, $k$ est également ${\rm GL}_m (\mathcal{D})$-équivariante et est une application affine.\\
Pour tout sommet $s$ de $X_{\mathbb{F}}$, notons $\mathcal{K}_s$ son stabilisateur. 
Alors, pour tout $g \in \mathcal{K}_s, g.s = s$, donc $g.k(s) = k(g.s) = k(s)$.
Ainsi, $\mathcal{K}_s$ est aussi le stabilisateur de $k(s)$, d'où $k(s) = s$.\\
Finalement, $k$ est l'identité sur les sommets. Le fait qu'elle soit affine nous permet de conclure 
que $k=Id$, et par suite $j=l$.

\end{démo}

\begin{propo}\label{PropositionCaracterisationInjection}
Soient $A_{\mathbb{F}}$ et $A_{\mathbb{K}}$ des appartements des immeubles de 
Bruhat Tits $X_{\mathbb{F}}$ et $X_{\mathbb{K}}$ respectivement. 
Supposons donnée une application affine injective 
$j : A_{\mathbb{F}} \hookrightarrow A_{\mathbb{K}}$. 
Soit $T$ le tore associé à l'appartement $A_{\mathbb{F}}$ et 
$N(T)$ son normalisateur dans ${\rm GL}_m (\mathcal{D})$.  
On suppose que $j$ est $N(T)$-équivariante, qu'il existe un sommet $s_0$ de $A_{\mathbb{F}}$ 
tel que $j(s_0)$ est fixe par l'action du 
groupe de Galois ${\rm Gal} (\mathbb{K} / \mathbb{F})$ et:
$$
{\rm Stab}_{{\rm GL}_m (\mathcal{D})} (s_0) \subseteq 
{\rm Stab}_{{\rm GL}_{\mu} (\Delta)} (j(s_0))
$$
Puis on définit $j : X_{\mathbb{F}} \rightarrow X_{\mathbb{K}}$ par $j(g.x) = g.j(x)$ pour tout $g$ dans 
${\rm GL}_m (\mathcal{D})$ et tout $x$ dans $A_{\mathbb{F}}$.\\
Alors l'application $j$ est bien définie et est l'injection naturelle entre les deux immeubles.

\end{propo}

\begin{démo}
Nous devons vérifier que $j$ est bien définie sur 
$X_{\mathbb{F}} = \bigcup_{g \in {\rm GL}_m (\mathcal{D})} g.A_{\mathbb{F}}$ 
(alors elle sera bien ${\rm GL}_m (\mathcal{D})$-équivariante), que $j$ 
est affine et que pour tout $x$ dans $X_{\mathbb{F}}$, on a:
$$
{\rm Stab}_{{\rm GL}_m (\mathcal{D})} (x) \subseteq {\rm Stab}_{{\rm GL}_{\mu} (\Delta)} (j(x)) \, \, 
\text{et} \, \, j(x) \in X_{\mathbb{K}}^{{\rm Gal} (\mathbb{K} / \mathbb{F})}
$$

\begin{itemize}
\item[$\ast$] Nous allons commencer par montrer que pour tout $x$ dans $A_{\mathbb{F}}$, 
$j(x)$ est stable sous l'action du groupe de Galois 
${\rm Gal} (\mathbb{K} / \mathbb{F})$ et 
${\rm Stab}_{{\rm GL}_m (\mathcal{D})} (x) \subseteq {\rm Stab}_{{\rm GL}_{\mu} (\Delta)} (j(x))$.
Soit $x$ dans $A_{\mathbb{F}}$. 
\begin{itemize}
\item[$\ast$] Considérons tout d'abord le cas où $x$ est un sommet de $A_{\mathbb{F}}$. 
Comme $N(T)$ agit transitivement sur les sommets de $A_{\mathbb{F}}$, il existe $n$ dans $N(T)$ 
tel que $x = n.s_0$.
De plus, comme $j$ est $N(T)$-équivariante sur $A_{\mathbb{F}}$, $j(x)  =j(n.s_0) = n. j(s_0)$, ainsi:
\begin{eqnarray*}
{\rm Stab}_{{\rm GL}_m (\mathcal{D})} (x) 
& = & {\rm Stab}_{{\rm GL}_m (\mathcal{D})} (n.s_0) \\
& = & n {\rm Stab}_{{\rm GL}_m (\mathcal{D})} (s_0) n^{-1} 
\subseteq n {\rm Stab}_{{\rm GL}_{\mu} (\Delta)} (j(s_0)) n^{-1} = {\rm Stab}_{{\rm GL}_{\mu} (\Delta)} (j(x)) 
\end{eqnarray*}
De plus, si $\gamma \in {\rm Gal} (\mathbb{K} / \mathbb{F})$, on a 
$\gamma (j(x)) = \gamma (n) \gamma (j(s_0))$ 
or $\gamma (n) = n$ car $n \in {\rm GL}_m (\mathcal{D})$ et par hypothèse $\gamma (j(s_0)) = j(s_0)$, donc 
$\gamma (j(x)) = j(x)$.
\item[$\ast$] Supposons désormais que $x$ est un point quelconque de l'appartement $A_{\mathbb{F}}$.
Soient $x_0, x_1, \cdots, x_r$ des sommets de $A_{\mathbb{F}}$ tels que 
$\tau = \{ x_0, x_1, \cdots, x_r \}$ est l'unique simplexe de l'appartement $A_{\mathbb{F}}$ 
contenant $x$ dans son intérieur. Il existe alors des réels strictement positifs 
$\lambda_0, \cdots, \lambda_r$ tels que $\sum_{i=0}^r \lambda_i = 1$ et 
$x = \sum_{i=0}^r \lambda_i x_i$.\\
Soit $g$ dans ${\rm Stab}_{{\rm GL}_m (\mathcal{D})} (x)$. Nous allons vérifier que 
$g \in {\rm Stab}_{{\rm GL}_{\mu} (\Delta)} (j(x))$. 
Remarquons tout d'abord que 
$g. \tau = \{ g.x_0, g.x_1, \cdots, g.x_r \}$ est un simplexe de $A_{\mathbb{F}}$ contenant 
$x = g.x = \sum_{i=0}^r \lambda_i g.x_i$ dans son intérieur (car tous les $\lambda_i$ sont non nuls). 
Par unicité, $g. \tau = \tau$ et $g$ appartient à 
${\rm Stab}_{{\rm GL}_m (\mathcal{D})} (\tau)$. D'après \cite{BushnellFröhlich} (partie 1.2 et 1.3), 
on a ${\rm Stab}_{{\rm GL}_m (\mathcal{D})} (\tau) = \langle \Pi_{\tau} \rangle \mathcal{A} (\tau)^{\times}$ 
où $\mathcal{A} (\tau)^{\times} = \cap_{i=0}^r {\rm Stab}_{{\rm GL}_m (\mathcal{D})} (x_i)$ et 
$\Pi_{\tau}$ permute les sommets $x_i$. De plus, on peut choisir $\Pi_{\tau}$ dans $N (T)$ 
(\cite{BushnellFröhlich} loc.cit.). 
On en déduit qu'il existe $l$ dans $\mathbb{Z}$ et $h$ dans $\mathcal{A} (\tau)^{\times}$ tel que 
$g = \Pi_{\tau}^l h$. Ainsi:
$$
g.x = x = \Pi_{\tau}^l h.x 
= \Pi_{\tau}^l . (\sum_{i=0}^r \lambda_i h.x_i) = \Pi_{\tau}^l.x
$$
Par suite, $\Pi_{\tau}^l.x = x$. D'après ce qui précède, $h \in \cap_{i=0}^r {\rm Stab}_{{\rm GL}_m (\mathcal{D})} (x_i) 
\subseteq \cap_{i=0}^r {\rm Stab}_{{\rm GL}_{\mu} (\Delta)} (j (x_i))$. 
En utilisant la $N(T)$-équivariance de $j$ et le fait que 
$\Pi_{\tau}^l \in N(T)$, on a:
$$
g. j(x) = \Pi_{\tau}^l h.j (x) 
= \Pi_{\tau}^l . (\sum_{i=0}^r \lambda_i h.j (x_i)) = \Pi_{\tau}^l.j (x)  
= j (\Pi_{\tau}^l.x) = j (x)
$$
Enfin, étant donné que $j(x) = \sum_{i=0}^r \lambda_i j(x_i)$ et que ${\rm Gal} (\mathbb{K}/ \mathbb{F})$ 
agit de façon affine en fixant l'image de chaque sommet de $A_{\mathbb{F}}$, il est immédiat que pour tout 
$\gamma$ dans ${\rm Gal} (\mathbb{K}/ \mathbb{F})$, $\gamma (j(x)) = j(x)$. 
\end{itemize}

\item[$\ast$] 
Soit $x \in X_{\mathbb{F}}$.
Soient $x_0, y_0 \in A_{\mathbb{F}}$ et $h, g \in {\rm GL}_m (\mathcal{D})$ tels que 
$x = g. x_0 = h.y_0$. 
On a $y_0 = h^{-1} g.x_0$. Comme $N(T)$ agit transitivement sur les sommets de $A_{\mathbb{F}}$, 
il existe $n$ dans $N(T)$ tel que $y_0 = n.x_0$.
On en déduit que $h^{-1} g.x_0 = n.x_0$ donc 
$n^{-1} h^{-1} g$ appartient à ${\rm Stab}_{{\rm GL}_m (\mathcal{D})} (x_0) 
\subseteq {\rm Stab}_{{\rm GL}_{\mu} (\Delta)} (j(x_0))$. 
Ainsi, par $N(T)$-équivariance de $j$ sur $A_{\mathbb{F}}$, il est clair que 
$g.j(x_0) = h.j(y_0)$. 
Par suite, $j$ est bien définie.

\item[$\ast$] On vérifie facilement que ${\rm Im} (j) \subseteq X_{\mathbb{K}}^{{\rm Gal} (\mathbb{K} / \mathbb{F})}$. 
En effet, si $x \in X_{\mathbb{F}}$, $x = g.x_0$ avec $g \in {\rm GL}_m (\mathcal{D})$ (fixé par 
${\rm Gal} (\mathbb{K} / \mathbb{F})$) et 
$x_0 \in A_{\mathbb{F}}$ (lui aussi fixé par ${\rm Gal} (\mathbb{K} / \mathbb{F})$). 

\item[$\ast$] En utilisant le fait que 
l'action de ${\rm GL}_m (\mathcal{D})$ est affine, que $j$ est affine sur $A_{\mathbb{F}}$ 
et ${\rm GL}_m (\mathcal{D})$-équivariante, on montre facilement que l'image par $j$ du barycentre de deux points 
de $X_{\mathbb{F}}$ est le barycentre de l'image de ces deux points, et donc que $j$ est affine.

\end{itemize}

\end{démo}

\noindent
\textbf{Dans toute la suite du chapitre, on fixe $\pi \in \mathcal{R}_0^2 ({\rm GL}_{\mu} (\Delta))$, 
une représentation cuspidale de niveau $0$,  
image d'une cuspidale de niveau $0$ de ${\rm GL}_n (\mathbb{K})$ par la correspondance de 
Jacquet-Langlands, et on garde les notations de \ref{NotationsRepresentationCuspidaleParametreSilbZink}.}

\section{Conditions de distinction lorsque $d$ est impair.}\label{PartieCNSDistinctionImpair}

\noindent
Dans cette partie, on suppose que $d$, l'indice de l'algèbre à division $\mathcal{D}$, est impair.


\begin{rmq}
Puisque $d$ est impair, on remarque que $d=\delta$, $m = \mu$ et 
$\mathcal{D} \otimes_{\mathbb{F}} \mathbb{K} \simeq \Delta$. 
On pourra donc supposer que:
$$
\mathbb{F} \subseteq \mathbb{K} \subseteq \mathcal{D} \subseteq \Delta
$$
\end{rmq}

\begin{nota}
On fixe $\varpi_{\mathbb{F}}$, $\varpi_{\mathbb{K}}$, $\varpi_{\mathcal{D}}$ et $\varpi_{\Delta}$ 
des uniformisantes de $\mathbb{F}$, $\mathbb{K}$, $\mathcal{D}$ et $\Delta$ respectivement telles que 
$\varpi_{\mathcal{D}}^d = \varpi_{\mathbb{F}}$ et $\varpi_{\Delta}^d = \varpi_{\mathbb{K}}$. 
\end{nota}

\subsection{Cas où l'extension $\mathbb{K} / \mathbb{F}$ est non ramifiée.}

On suppose dans toute cette partie que l'extension $\mathbb{K} / \mathbb{F}$ est non ramifiée. 
On peut donc identifier les uniformisantes de $\mathbb{F}$ et de $\mathbb{K}$. 
On en déduit directement la propriété suivante :

\begin{propr}
Pour tout $x$ dans $\mathbb{F}^{\times}$, $v_{\mathbb{K}} (x) = v_{\mathbb{F}} (x)$ et  
pour tout $y$ dans $\mathcal{D}^{\times}$, $v_{\Delta} (y) = v_{\mathcal{D}} (y)$.
\end{propr}

\begin{propr}\label{StabiliteSousGroupeOuvertCompact}
L'extension de corps fini $k_{\Delta} / k_{\mathcal{D}}$ est une extension quadratique.
De plus, on a $\mathcal{O}_{\mathcal{D}} \subseteq \mathcal{O}_{\Delta}$ et le sous-groupe ouvert compact maximal 
${\rm GL}_m (\mathcal{O}_{\Delta})$ est stable sous l'action du groupe de Galois  
$\langle \sigma \rangle = {\rm Gal} (\mathbb{K}/\mathbb{F})$.
\end{propr}

\begin{démo}
Il est clair que $k_{\Delta} / k_{\mathcal{D}}$ est une extension quadratique.\\
Pour montrer les autres propriétés, commençons par quelques rappels 
sur la structure de l'algèbre à division $\mathcal{D}$.
On peut fixer $\mathbb{L} / \mathbb{F}$ une extension non ramifiée de degré $d$ contenue dans $\mathcal{D}$ 
qui soit normalisée par $\varpi_{\mathcal{D}}$, alors $k_{\mathcal{D}} \simeq k_{\mathbb{L}}$. 
Soit $\varphi : \mathbb{L} \rightarrow \mathbb{L}, l \mapsto \varpi_{\mathcal{D}} l \varpi_{\mathcal{D}}^{-1}$.
Alors $\varphi$ induit un générateur $\overline{\varphi}$ du groupe de Galois 
${\rm Gal} (k_{\mathcal{D}} / k_{\mathbb{F}})$.
Posons $\Lambda = \mathbb{L} \otimes_{\mathbb{F}} \mathbb{K}$,
alors $\Lambda$ est un corps car $\mathbb{L}$ et $\mathbb{K}$ sont linéairement disjoints sur $\mathbb{F}$ 
(i.e n'ont pas de $\mathbb{F}$-sous-corps (distincts de $\mathbb{F}$) isomorphes puisque 
$d$ et $2$ sont premiers entre eux).
De plus $\Lambda$ est une extension non ramifiée de degré $d$ de $\mathbb{K}$. 
On pose:
$$
\widetilde{\varphi} : \Lambda \rightarrow \Lambda, \, l \otimes k \mapsto \varphi (l) \otimes k
$$
On vérifie alors que 
$\langle \widetilde{\varphi} \rangle = {\rm Gal} (\Lambda / \mathbb{K})$ 
et que $\widetilde{\varphi}$ est induit par la conjugaison par $\varpi_{\Delta}$.
On déduit de tout ceci que $\{ 1, \varpi_{\Delta}, \cdots, \varpi_{\Delta}^{d-1} \}$ est une $\Lambda$-base de $\Delta$ 
qui engendre $\mathcal{O}_{\Delta}$ comme $\mathcal{O}_{\Lambda}$-module à gauche.
De même, on sait que $\{ 1, \varpi_{\mathcal{D}}, \cdots, \varpi_{\mathcal{D}}^{d-1} \}$ 
est une $\mathbb{L}$-base de $\mathcal{D}$ 
qui engendre $\mathcal{O}_{\mathcal{D}}$ comme $\mathcal{O}_{\mathbb{L}}$-module à gauche.
On vérifie facilement que $\mathcal{O}_{\mathbb{L}} \otimes_{\mathcal{O}_{\mathbb{F}}} \mathcal{O}_{\mathbb{K}}$ 
est un sous-anneau de $\Lambda$ et aussi un réseau (car $\mathcal{O}_{\mathbb{L}}$ et $\mathcal{O}_{\mathbb{K}}$ 
sont des $\mathcal{O}_{\mathbb{F}}$-réseaux). Il s'agit donc d'un sous-anneau compact de $\Lambda$.
Or $\Lambda$ possède un unique sous-anneau compact maximal, $\mathcal{O}_{\Lambda}$, d'où:
$$
\mathcal{O}_{\mathbb{L}} \otimes_{\mathcal{O}_{\mathbb{F}}} \mathcal{O}_{\mathbb{K}} 
\subseteq \mathcal{O}_{\Lambda} 
\Rightarrow \mathcal{O}_{\mathbb{L}} \otimes 1 
\subseteq \mathcal{O}_{\Lambda} \subseteq \mathcal{O}_{\Delta}
$$
Puisque $\mathcal{O}_{\mathcal{D}}$ est engendré par $\mathcal{O}_{\mathbb{L}}$ et $\varpi_{\mathcal{D}}$, 
on a bien l'inclusion 
$\mathcal{O}_{\mathcal{D}} \subseteq \mathcal{O}_{\Delta}$.\\ 
De plus, le groupe de Galois $\langle \sigma \rangle$ agit naturellement sur 
$\Delta = \mathcal{D} \otimes_{\mathbb{F}} \mathbb{K}$, 
et donc sur $\Lambda = \mathbb{L} \otimes_{\mathbb{F}} \mathbb{K}$ par:
$$ 
\sigma . (l \otimes k) = l \otimes \sigma (k)
$$
pour $l$ dans $\mathbb{L}$ et $k$ dans $\mathbb{K}$. 
Ainsi $\sigma (\Lambda) = \Lambda$  et on peut considérer $\widetilde{\sigma} = \sigma_{\vert \Lambda}$.
On vérifie que $\Lambda / \mathbb{L}$ est une extension galoisienne 
non ramifiée de degré $2$ et $\widetilde{\sigma}$ est un générateur de ${\rm Gal} (\Lambda / \mathbb{L})$ 
(car non trivial sur $\Lambda$).
On a:
$$
\widetilde{\sigma} (\varpi_{\Lambda}) = \widetilde{\sigma} (\varpi_{\mathbb{L}}) = \varpi_{\mathbb{L}} = \varpi_{\Lambda}
$$
On en déduit que $\widetilde{\sigma} (\mathcal{O}_{\Lambda}) = \mathcal{O}_{\Lambda}$.
De plus $\sigma. \varpi_{\Delta}^{j} = \varpi_{\Delta}^{j}$ pour tout $j \in \{ 0, \cdots, d-1 \}$. 
Comme $\{ 1, \varpi_{\Delta}, \cdots, \varpi_{\Delta}^{d-1} \}$ est une $\Lambda$-base de $\Delta$ 
qui engendre $\mathcal{O}_{\Delta}$ comme $\mathcal{O}_{\Lambda}$-module à gauche:
$$
\mathcal{O}_{\Delta} 
= \mathcal{O}_{\Lambda}.1 \oplus \mathcal{O}_{\Lambda}.\varpi_{\Delta} \oplus 
\cdots \oplus \mathcal{O}_{\Lambda}.\varpi_{\Delta}^{d-1}
$$
On a 
$\sigma (\mathcal{O}_{\Delta}) = \mathcal{O}_{\Delta}$ et pour tout $k$ dans $\mathbb{Z}$,  
$\sigma (\mathcal{P}_{\Delta}^{k}) = \mathcal{P}_{\Delta}^{k}$. 
Par conséquent $\sigma$ stabilise le compact maximal ${\rm GL}_m (\mathcal{O}_{\Delta})$.

\end{démo}

\subsubsection{Explicitation des injections d'immeubles.}

\begin{propr}\label{InjectionImmeubleNRImpair}
Soit $A_{\mathbb{F}}$ (resp. $A_{\mathbb{K}}$) l'appartement standard de $X_{\mathbb{F}}$ 
(resp. $X_{\mathbb{K}}$) que l'on identifie à l'espace affine 
$\mathbb{R}^m / \mathbb{R} (1, \cdots, 1)$ et dont l'ensemble des sommets s'identifie à 
$\mathbb{Z}^m / \mathbb{Z} (1, \cdots, 1)$. 
L'injection naturelle entre les deux immeubles est donnée par :
$$
j : A_{\mathbb{F}} \rightarrow A_{\mathbb{K}}, \, 
\overline{(x_1, \cdots, x_m)} \mapsto \overline{(x_1, \cdots, x_m)}
$$
\end{propr}

\begin{démo}
On définit $j$, injection de $A_{\mathbb{F}}$ sur $A_{\mathbb{K}}$ par :
$$
j(\overline{(x_1, \cdots, x_m)}) = \overline{(x_1, \cdots, x_m)}
$$
pour tout $\overline{(x_1, \cdots, x_m)}$ dans $A_{\mathbb{F}}$. 
Montrons que $j$ est l'injection naturelle des immeubles.
On sait que l'on peut supposer que $\mathcal{D} \subseteq \Delta$ 
et que l'on a, pour tout $x \in \mathcal{D}^{\times}$, $v_{\Delta} (x) = v_{\mathcal{D}} (x)$.
Soit $T$ le tore maximal déployé diagonal de ${\rm GL}_m (\mathcal{D})$ (tore associé à l'appartement standard):
$$
T = \left\{ \left( \begin{array}{cccc}
t_1 & 0 & \cdots & 0\\
0 & \ddots &  & \vdots\\
\vdots & & \ddots & 0\\
0 & \cdots & 0 & t_m\\
            \end{array} \right) : t_i \in \mathbb{F}^{\times} \right\}
$$
Alors $N(T) = T_0 \mathcal{S}_m$
où:
$$
T_0 = \left\{ \left( \begin{array}{cccc}
t_1 & 0 & \cdots & 0\\
0 & \ddots &  & \vdots\\
\vdots & & \ddots & 0\\
0 & \cdots & 0 & t_m\\
            \end{array} \right) : t_i \in \mathcal{D}^{\times} \right\}
$$
et $\mathcal{S}_m$ désigne l'ensemble des matrices de permutation : si $\tau$ est une permutation de 
$\{ 1, \cdots, m \}$, on lui associe la matrice $P_{\tau} = [p_{i,j}]$ telle que $p_{i,j} = 1$ 
si et seulement si $i = \tau (j)$ et $p_{i,j} = 0$ sinon.\\
Soit $g \in N(T)$ tel que $g = t P_{\tau}$ où $t$ est une matrice diagonale 
${\rm diag} (t_1, \cdots, t_m)$ et $\tau$ une permutation de $\lbrace 1, \cdots, m \rbrace$.
Alors, pour tout $\overline{(x_1, \cdots, x_m)}$ dans $A_{\mathbb{F}}$:
$$
g. \overline{(x_1, \cdots, x_m)} = 
\overline{(v_{\mathcal{D}} (t_1) + x_{\tau (1)}, \cdots, v_{\mathcal{D}} (t_m) + x_{\tau (m)})}
$$
et si $g$ est vu comme élément de ${\rm GL}_m (\Delta)$, on a, pour tout 
$\overline{(y_1, \cdots, y_m)}$ dans $A_{\mathbb{K}}$:
$$
g. \overline{(y_1, \cdots, y_m)} = 
\overline{(v_{\Delta} (t_1) + y_{\tau(1)}, \cdots, v_{\Delta} (t_m) + y_{\tau (m)})}
$$
On en déduit facilement que $j$ est bien $N(T)$-équivariante sur l'appartement $A_{\mathbb{F}}$.\\
De plus, si l'on fixe 
$
s_0 = [ \mathcal{O}_{\mathcal{D}} \oplus \cdots \oplus \mathcal{O}_{\mathcal{D}} ] 
= \overline{(0, \cdots, 0)}
$
un sommet de $A_{\mathbb{F}}$. Son image par $j$ est 
$
j(s_0) = \overline{(0, \cdots, 0)} 
= [ \mathcal{O}_{\Delta} \oplus \cdots \oplus \mathcal{O}_{\Delta} ]
$.\\
 Ainsi ${\rm Stab}_{{\rm GL}_m (\mathcal{D})} (s_0 ) 
= \langle \varpi_{\mathcal{D}} \rangle {\rm GL}_m (\mathcal{O}_{\mathcal{D}})$
et ${\rm Stab}_{{\rm GL}_m (\Delta)} (j(s_0) ) 
= \langle \varpi_{\Delta} \rangle {\rm GL}_m (\mathcal{O}_{\Delta})$.
On a donc clairement 
$
{\rm Stab}_{{\rm GL}_m (\mathcal{D})} (s_0 ) \subseteq 
{\rm Stab}_{{\rm GL}_m (\Delta)} (j(s_0) )
$. 
De plus :
$$
\sigma ({\rm GL}_m (\mathcal{O}_{\Delta})) = {\rm GL}_m (\mathcal{O}_{\Delta})
$$
donc $\sigma ({\rm Stab}_{{\rm GL}_m (\Delta)} (j(s_0) )) = {\rm Stab}_{{\rm GL}_m (\Delta)} (j(s_0) )$.
On en déduit que le sommet $j(s_0)$ est fixé par l'action du groupe de Galois  
${\rm Gal} (\mathbb{K} / \mathbb{F})$. 
On conclut en utilisant la proposition \ref{PropositionCaracterisationInjection}.

\end{démo}

\begin{propr}\label{DefinitionSommet}
Notons $C_0$ la chambre standard de $X_{\mathbb{F}}$ et 
$\widetilde{C}_0$ la chambre standard de $X_{\mathbb{K}}$.
Alors $j(C_0) = \widetilde{C}_0$.
\end{propr}

\begin{démo}
Notons $\{ s_0, s_1, \cdots, s_{m-1} \}$ les sommets de $C_0$, de sorte que pour tout $i$ 
dans $\{ 0, \cdots, m-1 \}$, on a:
$$
s_i = [\mathcal{O}_{\mathcal{D}} \oplus \cdots \oplus \mathcal{O}_{\mathcal{D}}
 \oplus \underbrace{\mathcal{P}_{\mathcal{D}} \oplus \cdots \oplus \mathcal{P}_{\mathcal{D}}}_i]
$$
De même, on note $\{ S_0, S_1, \cdots, S_{m-1} \}$ les sommets de $\widetilde{C}_0$, tels que pour tout $i$ 
dans $\{ 0, \cdots, m-1 \}$:
$$
S_i = [\mathcal{O}_{\Delta} \oplus \cdots \oplus \mathcal{O}_{\Delta}
 \oplus \underbrace{\mathcal{P}_{\Delta} \oplus \cdots \oplus \mathcal{P}_{\Delta}}_i]
$$
L'image de $C_0$, notée $j(C_0)$, est l'image par $j$ de l'enveloppe convexe des points 
$\{ s_0, s_1, \cdots, s_{m-1} \}$. Comme $j$ est affine, il s'agit de 
l'enveloppe convexe des points $\{ j(s_0), j(s_1), \cdots, j(s_{m-1}) \}$.
Pour montrer que $j(C_0) = \widetilde{C}_0$, il suffit donc de vérifier que:
$$
\{ j(s_0), j(s_1), \cdots, j(s_{m-1}) \} 
= \{ S_0, S_1, \cdots, S_{m-1} \}
$$
Pour tout $i$ dans $\{ 0, \cdots, m-1 \}$, on a 
$
s_i = h_i.s_0 
$
où $h_i = {\rm diag} (1, \cdots, 1, \underbrace{\varpi_{\mathcal{D}}, \cdots, \varpi_{\mathcal{D}}}_i)$ 
est une matrice diagonale de ${\rm GL}_m (\mathcal{D})$.  
Par ${\rm GL}_m (\mathcal{D})$-équivariance de $j$, on a:
$$
j(s_i) = h_i.j(s_0) = S_i
$$

\end{démo}

\subsubsection{Conditions de distinction.}

On utilisera dans toute la suite le résultat suivant qui est en fait un cas particulier 
de la proposition 5.20 page 94 
de \cite{HakimMurnaghan2} :

\begin{theo}\label{ResultatHakimMurnaghan}
(\cite{HakimMurnaghan2})
Par formule de restriction de Mackey et Réciprocité de Frobenius pour l'induite compacte, 
on obtient un isomorphisme de $\mathbb{C}$-espaces vectoriels :
$$
{\rm Hom}_{{\rm GL}_m (\mathcal{D})} ({\rm c-Ind}_{\mathcal{K}}^{{\rm GL}_{\mu} (\Delta)} \pi_0, \mathds{1}) 
\simeq \bigoplus_{g \in {\rm GL}_m (\mathcal{D}) \backslash {\rm GL}_{\mu} (\Delta) / \mathcal{K}} 
{\rm Hom}_{g^{-1} {\rm GL}_m (\mathcal{D}) g \cap \mathcal{K}} (\pi_0, \mathds{1}), 
$$
où $\mathcal{K} = \langle \varpi_{\mathbb{K}} \rangle {\rm GL}_{\mu} (\mathcal{O}_{\Delta})$. 
Si $g \in {\rm GL}_m (\mathcal{D}) \backslash {\rm GL}_{\mu} (\Delta) / \mathcal{K}$ vérifie:
$$
{\rm Hom}_{g^{-1} {\rm GL}_m (\mathcal{D}) g \cap \mathcal{K}} (\pi_0, \mathds{1}) \neq 0,
$$ 
alors nécessairement 
$\sigma (g \mathcal{K} g^{-1}) = g \mathcal{K} g^{-1}$.
\end{theo}

\begin{propo}\label{RaisonnementDistinction}
On a:
$$
{\rm Hom}_{{\rm GL}_m (\mathcal{D})} (\pi, \mathds{1}) 
\simeq {\rm Hom}_{{\rm GL}_m (\mathcal{D}) \cap \mathcal{K}} (\pi_0, \mathds{1}) 
\, \, \text{où} \, \, 
\mathcal{K} = \langle \varpi_{\mathbb{K}} \rangle {\rm GL}_m (\mathcal{O}_{\Delta})
$$
Par conséquent, la représentation $\pi$ est ${\rm GL}_m (\mathcal{D})$-distinguée si et seulement si:
$$
{\rm Hom}_{{\rm GL}_m (\mathcal{D}) \cap \mathcal{K}} (\pi_0, \mathds{1}) \neq 0 
$$
\end{propo}

\begin{démo}
Notons $K = {\rm GL}_m (\mathcal{O}_{\Delta})$.
Alors $j(s_0)$ peut être vu comme une chaîne de réseaux (on fera souvent cette identification par la suite) 
et $K = \mathcal{A} (j(s_0))^{\times}$ est le sous-groupe parahorique fixateur de $j(s_0)$. 
Le sous-groupe $K$ est stable par l'action
 du groupe de Galois $\langle \sigma \rangle = {\rm Gal} (\mathbb{K} / \mathbb{F})$ sur 
les coefficients.
Ceci signifie que le sommet $j(s_0)$ est fixé par l'action de $\langle \sigma \rangle$ 
(rappelons que le stabilisateur d'un point de l'immeuble caractérise ce point, on pourra se 
référer à \cite{BushnellFröhlich} page 215). 
D'après \ref{ResultatHakimMurnaghan}, les seuls éléments $g$ dans 
${\rm GL}_m (\mathcal{D}) \backslash {\rm GL}_m (\Delta) / \mathcal{K}$ 
qui contribuent à la distinction sont ceux qui vérifient les conditions suivantes:
$$
\sigma (g K g^{-1}) = g K g^{-1} \, \, \text{et} \, \, 
{\rm Hom}_{g^{-1} {\rm GL}_m (\mathcal{D}) g \cap \mathcal{K}} (\pi_0, \mathds{1}) \neq 0
$$
On remarque que $g K g^{-1} = \mathcal{A} (g.j(s_0))^{\times}$,
ainsi $\sigma (g K g^{-1}) = g K g^{-1}$ si et seulement si 
le sommet $g. j(s_0)$ est stable sous l'action de $\langle \sigma \rangle$. 
D'après l'article \cite{PrasadYu}, théorème 1.9 page 555,  
les sommets de $X_{\mathbb{K}}$ fixes sous l'action de $\langle \sigma \rangle$ sont exactement les sommets de 
$X_{\mathbb{K}}$ qui sont dans $j(X_{\mathbb{F}})$.
On déduit de tout ceci que $\pi$ est distinguée si et seulement s'il existe $g$ dans 
${\rm GL}_m (\mathcal{D}) \backslash {\rm GL}_m (\Delta) / \mathcal{K}$ 
tel que $g.j(s_0)$ appartient à $j(X_{\mathbb{F}})$ et:
$$
{\rm Hom}_{g^{-1} {\rm GL}_m (\mathcal{D}) g \cap \mathcal{K}} (\pi_0, \mathds{1}) \neq 0
$$
Chaque ${\rm GL}_m (\mathcal{D})$-orbite d'un point de $X_{\mathbb{F}}$ a un représentant dans la chambre standard $C_0$.
Donc, puisque l'on cherche $g$ modulo ${\rm GL}_m (\mathcal{D})$, on se ramène à la recherche des éléments 
$g$ dans ${\rm GL}_m (\Delta)$ tels que $g.j(s_0)$ est un sommet de $j(C_0)$.
Or, on a vu que les sommets de $X_{\mathbb{K}}$ qui sont dans l'image de $C_0$ sont exactement:
$$
\lbrace j(s_0), \cdots, j(s_{m-1}) \rbrace
$$
Par ${\rm GL}_m (\mathcal{D})$-équivariance de $j$, comme les sommets de $C_0$ sont dans la 
même ${\rm GL}_m (\mathcal{D})$-orbite, quitte à multiplier $g$ à gauche par un élément de 
${\rm GL}_m (\mathcal{D})$, on peut se ramener au sommet $j(s_0)$.
De plus, on a vu précédemment que le sommet $j(s_0)$ est stable sous l'action de $\sigma$ car l'ordre héréditaire 
correspondant, $\mathcal{A} (j(s_0))^{\times} = {\rm GL}_m (\mathcal{O}_{\Delta})$ est stable sous l'action de $\sigma$.
On a donc bien le résultat annoncé.

\end{démo}

\begin{propo}\label{PropositionIntersection}
La représentation $\pi$ est ${\rm GL}_m (\mathcal{D})$-distinguée si et seulement si 
$\chi$ est trivial sur $\mathbb{F}^{\times}$ et:
$$
{\rm Hom}_{{\rm GL}_m (k_{\mathcal{D}})} (\overline{\gamma}_0, \mathds{1}) \neq 0
$$
\end{propo}

\begin{démo}
D'après ce qui précède, on sait que $\pi$ est ${\rm GL}_m (\mathcal{D})$-distinguée si et seulement si:
$$
{\rm Hom}_{{\rm GL}_m (\mathcal{D}) \cap \mathcal{K}} (\pi_0, \mathds{1}) \neq 0
$$
Il nous suffit donc de déterminer l'intersection 
${\rm GL}_m (\mathcal{D}) \cap 
\langle \varpi_{\mathbb{K}} \rangle {\rm GL}_m (\mathcal{O}_{\Delta})$.\\
On fixe $\Lambda / \mathbb{K}$ extension non ramifiée de degré $d$ contenue dans $\Delta$.
Il existe alors une injection $\Delta \subseteq {\rm M}_{d} (\Lambda)$ telle que, 
via cette injection, $\varpi_{\Delta}$ s'identifie à la matrice:
$$
w_0 = \left( \begin{array}{ccccc}
0 & 1 & 0 & \cdots & 0\\
\vdots & 0& 1 & & \vdots\\
\vdots & \vdots &\ddots & \ddots & 0\\
0 & \vdots &  & \ddots & 1\\
\varpi_{\mathbb{K}} & 0 & \cdots & \cdots & 0\\
\end{array} \right)
$$
(on pourra se référer à la démonstration du théorème 14.6 du chapitre 3 dans \cite{Reiner}). 
Nous allons travailler dans la $\mathbb{K}$-algèbre 
$A_0 = {\rm M}_{m} (\Delta) \otimes_{\mathbb{K}} \Lambda$ pour calculer la norme réduite 
des éléments de ${\rm M}_{m} (\Delta)$.
L'élément $\varpi_{\Delta} I_{m} \otimes 1 \in A_0$ s'identifie à la matrice diagonale par blocs:
$$
W_0 = {\rm diag} (\underbrace{w_0, \cdots, w_0}_{m})
$$
Il est clair que 
${\rm det} (w_0) = (-1)^{d+1} \varpi_{\mathbb{K}} = \varpi_{\mathbb{K}}$ (développement par rapport à la première colonne).
On en déduit que 
${\rm Nrd} (\varpi_{\Delta} I_{m}) = {\rm det} (W_0) = \varpi_{\mathbb{K}}^{m}$
et:
$$
{\rm Nrd} (\varpi_{\mathbb{K}} I_{m}) = {\rm Nrd} ((\varpi_{\Delta} I_{m})^{d}) 
= \varpi_{\mathbb{K}}^{dm} = \varpi_{\mathbb{K}}^{n}
$$
Ainsi, $v_{\mathbb{K}} ({\rm Nrd} (\varpi_{\mathbb{K}} I_{m})) = n$,  
$v_{\mathbb{K}} ({\rm Nrd} (\varpi_{\Delta} I_{m})) = m$ et, puisque  
${\rm Nrd} ({\rm GL}_{m} (\mathcal{O}_{\Delta})) \subseteq \mathcal{O}_{\mathbb{K}}^{\times}$, 
on a:
$$
 \langle \varpi_{\mathbb{K}} \rangle {\rm GL}_{m} (\mathcal{O}_{\Delta}) 
 = \lbrace g \in  \langle \varpi_{\Delta} \rangle {\rm GL}_{m} (\mathcal{O}_{\Delta}) : 
n \, \, \text{divise} \, \, v_{\mathbb{K}} ({\rm Nrd} (g)) \rbrace 
$$
On montre de même que:
$$
 \langle \varpi_{\mathbb{F}} \rangle {\rm GL}_{m} (\mathcal{O}_{\mathcal{D}}) 
 = \lbrace g \in  \langle \varpi_{\mathcal{D}} \rangle {\rm GL}_m (\mathcal{O}_{\mathcal{D}}) : 
n \, \, \text{divise} \, \, v_{\mathbb{F}} ({\rm Nrd} (g)) \rbrace 
$$
On a:
$$
{\rm GL}_m (\mathcal{D}) \cap \langle \varpi_{\Delta} \rangle 
{\rm GL}_m (\mathcal{O}_{\Delta}) 
= {\rm GL}_m (\mathcal{D}) \cap {\rm Stab}_{{\rm GL}_m (\Delta)} (j(s_0)) 
= {\rm Stab}_{{\rm GL}_m (\mathcal{D})} (s_0)
$$
Ainsi ${\rm GL}_m (\mathcal{D}) \cap \langle \varpi_{\Delta} \rangle 
{\rm GL}_m (\mathcal{O}_{\Delta}) 
= \langle \varpi_{\mathcal{D}} \rangle 
{\rm GL}_m (\mathcal{O}_{\mathcal{D}})$. 
On en déduit que $g \in {\rm GL}_m (\mathcal{D}) \cap \langle \varpi_{\mathbb{K}} \rangle 
{\rm GL}_m (\mathcal{O}_{\Delta})$ si et seulement si $g \in \langle \varpi_{\mathcal{D}} \rangle 
{\rm GL}_m (\mathcal{O}_{\mathcal{D}})$ et $v_{\mathbb{F}} ({\rm Nrd} (g)) = v_{\mathbb{K}} ({\rm Nrd} (g))$ 
est divisible par $n$. 
Finalement ${\rm GL}_m (\mathcal{D}) \cap \langle \varpi_{\mathbb{K}} \rangle 
{\rm GL}_m (\mathcal{O}_{\Delta}) 
= \langle \varpi_{\mathbb{F}} \rangle 
{\rm GL}_m (\mathcal{O}_{\mathcal{D}})$.
Par conséquent, $\pi$ est ${\rm GL}_m (\mathcal{D})$-distinguée si et seulement si:
$$
{\rm Hom}_{\langle \varpi_{\mathbb{F}} \rangle 
{\rm GL}_m (\mathcal{O}_{\mathcal{D}}) } (\pi_0, \mathds{1}) \neq 0
$$
Supposons que:
$$
{\rm Hom}_{\langle \varpi_{\mathbb{F}} \rangle 
{\rm GL}_m (\mathcal{O}_{\mathcal{D}}) } (\pi_0, \mathds{1}) \neq 0
$$
En écrivant les conditions de distinction sur les éléments du centre, on vérifie que 
$\chi$, le caractère central de $\pi_0$, est trivial sur $\mathbb{F}^{\times}$ et donc:
$$
{\rm Hom}_{ \langle \varpi_{\mathbb{F}} \rangle {\rm GL}_m (\mathcal{O}_{\mathcal{D}}) } (\pi_0, \mathds{1}) 
\simeq  {\rm Hom}_{ {\rm GL}_m (\mathcal{O}_{\mathcal{D}}) } (\pi_0, \mathds{1}) 
\simeq {\rm Hom}_{{\rm GL}_m (k_{\mathcal{D}}) } (\overline{\gamma}_0, \mathds{1}) \neq 0
$$
La réciproque est immédiate.

\end{démo}

\begin{theo}\label{CritereDistinctionNRdImpair}
La représentation $\pi$ est ${\rm GL}_m (\mathcal{D})$-distinguée si et seulement si 
$n$ est impair, $\chi$ est trivial sur $\mathbb{F}^{\times}$ et:
$$
\overline{\chi}^{-1} \simeq_{k_{\Delta}} \overline{\chi} \circ \widetilde{\widetilde{\sigma}} 
\, \, \text{où} \, \, \langle \widetilde{\widetilde{\sigma}} \rangle = {\rm Gal} (k_{\Delta,m} / k_{\mathcal{D}})
$$
De plus, si $\pi$ est ${\rm GL}_m (\mathcal{D})$-distinguée, on a:
$$
{\rm dim}_{\mathbb{C}} ({\rm Hom}_{{\rm GL}_m (\mathcal{D})} (\pi, \mathds{1})) = 1
$$
\end{theo}

Pour montrer ce résultat, nous utilisons la paramétrisation de Green pour les représentations cuspidales du groupe 
linéaire sur un corps fini (cf. \cite{Green}) que nous rappelons ici :

\begin{theo}\label{ParametrisationGreen} (\cite{Green}) 
Soit $k$ un corps fini. Alors, les représentations cuspidales irréductibles de 
${\rm GL}_f (k)$ sont en bijection avec les caractères $k$-réguliers de $l^{\times}$ où $l$ est une extension de degré $f$ de $k$. 
Si $\delta : l^{\times} \rightarrow \mathbb{C}^{\times}$ est un caractère $k$-régulier, la représentation cuspidale 
$\pi_{\delta}$ correspondant au caractère $\delta$ est caractérisée par le fait que pour tout élément elliptique régulier 
$x$ dans ${\rm GL}_f (k)$, on a:
$$
{\rm tr} (\pi_{\delta} (x)) = (-1)^{f-1} 
\sum_{\beta \in {\rm Gal} (l/k)} \delta (x^{\beta})
$$
(où ${\rm tr} (\pi_{\delta} (x))$ est la trace de $\pi_{\delta} (x)$).\\
\end{theo}

Démontrons à présent le théorème \ref{CritereDistinctionNRdImpair} :\\

\newpage
\begin{démo}

\begin{itemize}
\item[$\ast$] Supposons que $\pi$ est ${\rm GL}_m (\mathcal{D})$-distinguée.
Alors, d'après la proposition précédente, on sait que $\chi$ est trivial sur 
$\mathbb{F}^{\times}$ et :
$$
{\rm Hom}_{{\rm GL}_m (k_{\mathcal{D}})} (\overline{\gamma}_0, \mathds{1}) \neq 0
$$
L'extension $k_{\Delta} / k_{\mathcal{D}}$ est une extension quadratique. D'après 
\cite{Prasad} théorème 2 page 341, on~a :
$$
{\rm Hom}_{{\rm GL}_m (k_{\mathcal{D}})} (\overline{\gamma}_0, \mathds{1}) \neq 0 
\Leftrightarrow \overline{\gamma}_0^{\vee} \simeq \overline{\gamma}_0^{\widetilde{\sigma}}
$$
où $\langle \widetilde{\sigma} \rangle = {\rm Gal} (k_{\Delta} / k_{\mathcal{D}})$ 
et $\overline{\gamma}_0^{\vee}$ est la contragrédiente de $\overline{\gamma}_0$.
Les représentations $\overline{\gamma}_0^{\vee}$ et $\overline{\gamma}_0^{\widetilde{\sigma}}$ 
sont équivalentes si et seulement si elles ont même paramètre de Green.
Le paramètre de Green associé à $\overline{\gamma}_0^{\vee}$ est 
$\overline{\chi}^{-1}$. Il nous reste à déterminer le paramètre de Green de la représentation 
$\overline{\gamma}_0^{\widetilde{\sigma}}$ en utilisant \ref{ParametrisationGreen}.\\
Soit $x \in {\rm GL}_{m} (k_{\Delta})$ un élément elliptique régulier. Nous allons expliciter un 
isomorphisme de corps entre $k_{\Delta} [x]$ et $k_{\Delta,m} = k_{\mathbb{K},n}$.
Fixons $v \in k_{\Delta,m}$ ayant même polynôme minimal (sur $k_{\Delta}$) que $x$.
Notons $\mu_{k_{\Delta}}^x$ (resp. $\mu_{k_{\Delta}}^v$) ces polynômes minimaux, alors 
$\mu_{k_{\Delta}}^x = \mu_{k_{\Delta}}^v$ par hypothèse. On fixe 
$\varphi : k_{\Delta} [x] \rightarrow k_{\Delta,m}$ un $k_{\Delta}$-isomorphisme de corps en imposant 
que $\varphi (x) = v$. D'après la formule des caractères, on a:
$$
{\rm tr} (\overline{\gamma}_0 (x)) 
= (-1)^{d-1} \sum_{\beta \in {\rm Gal} (k_{\Delta,m} / k_{\Delta})} \overline{\chi} (\beta \circ \varphi (x)) 
= \sum_{\beta \in {\rm Gal} (k_{\Delta,m} / k_{\Delta})} \overline{\chi} (\beta \circ \varphi (x))
$$
(car $d$ est impair). Fixons $\widetilde{\widetilde{\sigma}}$ un générateur du groupe de Galois 
${\rm Gal} (k_{\Delta,m} / k_{\mathcal{D}})$. La restriction de $\widetilde{\widetilde{\sigma}}$ 
à $k_{\Delta}$ ne peut pas être triviale, donc:
$$
\widetilde{\widetilde{\sigma}}_{\vert k_{\Delta}} = \widetilde{\sigma}
$$
Alors le polynôme $\widetilde{\widetilde{\sigma}} (\mu_{k_{\Delta}}^x)$ est un polynôme irréductible 
annulateur de $\widetilde{\sigma} (x)$, il s'agit du polynôme minimal de 
$\widetilde{\sigma} (x)$. On vérifie de même que $\widetilde{\widetilde{\sigma}} (\mu_{k_{\Delta}}^v)$ 
est le polynôme minimal de $\widetilde{\widetilde{\sigma}} (v)$. 
Par conséquent, $\widetilde{\widetilde{\sigma}} (v)$ et $\widetilde{\sigma} (x)$ ont même polynôme minimal 
(sur $k_{\Delta}$).
On fixe $\psi: k_{\Delta} [\widetilde{\sigma} (x)] \rightarrow k_{\Delta,m}$ un 
$k_{\Delta}$-isomorphisme de corps en imposant:
$$
\psi (\widetilde{\sigma} (x)) = \widetilde{\widetilde{\sigma}} (v)
$$
On a donc $\psi \circ \widetilde{\sigma} (x) = \widetilde{\widetilde{\sigma}} \circ \varphi (x)$.
On en déduit que:
\begin{eqnarray*}
{\rm tr} (\overline{\gamma}_0^{\widetilde{\sigma}} (x)) 
& = & {\rm tr} (\overline{\gamma}_0 ({\widetilde{\sigma}} (x))
     =  \sum_{\beta \in {\rm Gal} (k_{\Delta,m} / k_{\Delta})} \overline{\chi} 
           (\beta \circ \psi ({\widetilde{\sigma}} (x)))\\
& = & \sum_{\beta \in {\rm Gal} (k_{\Delta,m} / k_{\Delta})} \overline{\chi} 
                (\beta \circ \widetilde{\widetilde{\sigma}} (\varphi (x)))
       =  \sum_{\beta \in {\rm Gal} (k_{\Delta,m} / k_{\Delta})} \overline{\chi} 
                (\widetilde{\widetilde{\sigma}} \circ \beta (\varphi (x)))\\
& = & \sum_{\beta \in {\rm Gal} (k_{\Delta,m} / k_{\Delta})} 
       \overline{\chi} \circ \widetilde{\widetilde{\sigma}} (\beta \circ \varphi (x))
\end{eqnarray*}
(car $\widetilde{\widetilde{\sigma}} \in {\rm Gal} (k_{\Delta,m} / k_{\mathcal{D}}) 
\supseteq {\rm Gal} (k_{\Delta,m} / k_{\Delta})$ donc commute avec tous les éléments de 
${\rm Gal} (k_{\Delta,m} / k_{\Delta})$). On vérifie facilement que le caractère 
$\overline{\chi} \circ \widetilde{\widetilde{\sigma}}$ est $k_{\Delta}$-régulier. 
En effet, si $\alpha \in {\rm Gal} (k_{\Delta,m} / k_{\Delta})$ vérifie 
$\overline{\chi} \circ \widetilde{\widetilde{\sigma}} \circ \alpha = 
\overline{\chi} \circ \widetilde{\widetilde{\sigma}}$, comme $\alpha$ et 
$\widetilde{\widetilde{\sigma}}$ commutent, la $k_{\Delta}$-régularité de $\overline{\chi}$ 
impose que $\alpha = Id$.\\
Par conséquent, le paramètre de Green associé à la représentation $\overline{\gamma}_0^{\widetilde{\sigma}}$ 
est $\overline{\chi} \circ \widetilde{\widetilde{\sigma}}$.
Par suite, si $\pi$ est ${\rm GL}_m (\mathcal{D})$-distinguée, on a:
$$
\overline{\chi}^{-1} \simeq_{k_{\Delta}} \overline{\chi} \circ \widetilde{\widetilde{\sigma}}
$$
Il existe donc 
$\alpha$ dans ${\rm Gal} (k_{\Delta,m} / k_{\Delta})$ tel que 
$\overline{\chi}^{-1} \circ \alpha = \overline{\chi} \circ \widetilde{\widetilde{\sigma}}$.
Soit $\delta = \widetilde{\widetilde{\sigma}} \circ \alpha^{-1}$, alors 
$\overline{\chi} \circ \delta^2 = \overline{\chi}$.
De plus, $\delta^2 = \widetilde{\widetilde{\sigma}}^2 \circ \alpha^{-2}$ avec 
$\langle \widetilde{\widetilde{\sigma}}^2 \rangle = {\rm Gal} (k_{\Delta,m} / k_{\Delta})$. 
On en déduit que $\delta^2 \in {\rm Gal} (k_{\Delta,m} / k_{\Delta})$.
Comme $\overline{\chi}$ est $k_{\Delta}$-régulier, on a forcément:
$$
\delta^2 = Id \Rightarrow \widetilde{\widetilde{\sigma}}^2 = \alpha^{2}
$$
Puisque $\alpha \in \langle \widetilde{\widetilde{\sigma}}^2 \rangle 
= {\rm Gal} (k_{\Delta,m} / k_{\Delta})$, il existe $i \in \mathbb{Z}$ tel que:
$$
\alpha = \widetilde{\widetilde{\sigma}}^{2i} \Rightarrow 
\widetilde{\widetilde{\sigma}}^2 = \widetilde{\widetilde{\sigma}}^{4i} 
\Rightarrow (\widetilde{\widetilde{\sigma}}^2)^{2i-1} = Id
$$
Ainsi, $\widetilde{\widetilde{\sigma}}^2$ est d'ordre impair.
Comme $\widetilde{\widetilde{\sigma}}^2$ est d'ordre $m$, $m$ est nécessairement impair et 
$n= dm$ est aussi impair (car $d$ est par hypothèse impair).
Donc, si $\pi$ est distinguée, $n$ est impair.
De plus, si $\pi$ est distinguée, on a d'après \ref{RaisonnementDistinction} :
$$
{\rm Hom}_{{\rm GL}_m (\mathcal{D})} (\pi, \mathds{1}) 
\simeq {\rm Hom}_{\langle \varpi_{\mathbb{F}} \rangle {\rm GL}_m (\mathcal{O}_{\mathcal{D}})} (\pi_0, \mathds{1}) 
\simeq {\rm Hom}_{{\rm GL}_m (k_{\mathcal{D}})} (\overline{\gamma}_0, \mathds{1})
$$
Or, d'après l'article de D.Prasad \cite{Prasad}, théorème 3 page 341, l'espace 
${\rm Hom}_{{\rm GL}_m (k_{\mathcal{D}})} (\overline{\gamma}_0, \mathds{1})$ est de dimension au plus $1$, d'où le résultat.\\

\item[$\ast$] Supposons à présent que $n$ est impair, $\chi$ trivial sur $\mathbb{F}^{\times}$ et 
$\overline{\chi}^{-1} \simeq_{k_{\Delta}} \overline{\chi} \circ \widetilde{\widetilde{\sigma}}$. 
Alors les représentations $\overline{\gamma}_0^{\vee}$ et 
$\overline{\gamma}_0^{\widetilde{\sigma}}$ sont équivalentes (car ont même paramètre de Green) et
d'après \cite{Prasad}, théorème 2 page 341, on a:
$$
{\rm Hom}_{{\rm GL}_m (k_{\mathcal{D}})} (\overline{\gamma}_0, \mathds{1}) \neq 0
$$
La proposition précédente nous permet immédiatement de conclure que $\pi$ est bien distinguée dans ce cas.

\end{itemize}
\end{démo}

\subsection{Cas où l'extension $\mathbb{K} / \mathbb{F}$ est totalement ramifiée, modérément ramifiée.}

On suppose dans toute cette partie que $\mathbb{K} / \mathbb{F}$ est totalement ramifiée, modérément ramifiée. 
On fixe les uniformisantes $\varpi_{\mathbb{F}}$ et $\varpi_{\mathbb{K}}$ 
telles que $\varpi_{\mathbb{K}}^2 = \varpi_{\mathbb{F}}$. 
On a directement les propriétés suivantes : 

\begin{propr}
Pour tout $x \in \mathbb{F}^{\times}$, $v_{\mathbb{K}} (x) = 2 v_{\mathbb{F}} (x)$ et 
pour tout $x \in \mathcal{D}^{\times}$, $v_{\Delta} (x) = 2 v_{\mathcal{D}} (x)$.
\end{propr}

En utilisant les mêmes arguments que dans la démonstration de \ref{StabiliteSousGroupeOuvertCompact}, on montre que : 

\begin{propr}
On a $\mathcal{O}_{\mathcal{D}} \subseteq \mathcal{O}_{\Delta}$ et le sous-groupe ouvert compact maximal 
${\rm GL}_m (\mathcal{O}_{\Delta})$ est stable sous l'action du groupe de Galois de 
$\langle \sigma \rangle = {\rm Gal} (\mathbb{K}/\mathbb{F})$.
\end{propr}

\subsubsection{Explicitation des injections d'immeubles.}

\begin{propr}\label{InjectionImmeubleTRImpair}
Soient $A_{\mathbb{F}}$ et $A_{\mathbb{K}}$ les appartements standards de $X_{\mathbb{F}}$ et $X_{\mathbb{K}}$ 
que l'on identifie à l'espace affine $\mathbb{R}^m / \mathbb{R} (1, \cdots, 1)$ (et dont l'ensemble des sommets s'identifie 
à $\mathbb{Z}^m / \mathbb{Z} (1, \cdots, 1)$).
Alors, l'injection naturelle entre les immeubles, notée $j$, est donnée par:
$$
j : A_{\mathbb{F}} \rightarrow A_{\mathbb{K}}, \, 
\overline{(x_1, \cdots, x_m)} \mapsto \overline{(2 x_1, \cdots, 2 x_m)}
$$

\end{propr}

\begin{démo}
On définit $j$, une application de $A_{\mathbb{F}}$ dans $A_{\mathbb{K}}$ par:
$$
j(\overline{(x_1, \cdots, x_m)}) = \overline{(2 x_1, \cdots, 2 x_m)}
$$
pour tout $\overline{(x_1, \cdots, x_m)}$ dans $A_{\mathbb{F}}$. 
Puis on prolonge $j$ sur $X_{\mathbb{F}}$ par ${\rm GL}_m (\mathcal{D})$-équivariance.
On a déjà vu que pour tout $x$ dans $\mathcal{D}^{\times}$, 
$v_{\Delta} (x) = 2 v_{\mathcal{D}} (x)$,
et que le sous-groupe ouvert compact maximal 
${\rm GL}_m (\mathcal{O}_{\Delta})$ est stable sous l'action du groupe de Galois de 
$\langle \sigma \rangle = {\rm Gal} (\mathbb{K}/\mathbb{F})$.\\
Avec des calculs analogues à ceux de la démonstration de la proposition \ref{InjectionImmeubleNRImpair}, 
on montre que $j$ est bien l'injection naturelle des immeubles.

\end{démo}

\subsubsection{Conditions de distinction.}

On utilise un raisonnement analogue à celui de la démonstration de la propriété \ref{RaisonnementDistinction} 
pour montrer la proposition suivante :

\begin{propo}
La représentation $\pi$ est ${\rm GL}_m (\mathcal{D})$-distinguée si et seulement s'il existe $g$ dans 
${\rm GL}_m (\mathcal{D}) \backslash {\rm GL}_m (\Delta) / \mathcal{K}$ tel que 
$g. j(s_0)$ appartient à $j(C_0)$ et:
$$
{\rm Hom}_{g^{-1} {\rm GL}_m (\mathcal{D}) g \cap \mathcal{K}} (\pi_0, \mathds{1}) \neq 0
$$
où $\mathcal{K} = \langle \varpi_{\mathbb{K}} \rangle {\rm GL}_m (\mathcal{O}_{\Delta})$
et $C_0$ est la chambre standard dans $A_{\mathbb{F}}$.\\
\end{propo}

\begin{propo}\label{PropositionNbSommets}
On utilise les mêmes notations que dans \ref{DefinitionSommet} pour la chambre standard $C_0$. 
Il y a exactement $\frac{m(m-1)}{2} +m$ sommets de $X_{\mathbb{K}}$ dans $j(C_0)$ qui sont les sommets:
$$
j(s_i) = [\mathcal{O}_{\Delta} \oplus \cdots \oplus \mathcal{O}_{\Delta}
\oplus \underbrace{\mathcal{P}_{\Delta}^2 \oplus \cdots \mathcal{P}_{\Delta}^2}_i]
$$
pour $i$ dans $\lbrace 0, \cdots, m-1 \rbrace$
et tous les sommets de la forme suivante:
$$
t_{k,l} = [\underbrace{\mathcal{O}_{\Delta} \oplus \cdots \oplus \mathcal{O}_{\Delta}}_{m-l} \oplus 
\underbrace{\mathcal{P}_{\Delta} \oplus \cdots \oplus \mathcal{P}_{\Delta}}_{l-k} 
\oplus \underbrace {\mathcal{P}_{\Delta}^2 \oplus \cdots \oplus \mathcal{P}_{\Delta}^2}_{k}]
$$
où $0 \leq k < l \leq m-1$ (il s'agit des milieux des arêtes 
$[j(s_k), j(s_l)]$). 
\end{propo}

\begin{démo} 
Pour tout $i$ dans $\lbrace 0, 1, \cdots, m-1 \rbrace$, on a $s_i = h_i. s_0$
où $h_i$ est la matrice diagonale 
$h_i = {\rm diag} (1, \cdots, 1, \underbrace{\varpi_{\mathcal{D}}, \cdots, \varpi_{\mathcal{D}}}_i)$. 
Comme $j$ est ${\rm GL}_m (\mathcal{D})$-équivariante, on a:
$$
j(s_i) = h_i . j(s_0 ) 
=  [\mathcal{O}_{\Delta} \oplus\cdots \oplus \mathcal{O}_{\Delta} 
\oplus \underbrace{\mathcal{P}_{\Delta}^2 \oplus \cdots \oplus \mathcal{P}_{\Delta}^2}_i] 
= \overline{(0, \cdots, 0,\underbrace{2, \cdots, 2}_i)}
$$

Soit $t \in j(C_0)$ tel que $t$ est aussi un sommet de $X_{\mathbb{K}}$.
Comme $j$ est une application affine, 
il existe $(\lambda_0, \lambda_1, \cdots, \lambda_{m-1})$ dans $(\mathbb{R}^{+})^m$ tels que:
$$
\sum_{i=0}^{m-1} \lambda_i = 1 \, \, \text{et} \, \, 
t = \sum_{i=0}^{m-1} \lambda_i j(s_i)
$$
On en déduit que:
$$
t = \overline{(0, 2 \lambda_{m-1}, 2 (\lambda_{m-2} + \lambda_{m-1}), \cdots, 
2 (\lambda_1 + \lambda_2 + \cdots + \lambda_{m-1}))}
$$
Fixons $i_0$ le plus grand indice $k$ dans $\{ 0, \cdots, m-1 \}$ tel que 
$\lambda_k \neq 0$.
Alors, la première coordonnée non nulle de $t$ est 
$2 (\lambda_{i_0} + \lambda_{i_0 +1} + \cdots + \lambda_{m-1})
 = 2 \lambda_{i_0}$. 
Comme $t$ est un sommet de $X_{\mathbb{K}}$, toutes ses coordonnées sont entières, donc  
$\lambda_{i_0} \in \{ \frac{1}{2}, 1 \}$.
Si $\lambda_{i_0} = 1$, on a $t = j(s_{i_0})$. 
Sinon, $\lambda_{i_0} = \frac{1}{2}$ et il existe au 
moins un autre coefficient $\lambda_k$ qui est non nul.
Fixons $k_0$ le plus grand coefficient tel que $\lambda_{k_0} \neq 0$ et $k_0 < i_0$. 
Alors, la deuxième coordonnée non nulle de $t$ est $2 (\lambda_{i_0} + \lambda_{k_0}) = 1 + 2 \lambda_{k_0}$.  
De même, comme $t$ est un sommet, $1 + 2 \lambda_{k_0} \in \mathbb{Z}$ et donc 
$\lambda_{k_0} = \frac{1}{2}$.
Ainsi:
$$
t = \frac{1}{2} j(s_{i_0}) + \frac{1}{2} j(s_{k_0}) 
= \overline{(0, \cdots, 0, \underbrace{1, \cdots, 1}_{i_0})} 
+ \overline{(0, \cdots, 0, \underbrace{1, \cdots, 1}_{k_0})} 
= \overline{(0, \cdots, 0, \underbrace{1, \cdots, 1}_{i_0 - k_0}, \underbrace{2, \cdots, 2}_{k_0})}
$$
Avec les notations de la proposition, il s'agit du sommet $t_{k_0, i_0}$, qui est aussi le milieu de 
l'arête $[ j(s_{i_0}), j(s_{k_0}) ]$.

\end{démo}

\begin{propo}
\begin{itemize}
\item[$\ast$] Si $m=1$, alors $\pi$ n'est pas ${\rm GL}_m (\mathcal{D})$-distinguée.
\item[$\ast$] Supposons que $m \geq 2$. 
Pour $l$ dans $\{ 1, \cdots, m-1 \}$, notons 
$g_l = {\rm diag} (\underbrace{1, \cdots, 1}_{m-l}, 
\underbrace{\varpi_{\Delta}, \cdots, \varpi_{\Delta}}_l )$
de sorte que 
$g_l . j(s_0) = t_{0,l} 
= [\underbrace{\mathcal{O}_{\Delta} \oplus \cdots \oplus \mathcal{O}_{\Delta}}_{m-l} \oplus 
\underbrace{\mathcal{P}_{\Delta} \oplus \cdots \oplus \mathcal{P}_{\Delta}}_{l}]$.\\
Alors, $\pi$ est ${\rm GL}_m (\mathcal{D})$-distinguée si et seulement s'il existe 
$l$ dans $\{ 1, \cdots, m-1 \}$ tel que:
$$
{\rm Hom}_{g_l^{-1} {\rm GL}_m (\mathcal{D}) g_l \cap \mathcal{K}} (\pi_0, \mathds{1}) \neq 0
$$
\end{itemize}
\end{propo}

\begin{démo}
On cherche tous les points $g.j(s_0)$ qui sont dans l'image de $C_0$.
Quitte à multiplier à gauche par un élément de ${\rm GL}_m (\mathcal{D})$, on peut s'intéresser uniquement 
au sommet $j(s_0)$ et, si $m \geq 2$,  
aux sommets qui sont sur une arête dont l'une des extrémités est $j(s_0)$.\\
Pour $j(s_0)$ on choisit $g= {\rm Id}$. 
On remarque que 
${\rm GL}_m (\mathcal{O}_{\mathcal{D}}) \subseteq {\rm GL}_m (\mathcal{D}) \cap \mathcal{K}$. Par suite:
$$
{\rm Hom}_{{\rm GL}_m (\mathcal{D}) \cap \mathcal{K}} (\pi_0, \mathds{1}) 
\subseteq {\rm Hom}_{{\rm GL}_m (\mathcal{O}_{\mathcal{D}})} (\pi_0, \mathds{1}) 
\simeq {\rm Hom}_{{\rm GL}_m (k_{\mathcal{D}})} (\overline{\gamma}_0, \mathds{1})
$$
La représentation $\overline{\gamma}_0$ étant une représentation irréductible 
de ${\rm GL}_m (k_{\mathcal{D}}) = {\rm GL}_m (k_{\Delta})$, l'espace d'entrelacement 
${\rm Hom}_{{\rm GL}_m (k_{\mathcal{D}})} (\overline{\gamma}_0, \mathds{1})$ est non trivial 
si et seulement si $\overline{\gamma}_0 = \mathds{1}$. 
Supposons que $\overline{\gamma}_0 = \mathds{1}$. Puisque $\overline{\gamma}_0$ est une représentation cuspidale, 
on a forcément $m=1$. Par conséquent, si $m=1$, 
$\overline{\gamma}_0 = \overline{\chi}$ est un caractère $k_{\mathbb{K}}$-régulier de 
$k_{\mathbb{K},n}^{\times}$ avec $n \geq 2$. Donc $\overline{\chi}$ ne peut pas être le 
caractère trivial. On déduit de tout ceci que, pour tout $m \in \mathbb{N}^{\times}$, on a :
$$
{\rm Hom}_{{\rm GL}_m (k_{\mathcal{D}})} (\overline{\gamma}_0, \mathds{1}) = 0
$$
et, si $m=1$, la représentation $\pi$ n'est pas distinguée.\\
Supposons à présent que $m \geq 2$. 
Il y a donc seulement $m-1$ sommets à considérer et il s'agit des sommets  
$g_1.j(s_0), g_2.j(s_0), \cdots, g_{m-1}.j(s_0)$.

\end{démo}

Rappelons un théorème démontré par Lusztig dans \cite{Lusztig} (théorème 3.4 page 62). 
Pour cela, nous allons introduire quelques notations afin de pouvoir énoncer ce résultat.
\begin{nota}
Notons $G$ un groupe réductif connexe défini sur un corps fini $\mathbb{F}_q$ de cardinal $q$ impair. 
Soit $\theta : G \rightarrow G$ une involution définie sur $\mathbb{F}_q$. 
Soit $K$ un sous-groupe fermé de $G^{\theta}$, l'ensemble des points fixes de $G$ sous l'action de $\theta$, 
tel que $K$ est défini sur $\mathbb{F}_q$ et contient $(G^{\theta})^{\circ}$, la composante de l'unité de 
$G^{\theta}$. 
Soit $T$ un tore maximal de $G$ défini sur $\mathbb{F}_q$ et 
$\lambda : T(\mathbb{F}_q) \rightarrow \overline{\mathbb{Q}_l}^{\times}$ un caractère. 
Pour tout $h$ dans $G$, on définit $\lambda^h : (h^{-1} T h)(\mathbb{F}_q) \rightarrow \overline{\mathbb{Q}_l}^{\times}, 
x \mapsto \lambda (h x h^{-1})$. On note 
$\Xi = \{ f \in G : \theta (f^{-1} T f) = f^{-1} T f \}$ vu comme un ensemble de $T-K$ doubles classes. 
Si $T^{'}$ est un tore maximal de $G$ défini sur $\mathbb{F}_q$ et $\theta$-stable, on note 
$Z$ le centralisateur de $(T^{'} \cap K)^{\circ}$ et, pour tout $t$ dans $T^{'} \cap K$, 
$Z_t = Z \cap (Z_{G} (t))^{\circ}$. On définit alors l'application $\varepsilon_{T^{'}}$ par:
$$
\varepsilon_{T^{'}} : T^{'} \cap K \rightarrow \{ -1, 1 \}, \, 
t \mapsto \sigma (Z) \sigma (Z_t)
$$
où, pour tout 
groupe réductif connexe $M$ défini sur $\mathbb{F}_q$, $\sigma (M)$ vaut $-1$ si le $\mathbb{F}_q$-rang de $M$ est 
impair et vaut $1$ sinon. 
On définit enfin l'ensemble $\widetilde{\Xi}$:
$$
\widetilde{\Xi} = \{ h \in \Xi : \lambda^h_{\vert (h^{-1} T h \cap K)(\mathbb{F}_q)} = \varepsilon_{h^{-1} T h} \}
$$
\end{nota}

\begin{theo}\label{TheoremeLusztig} (\cite{Lusztig})
Soit $G$ un groupe réductif connexe défini sur un corps fini $\mathbb{F}_q$ de cardinal $q$ impair. 
Soit $\theta : G \rightarrow G$ une involution définie sur $\mathbb{F}_q$. 
On fixe $K$ un sous-groupe fermé de $G^{\theta}$, l'ensemble des points fixes de $G$ sous l'action de $\theta$, 
tel que $K$ est défini sur $\mathbb{F}_q$ et contient $(G^{\theta})^{\circ}$, la composante de l'unité de 
$G^{\theta}$. 
Soit $T$ un tore maximal de $G$ défini sur $\mathbb{F}_q$ et 
$\lambda : T(\mathbb{F}_q) \rightarrow \overline{\mathbb{Q}_l}^{\times}$ un caractère 
(où $l$ est un nombre premier ne divisant pas $q$). 
Enfin, on note ${\rm R}_T^{\lambda}$ la représentation virtuelle de $G(\mathbb{F}_q)$ attachée à $(T, \lambda)$. 
Alors:
$$
\frac{1}{\vert K(\mathbb{F}_q) \vert} \sum_{g \in K (\mathbb{F}_q)} {\rm Tr} (g, {\rm R}_T^{\lambda}) 
= \sum_{f \in \Xi} r(f)
$$
où  
$r(f)$ est un entier appartenant à $\{ -1, 0, 1 \}$ qui est non nul si et seulement si 
$f \in \widetilde{\Xi}$. Dans ce cas, on a:
$$
r(f) = \sigma (T) \sigma (Z (((f^{-1} T f)\cap K)^{\circ}))
$$
\end{theo}

\begin{lem}\label{LemmeDistinctionLevi}
Si la représentation $\overline{\gamma}_0$ est distinguée par un sous-groupe de 
Lévi $L_l$ isomorphe à ${\rm GL}_{m-l} (k_{\Delta}) \times {\rm GL}_{l} (k_{\Delta})$, alors $m$ est pair et 
$l = \frac{m}{2}$.
\end{lem}

\begin{démo}
On peut voir $k_{\Delta,m}$ comme un $k_{\Delta,m}$-espace vectoriel de dimension $1$ et l'identifier à $k_{\Delta}^m$. 
Il existe donc une injection de $k_{\Delta}$-algèbres:
$$
\psi : k_{\Delta,m} \hookrightarrow {\rm M}_m (k_{\Delta})
$$
On note $T$ l'image de $k_{\Delta,m}^{\times}$ dans ${\rm GL}_m (k_{\Delta})$ via l'application $\psi$. 
Soit $l$ dans $\{ 1, \cdots, m-1\}$. 
Soit $L_l$ le sous-groupe de Lévi :
$$
L_l = \left( \begin{array}{c|c}
{\rm GL}_{m-l} (k_{\Delta}) & 0\\
\hline 
0 & {\rm GL}_l (k_{\Delta})\\
            \end{array} \right)
$$
On peut voir $( {\rm GL}_m (k_{\Delta}), L_l )$ comme un espace symétrique. 
En effet, posons:
$$
w_l =  \left( \begin{array}{c|c}
I_{m-l} & 0\\
\hline 
0 & - I_l \\
            \end{array} \right)
$$
et $\tau : {\rm GL}_m (k_{\Delta}) \rightarrow {\rm GL}_m (k_{\Delta}), \, x \mapsto w_l x w_l^{-1}$. 
Il est clair que $w_l^{-1} = w_l$, donc $\tau$ est une involution. 
Un calcul rapide nous montre que les points fixes de $\tau$ sont exactement les éléments de $L_l$.\\
On définit $\Xi = \lbrace g \in {\rm GL}_m (k_{\Delta}) : \tau (g^{-1} T g) = g^{-1} T g \rbrace$.
Supposons que $\overline{\gamma}_0$ est $L_l$-distinguée, alors d'après l'article de Lusztig \cite{Lusztig} 
(cf. le théorème \ref{TheoremeLusztig}, ici $\overline{\gamma}_0$ correspond à $R_T^{\overline{\chi}}$), 
l'ensemble $\Xi$ est nécessairement non vide. 
Il existe donc $g$ dans ${\rm GL}_m (k_{\Delta})$, tel que $g^{-1} T g$ est $\tau$-stable. 
On peut voir l'application $\tau : {\rm M}_m (k_{\Delta}) \rightarrow {\rm M}_m (k_{\Delta}), \, x \mapsto w_l x w_l^{-1}$ comme 
un automorphisme de $k_{\Delta}$-algèbres et la restriction de $\tau$ à ${\rm GL}_m (k_{\Delta})$ comme un automorphisme de 
groupes. 
Posons $l_1 = g^{-1} T g \cup \lbrace 0 \rbrace = g^{-1} k_{\Delta,m} g$. 
Alors $l_1$ est une sous-$k_{\Delta}$-algèbre de 
${\rm M}_m (k_{\Delta})$ isomorphe à $k_{\Delta,m}$, donc est également un corps. 
Ainsi $\tau_{\vert l_1} : l_1 \rightarrow l_1$ est un automorphisme de corps ($k_{\Delta}$-linéaire) d'ordre $1$ ou $2$. 
Notons $l_1^{\tau}$ les points fixes de $\tau_{\vert l_1}$. Ou bien $l_1^{\tau} = l_1$ ou bien 
il existe $l_0$ sous-extension quadratique (de corps) de $l_1$ telle que $l_1^{\tau} = l_0$ (dans ce cas, puisque 
$k_{\Delta} \subseteq l_1$ et que les éléments de $k_{\Delta}$ sont fixés par $\tau$, on a 
$k_{\Delta} \subseteq l_0 \subseteq l_1$).\\
Dans les deux cas, il existe un corps $l_0$ tel que $k_{\Delta} \subseteq l_0 \subseteq l_1$, $l_1^{\tau} = l_0$ et 
$[l_1 : l_0 ] \leq 2$. 
Puisque $\tau (x) = x$ pour tout $x$ dans $l_0^{\times}$, on a $l_0^{\times} \subseteq L_l$ et 
$l_0$ peut être vu comme un sous-anneau 
(et même une sous-$k_{\Delta}$-algèbre) de ${\rm M}_{m-l} (k_{\Delta}) \times {\rm M}_l (k_{\Delta})$. 
Soit $\psi_0 : l_0 \hookrightarrow {\rm M}_{m-l} (k_{\Delta}) \times {\rm M}_l (k_{\Delta})$ 
une injection de $k_{\Delta}$-algèbres (unitaires).
Posons $p_1 : {\rm M}_{m-l} (k_{\Delta}) \times {\rm M}_l (k_{\Delta}) \rightarrow {\rm M}_{m-l} (k_{\Delta}), \, 
(X_1, X_2) \mapsto X_1$. 
De même, on définit l'application 
$p_2 : {\rm M}_{m-l} (k_{\Delta}) \times {\rm M}_l (k_{\Delta}) \rightarrow {\rm M}_{l} (k_{\Delta}), \, 
(X_1, X_2) \mapsto X_2$. Il s'agit de deux 
morphismes surjectifs de $k_{\Delta}$-algèbres unitaires. 
Alors $p_1 \circ \psi_0 : l_0 \rightarrow {\rm M}_{m-l} (k_{\Delta})$ est un morphisme de 
$k_{\Delta}$-algèbres unitaires, injectif (car $l_0$ est un corps). Il en est de même pour $p_2 \circ \psi_0$.\\
Ainsi, on peut supposer que $l_0 \subseteq {\rm M}_{m-l} (k_{\Delta})$, et donc $[l_0 : k_{\Delta}] \leq m-l$.
On montre de même que $[l_0 : k_{\Delta}] \leq l$.\\
Maintenant, si $l_1 = l_0$, on a $[l_0 : k_{\Delta}] = m \leq l$ et $m \leq m-l$, ce qui est impossible.
On en déduit que $l_1 \neq l_0$, donc $l_1 / l_0$ est une extension quadratique et $m$ est pair. 
De plus, $[l_0 : k_{\Delta}] = m/2 \leq l$ et $m/2 \leq m-l$ donc $m-l = l = \frac{m}{2}$.

\end{démo}

\begin{propo}\label{NonDistinctionCuspidalTotRamnImpair}
Si $n$ est impair alors $\pi$ n'est pas ${\rm GL}_m (\mathcal{D})$-distinguée.
\end{propo}

\begin{démo}
Supposons que $n$ est impair, alors nécessairement $m$ est également impair. 
Si $m=1$, on a déjà vu que $\pi$ n'est pas distinguée.\\
Supposons à présent que $m \geq 2$. 
Pour regarder si $\pi$ est ou non distinguée, il y a $m-1$ intersections 
à calculer:  $g_l^{-1} {\rm GL}_m (\mathcal{D}) g_l \cap \mathcal{K}$ 
où $l \in \lbrace 1, \cdots, m-1 \rbrace$. 
On sait que ${\rm GL}_m (\mathcal{D}) \cap g_l (\langle \varpi_{\mathbb{K}} \rangle 
{\rm GL}_m (\mathcal{O}_{\Delta}))g_l^{-1}$ est inclus dans  
${\rm GL}_m (\mathcal{D}) \cap g_l (\langle \varpi_{\Delta} \rangle 
{\rm GL}_m (\mathcal{O}_{\Delta}))g_l^{-1}$, 
or $\langle \varpi_{\Delta} \rangle 
{\rm GL}_m (\mathcal{O}_{\Delta}) = \mathcal{K}_{{\rm GL}_m (\Delta)} (j(s_0))$ 
(on utilise \cite{BushnellFröhlich} ainsi que les notations \ref{NotationsImmeuble}), donc:
$$
{\rm GL}_m (\mathcal{D}) \cap g_l (\langle \varpi_{\mathbb{K}} \rangle 
{\rm GL}_m (\mathcal{O}_{\Delta}))g_l^{-1} \subseteq 
{\rm GL}_m (\mathcal{D}) \cap \mathcal{K}_{{\rm GL}_m (\Delta)} (g_l.j(s_0))
$$
Nous allons commencer par déterminer, pour $l$ dans $\lbrace 1, \cdots, m-1 \rbrace$, l'intersection suivante:
$$
{\rm GL}_m (\mathcal{D})  \cap \mathcal{K}_{{\rm GL}_m (\Delta)} (g_l.j(s_0))
$$
Comme $g_l.j(s_0) \in j(X_{\mathbb{F}})$, il existe $t_l \in X_{\mathbb{F}}$ tel que 
$g_l.j(s_0) = j(t_l)$. 
En utilisant l'injectivité et la ${\rm GL}_m (\mathcal{D})$-équivariance de $j$, on vérifie facilement que:
$$
{\rm GL}_m (\mathcal{D})  \cap \mathcal{K}_{{\rm GL}_m (\Delta)} (j(t_l)) 
= \mathcal{K}_{{\rm GL}_m (\mathcal{D})} (t_l)
$$
Nous avons vu lors de la démonstration de la proposition \ref{PropositionNbSommets}  
que $t_l$ est en fait le milieu de l'arête $[s_0, s_l ]$.
Nous devons à présent déterminer le stabilisateur de $t_l$.
Comme $t_l$ est le milieu de l'arête $[s_0, s_l]$, on a $\mathcal{K}_{{\rm GL}_m (\mathcal{D})} (t_l) 
= \mathcal{K}_{{\rm GL}_m (\mathcal{D})} (\mathcal{L})$
où $\mathcal{L}$ est la chaîne de période $2$: $\mathcal{L} = s_0 \cup s_l$
avec:
$$
s_0 = [L_0] = [\mathcal{O}_{\mathcal{D}} \oplus \cdots \oplus \mathcal{O}_{\mathcal{D}}]
$$
et:
$$
s_l = [L_l] = [\underbrace{\mathcal{O}_{\mathcal{D}} \oplus \cdots \oplus \mathcal{O}_{\mathcal{D}}}_{m-l} 
\oplus \underbrace{\mathcal{P}_{\mathcal{D}} \oplus \cdots \oplus \mathcal{P}_{\mathcal{D}}}_{l} ]
$$
Ainsi, $\mathcal{L} = (M_k)_{k \in \mathbb{Z}}$ où:
$$
M_ 0 = L_0 = \mathcal{O}_{\mathcal{D}} \oplus \cdots \oplus \mathcal{O}_{\mathcal{D}}
$$
$$
M_1 = L_l = \underbrace{\mathcal{O}_{\mathcal{D}} \oplus \cdots \oplus \mathcal{O}_{\mathcal{D}}}_{m-l} 
\oplus \underbrace{\mathcal{P}_{\mathcal{D}} \oplus \cdots \oplus \mathcal{P}_{\mathcal{D}}}_{l} 
$$
et pour tout $k$ dans $\mathbb{Z}$, $M_{k+2} = \varpi_{\mathcal{D}} M_k$.
Soit $(d_1, d_2)$ une partition associée à $\mathcal{L}$, alors:
$$
d_1 = {\rm dim}_{k_{\mathcal{D}}} (M_1 / M_2) = m-l \, \, \text{et} \, \, 
d_2 = {\rm dim}_{k_{\mathcal{D}}} (M_0 / M_1) = l
$$
On en déduit que:
$$
\mathcal{A} (\mathcal{L})^{\times} = \left( \begin{array}{c|c}
{\rm M}_{m-l, m-l} (\mathcal{O}_{\mathcal{D}}) & {\rm M}_{m-l, l} (\mathcal{O}_{\mathcal{D}}) \\
\hline
{\rm M}_{l, m-l} (\mathcal{P}_{\mathcal{D}})  & {\rm M}_{l, l} (\mathcal{O}_{\mathcal{D}}) \\
\end{array} \right)^{\times}
$$
et $\mathcal{K}_{{\rm GL}_m (\mathcal{D})} (\mathcal{L}) 
= \langle \Pi_{\mathcal{L}} \rangle \mathcal{A} (\mathcal{L})^{\times}$
où pour tout $i$ dans $\mathbb{Z}$, $\Pi_{\mathcal{L}} M_i = M_{i+ \nu}$
avec $\nu \in \mathbb{N}^{\times}$ le plus petit possible.
Pour tout $i$ dans $\mathbb{Z}$, $\varpi_{\mathcal{D}} M_i = M_{i+ 2}$
donc $\nu \in \{1,2 \}$. 
Or, si $\Pi_{\mathcal{L}} M_i = M_{i+ \nu}$, on a $d_i = d_{i+ \nu}$ ($(d_i)$ est $\nu$-périodique), donc
 si $\nu = 1$, on a $d_2 = d_1$ et $n = d_1 + d_2 = 2 d_1$, ce qui est impossible car on est dans le cas où 
$n$ est impair.
On en déduit que $\nu = 2$ et que l'on peut choisir $\Pi_{\mathcal{L}} = \varpi_{\mathcal{D}} {\rm Id}$.
Finalement:
$$
\mathcal{K}_{{\rm GL}_m (\mathcal{D})} (\mathcal{L}) 
= \langle \varpi_{\mathcal{D}} \rangle \mathcal{A} (\mathcal{L})^{\times} 
$$

On va à présent caractériser $g_l (\langle \varpi_{\mathbb{K}} \rangle {\rm GL}_m (\mathcal{O}_{\Delta})) g_l^{-1}$ 
dans $g_l (\langle \varpi_{\Delta} \rangle {\rm GL}_m (\mathcal{O}_{\Delta})) g_l^{-1}$. 
On remarque que pour tout $x$ dans ${\rm GL}_m (\mathcal{O}_{\Delta})$, on a:
$$
{\rm Nrd} (g_l x g_l^{-1}) = {\rm Nrd} (x) \in \mathcal{O}_{\mathbb{K}}^{\times}
$$
De plus:
$$
g_l (\langle \varpi_{\mathbb{K}} \rangle 
{\rm GL}_m (\mathcal{O}_{\Delta}))g_l^{-1} 
= \langle \varpi_{\mathbb{K}} \rangle 
g_l {\rm GL}_m (\mathcal{O}_{\Delta})g_l^{-1}
$$
et:
$$
g_l (\langle \varpi_{\Delta} \rangle 
{\rm GL}_m (\mathcal{O}_{\Delta}))g_l^{-1} 
= \langle \varpi_{\Delta} \rangle 
g_l {\rm GL}_m (\mathcal{O}_{\Delta})g_l^{-1}
$$
Des calculs d'intersections analogues à ceux de la démonstration de la proposition 
\ref{PropositionIntersection}, nous permettent de montrer que:
$$
v_{\mathbb{K}} ({\rm Nrd} (\varpi_{\Delta} I_m)) 
= m = v_{\mathbb{F}} ({\rm Nrd} (\varpi_{\mathcal{D}} I_m)) \, \, \text{et} \, \,
v_{\mathbb{K}} ({\rm Nrd} (\varpi_{\mathbb{K}} I_m)) 
= n = v_{\mathbb{F}} ({\rm Nrd} (\varpi_{\mathbb{F}} I_m)) 
$$
On en déduit que:
$$
\langle \varpi_{\mathbb{K}} \rangle 
g_l {\rm GL}_m (\mathcal{O}_{\Delta})g_l^{-1} 
= \{ g \in \mathcal{K}_{{\rm GL}_m (\Delta)} (j(t_l)) : n \, \, 
\text{divise} \, \, v_{\mathbb{K}} ({\rm Nrd} (g)) \}
$$
Par conséquent, $x$ appartient à ${\rm GL}_m (\mathcal{D}) \cap \langle \varpi_{\mathbb{K}} \rangle 
g_l {\rm GL}_m (\mathcal{O}_{\Delta})g_l^{-1}$ si et seulement si $x$ appartient à la fois à 
${\rm GL}_m (\mathcal{D})$ et $\mathcal{K}_{{\rm GL}_m (\Delta)} (j(t_l))$ et 
$n$ divise $v_{\mathbb{K}} ({\rm Nrd} (g))$, i.e  
si $x$ appartient à $\mathcal{K}_{{\rm GL}_m (\mathcal{D})} (t_l)$ et 
$n$ divise $2 v_{\mathbb{F}} ({\rm Nrd} (g))$.\\
Avec le calcul précédent, on a $g \in {\rm GL}_m (\mathcal{D}) \cap \langle \varpi_{\mathbb{K}} \rangle 
g_l {\rm GL}_m (\mathcal{O}_{\Delta})g_l^{-1}$ si et seulement si 
$g$ appartient à $\langle \varpi_{\mathcal{D}} \rangle \mathcal{A} (\mathcal{L})^{\times}$
et $n$ divise $2 v_{\mathbb{F}} ({\rm Nrd} (g))$.
Si $g$ appartient à $\langle \varpi_{\mathcal{D}} \rangle \mathcal{A} (\mathcal{L})^{\times}$, 
$g = \varpi_{\mathcal{D}}^r u$ ($r \in \mathbb{Z}$ et $u \in \mathcal{A} (\mathcal{L})^{\times}$), alors 
$v_{\mathbb{F}} ({\rm Nrd} (g)) = v_{\mathbb{F}} ({\rm Nrd} (\varpi_{\mathcal{D}}^r)) 
= rm$ et $n$ divise $2 v_{\mathbb{F}} ({\rm Nrd} (g))$ si et seulement si $d$ divise $r$. 
On en déduit que 
${\rm GL}_m (\mathcal{D}) \cap \langle \varpi_{\mathbb{K}} \rangle 
g_l {\rm GL}_m (\mathcal{O}_{\Delta})g_l^{-1} 
= \langle \varpi_{\mathbb{F}} \rangle \mathcal{A} (\mathcal{L})^{\times}$.\\
Et donc 
$g_l^{-1} {\rm GL}_m (\mathcal{D}) g_l \cap \langle \varpi_{\mathbb{K}} \rangle 
{\rm GL}_m (\mathcal{O}_{\Delta}) = \langle \varpi_{\mathbb{F}} \rangle 
g_l^{-1} \mathcal{A} (\mathcal{L})^{\times} g_l$.
Un calcul rapide nous permet de constater que 
la réduction de $g_l^{-1} \mathcal{A} (\mathcal{L})^{\times} g_l$  dans 
${\rm GL}_m (k_{\mathcal{D}}) = {\rm GL}_m (k_{\Delta})$ peut 
être vu comme le sous-groupe de Lévi $L_l$:
$$
L_l = \left( \begin{array}{c|c}
{\rm GL}_{m-l} (k_{\Delta}) & 0\\
\hline
0 & {\rm GL}_l (k_{\Delta})\\
\end{array} \right)
$$
On en déduit que $\pi$ est distinguée si et seulement s'il existe $l$ dans 
$\{ 1, \cdots, m-1 \}$ tel que:
$$
{\rm Hom}_{\langle \varpi_{\mathbb{F}} \rangle 
g_l^{-1} \mathcal{A} (\mathcal{L})^{\times} g_l} (\pi_0, \mathds{1}) \neq 0
$$
i.e si $\chi$ est trivial sur $\mathbb{F}^{\times}$ et:
$$
{\rm Hom}_{g_l^{-1} \mathcal{A} (\mathcal{L})^{\times} g_l} (\pi_0, \mathds{1}) 
\simeq {\rm Hom}_{L_l} (\overline{\gamma}_0, \mathds{1}) \neq 0
$$
Or, comme $m$ est impair, d'après le lemme \ref{LemmeDistinctionLevi}, $\overline{\gamma}_0$ ne 
peut pas être distinguée par le sous-groupe de Lévi $L_l$. 
Par conséquent, $\pi$ ne peut pas être distinguée.

\end{démo}

\begin{propo}
Si $n$ est pair (alors $m$ est pair), on a:
$$
{\rm Hom}_{{\rm GL}_m (\mathcal{D})} (\pi, \mathds{1}) 
\simeq {\rm Hom}_{\langle \varpi_{\mathbb{K}} w \rangle L_{\frac{m}{2}}} (\pi_0, \mathds{1})
$$
avec:
$$
w = \left( \begin{array}{c|c}
0_{\frac{m}{2}} & I_{\frac{m}{2}} \\
\hline
I_{\frac{m}{2}}  & 0_{\frac{m}{2}} \\
\end{array} \right) \, \, \text{et} \, \, 
L_{\frac{m}{2}} = \left( \begin{array}{c|c}
{\rm GL}_{m/2} (k_{\Delta}) & 0 \\
\hline
0  & {\rm GL}_{m/2} (k_{\Delta}) \\
\end{array} \right)
$$
Par conséquent, la représentation $\pi$ est ${\rm GL}_m (\mathcal{D})$-distinguée si et seulement si 
la représentation $\pi_0$ est $\langle \varpi_{\mathbb{K}} w \rangle L_{\frac{m}{2}}$-distinguée.
\end{propo}

\begin{démo}
On suppose que $n$ est pair. Comme $n = md$ avec $d$ impair, $m$ est aussi pair.
Comme dans la proposition précédente, il y a $m-1$ intersections à calculer.
On remarque de même que si $l$ est dans $\lbrace 1, \cdots, m-1 \rbrace$ avec 
$l \neq \frac{m}{2}$, on a:
$$
{\rm Hom}_{g_l^{-1} {\rm GL}_m (\mathcal{D}) g_l \cap \mathcal{K} } 
(\pi_0, \mathds{1}) = {\rm Hom}_{\langle \varpi_{\mathbb{F}} \rangle L_l} (\pi_0, \mathds{1})
$$
Alors, comme précédemment, on montre que nécessairement:
$$
{\rm Hom}_{\langle \varpi_{\mathbb{F}} \rangle L_l} (\pi_0, \mathds{1}) = 0
$$
On déduit de tout ceci que $\pi$ est ${\rm GL}_m (\mathcal{D})$-distinguée si 
et seulement si:
$$
{\rm Hom}_{g_{m/2}^{-1} {\rm GL}_m (\mathcal{D}) g_{m/2} \cap \mathcal{K} } 
(\pi_0, \mathds{1}) \neq 0
$$
En utilisant un raisonnement analogue à \ref{RaisonnementDistinction} (utilisation de la formule de restriction de 
Mackey, de la réciprocité de Frobenius pour l'induction compacte ainsi que \cite{HakimMurnaghan2}, qui est rappelé au 
théorème \ref{ResultatHakimMurnaghan}) 
on montre l'isomorphisme de $\mathbb{C}$-espaces vectoriels:
$$
{\rm Hom}_{{\rm GL}_m (\mathcal{D})} (\pi, \mathds{1}) 
\simeq {\rm Hom}_{g_{m/2}^{-1} {\rm GL}_m (\mathcal{D}) g_{m/2} \cap \mathcal{K} } 
(\pi_0, \mathds{1})
$$
On commence par calculer 
${\rm GL}_m (\mathcal{D}) \cap g_{m/2} \mathcal{K} g_{m/2}^{-1}$.
On remarque que:
$$
{\rm GL}_m (\mathcal{D}) \cap g_{m/2} \mathcal{K} g_{m/2}^{-1} 
\subseteq {\rm GL}_m (\mathcal{D}) \cap \mathcal{K}_{{\rm GL}_m (\Delta)} (g_{m/2}. j(s_0)) 
= \mathcal{K}_{{\rm GL}_m (\mathcal{D})} (t_{m/2})
$$
En raisonnant de façon analogue, on vérifie que:
$$
{\rm GL}_m (\mathcal{D})  \cap \mathcal{K}_{{\rm GL}_m (\Delta)} (g_{m/2}.j(s_0)) 
= \langle \Pi_{\mathcal{L}} \rangle \mathcal{A} (\mathcal{L})^{\times}
$$
où:
$$
\mathcal{A} (\mathcal{L})^{\times} 
= \left( \begin{array}{c|c}
{\rm M}_{m/2} (\mathcal{O}_{\mathcal{D}}) & {\rm M}_{m/2} (\mathcal{O}_{\mathcal{D}}) \\
\hline
{\rm M}_{m/2} (\mathcal{P}_{\mathcal{D}})  & {\rm M}_{m/2} (\mathcal{O}_{\mathcal{D}}) \\
\end{array} \right)^{\times} \, \, \text{et} \, \, 
\Pi_{\mathcal{L}} = \left( \begin{array}{c|c}
0_{\frac{m}{2}} & I_{\frac{m}{2}} \\
\hline
\varpi_{\mathcal{D}} I_{\frac{m}{2}}  & 0_{\frac{m}{2}} \\
\end{array} \right)
$$
De plus,  
$v_{\mathbb{F}} ({\rm Nrd} (\Pi_{\mathcal{L}}^2)) 
= v_{\mathbb{F}} ({\rm Nrd} (\varpi_{\mathcal{D}} I_{m})) = m$ 
donc 
$v_{\mathbb{F}} ({\rm Nrd} (\Pi_{\mathcal{L}})) = \frac{m}{2}$.
Comme au premier cas, on a $x \in 
g_{m/2}^{-1} {\rm GL}_m (\mathcal{D}) g_{m/2} \cap \langle \varpi_{\mathbb{K}} \rangle 
{\rm GL}_m (\mathcal{O}_{\Delta})$ si et seulement si 
$x$ appartient à $\langle \Pi_{\mathcal{L}} \rangle \mathcal{A} (\mathcal{L})^{\times}$ 
et $n$ divise $2 v_{\mathbb{F}} ({\rm Nrd} (x))$.
Si $x \in \langle \Pi_{\mathcal{L}} \rangle \mathcal{A} (\mathcal{L})^{\times}$, il existe 
$r \in \mathbb{Z}$ et $u$ dans $\mathcal{A} (\mathcal{L})^{\times}$ tels que 
$x = \Pi_{\mathcal{L}}^r u$.
Ainsi 
$v_{\mathbb{F}} ({\rm Nrd} (x)) = r \times \frac{m}{2}$.
Par conséquent, $n$ divise $2 v_{\mathbb{F}} ({\rm Nrd} (x))$ si et seulement si 
$n$ divise $rm$, i.e si $\frac{n}{m} = d$ divise $r$.
On en déduit que:
$$
{\rm GL}_m (\mathcal{D}) \cap g_{m/2} \mathcal{K} g_{m/2}^{-1} 
= \langle \Pi_{\mathcal{L}}^d \rangle \mathcal{A} (\mathcal{L})^{\times}
$$
On a donc:
$$
g_{m/2}^{-1} {\rm GL}_m (\mathcal{D}) g_{m/2} \cap \langle \varpi_{\mathbb{K}} \rangle 
{\rm GL}_m (\mathcal{O}_{\Delta}) 
= (g_{m/2}^{-1} \langle \Pi_{\mathcal{L}}^d \rangle g_{m/2}) ( g_{m/2}^{-1} \mathcal{A} (\mathcal{L})^{\times} g_{m/2})
$$
Un calcul rapide nous montre que:
$$
\overline{g_{m/2}^{-1} \mathcal{A} (\mathcal{L})^{\times} g_{m/2}} 
\simeq L_{m/2} = \left( \begin{array}{c|c}
{\rm GL}_{m/2} (k_{\Delta}) & 0\\
\hline
0 & {\rm GL}_{m/2} (k_{\Delta})\\
\end{array} \right)
$$
De plus, comme 
$g_{m/2}^{-1} \Pi_{\mathcal{L}}^d g_{m/2} 
= (g_{m/2}^{-1} \Pi_{\mathcal{L}} g_{m/2})^d$
et que 
$g_{m/2}^{-1} \Pi_{\mathcal{L}} g_{m/2} = \varpi_{\Delta} w$
où:
$$
w = \left( \begin{array}{c|c}
0_{m/2} & I_{m/2}\\
\hline
I_{m/2} & 0_{m/2}\\
\end{array} \right)
$$
et $w^2 = I_m$, on a  
$(g_{m/2}^{-1} \Pi_{\mathcal{L}} g_{m/2})^d 
= \varpi_{\Delta}^d w^d = \varpi_{\mathbb{K}} (w^2)^{\frac{d-1}{2}} w 
= \varpi_{\mathbb{K}} w$. 
Ainsi:
$$
g_{m/2}^{-1} {\rm GL}_m (\mathcal{D}) g_{m/2} \cap \langle \varpi_{\mathbb{K}} \rangle 
{\rm GL}_m (\mathcal{O}_{\Delta}) = \langle \varpi_{\mathbb{K}} w \rangle 
( g_{m/2}^{-1} \mathcal{A} (\mathcal{L})^{\times} g_{m/2})
$$

Finalement, $\pi$ est ${\rm GL}_m (\mathcal{D})$-distinguée si et seulement si:
$$
{\rm Hom}_{\langle \varpi_{\mathbb{K}} w \rangle L_{\frac{m}{2}}} (\pi_0, \mathds{1}) \neq 0
$$
\end{démo}

\begin{nota}
On suppose que $n$ est pair (i.e que $m$ est pair).\\
On note $l$ le corps résiduel de $\mathbb{K}_n$ (extension de degré $m$ de $k_{\Delta}$) 
et $l_0$ celui de $\mathbb{K}_{n/2}$ (extension de degré $m/2$ de $k_{\Delta}$). 
On note $\Lambda = \mathbb{K}_d$ (resp. $\mathbb{L} = \mathbb{F}_d$), alors le corps résiduel $k_{\Lambda}$ de $\Lambda$ 
(resp. $k_{\mathbb{L}}$) s'identifie à $k_{\Delta}$ (resp. $k_{\mathcal{D}}$).  
On définit $\mathbb{E}$, extension totalement ramifiée de degré $2$ de $\mathbb{F}_{n/2}$,  
par le diagramme d'extensions de corps suivant:
$$
\xymatrix{ l \ar@{-}[d]_{2} \\
           l_0 \ar@{-}[d]_{m/2} \\
           k_{\mathcal{D}} = k_{\Delta} \ar@{-}[d]_{d} \\
           k = k_{\mathbb{F}} = k_{\mathbb{K}} }
\quad \quad \quad 
\xymatrix{ & \Lambda_m = \mathbb{K}_n & \\ 
         \Lambda_{m/2} = \mathbb{K}_{n/2} \ar@{-}[ru]_{{\rm n.r}}^2 \ar@{-}[rd]_2^{{\rm t.r}} \ar@{-}[d]_{n/2}^{{\rm n.r}} 
         & \mathbb{F}_n  \ar@{-}[u]_{{\rm t.r}}^2 \ar@{-}[d]_2^{{\rm n.r}} 
         & \mathbb{E} \ar@{-}[lu]_{{\rm n.r}}^2 \ar@{-}[ld]_2^{{\rm t.r}}\\
         \mathbb{K} \ar@{-}[rdd]_2^{{\rm t.r}} & \mathbb{F}_{n/2} \ar@{-}[d]_{m/2}^{{\rm n.r}} & \\
          & \mathbb{L} \ar@{-}[d]_{d}^{{\rm n.r}} & \\
          & \mathbb{F} & } 
$$
\end{nota}

\begin{lem}\label{DistinctionParLeviCuspidale}
On suppose que $m$ est pair. 
Alors, la représentation cuspidale $\overline{\gamma}_0$ est $L_{m/2}$-distinguée si 
et seulement si $\overline{\chi}$ est trivial sur $l_0^{\times}$. 
De plus, si $\overline{\gamma}_0$ est $L_{m/2}$-distinguée on a:
$$
{\rm dim}_{\mathbb{C}} ({\rm Hom}_{L_{m/2}} (\overline{\gamma}_0, \mathds{1})) = 1
$$
\end{lem}

\begin{démo}
Pour la démonstration, on peut se référer à l'article \cite{HakimMurnaghan1}, proposition 6.1 page 1874.

\end{démo}

\begin{rmq}
Une conséquence de ce lemme est que si $\pi$ est ${\rm GL}_m (\mathcal{D})$-distinguée, alors 
$m$ est pair, $\chi$ est trivial sur $\mathbb{F}^{\times}$ et $\overline{\chi}$ 
est trivial sur $l_0^{\times}$.
\end{rmq}

\begin{nota}
On fixe $\varepsilon$ dans $l_0^{\times}$ tel que 
$\sqrt{\varepsilon} \in l^{\times} \backslash l_0^{\times}$.
Ainsi 
$l = l_0 (\sqrt{\varepsilon})$ 
et $(1, \sqrt{\varepsilon})$ est une $l_0$-base de $l$.
Comme $l_0$ est une extension de degré $m/2$ de $k_{\Delta}$, on peut voir $l_0$ comme 
un $k_{\Delta}$-espace vectoriel de dimension $m/2$.
On a donc une injection:
$$
i: l_0 \hookrightarrow {\rm End}_{k_{\Delta}} (l_0) \simeq {\rm M}_{m/2} (k_{\Delta}), \, 
x \mapsto \varphi_x: y \mapsto xy
$$
On en déduit que $l_0^{\times}$ s'injecte dans le sous groupe de Lévi 
$L_{\frac{m}{2}}$ via:
$$
l_ 0 \rightarrow {\rm M}_{m/2} (k_{\Delta}) \times {\rm M}_{m/2} (k_{\Delta}), \, 
x \mapsto (i(x), i(x))
$$
Alors $\sqrt{\varepsilon}$ s'identifie à une matrice $\nu \in {\rm GL}_{m/2} (k_{\Delta})$ 
(et $\nu^2 \in i (l_0^{\times})$).\\
De même, on a une injection:
$$
i_0 : l \hookrightarrow {\rm M}_{m} (k_{\Delta}), \, 
x = x_1 + \sqrt{\varepsilon} x_2 \mapsto {\rm Mat}_{(1, \sqrt{\varepsilon})} (\varphi_x)
$$
où $\varphi_x$ est la multiplication à gauche par $x$. 
On en déduit que le tore $l^{\times}$ s'identifie à:
$$
\left\{ \left( \begin{array}{c|c}
X_1 & X_2 \nu^2 \\
\hline
X_2 & X_1 \\
\end{array} \right) : X_1, X_2 \in i(l_0) \right\}^{\times} 
$$
Posons:
$$
\eta = \left( \begin{array}{c|c}
0 & \nu^2 \\
\hline
I_{m/2} & 0 \\
\end{array} \right)
$$
On remarque que $\eta$ correspond à $\sqrt{\varepsilon}$ et que:
$$
\langle w \rangle L_{m/2} = \langle \eta \rangle L_{m/2} 
\Rightarrow \langle \varpi_{\mathbb{K}} w \rangle L_{m/2} = \langle \varpi_{\mathbb{K}} \eta \rangle L_{m/2} 
$$
De plus, on vérifie facilement que 
$\varpi_{\mathbb{K}} \eta$
est une uniformisante de $\mathbb{E}$.
On définit alors $\upsilon$ un caractère de $\langle \varpi_{\mathbb{K}} \eta \rangle L_{m/2}$ par:
$$
\forall x \in L_{m/2}, \forall s \in \mathbb{Z}, 
\upsilon ((\varpi_{\mathbb{K}} \eta)^s x) = \chi ((\varpi_{\mathbb{K}} \eta)^s) 
\chi_{\Lambda_m / \mathbb{E}} ((\varpi_{\mathbb{K}} \eta)^s)
$$
où $\chi_{\mathbb{K}_n / \mathbb{E}} = \chi_{\Lambda_m / \mathbb{E}}$ est le caractère quadratique de 
$\mathbb{K}_n / \mathbb{E}$.
Comme $\mathbb{K}_n / \mathbb{E}$ est non ramifiée, la norme ${\rm N}_{\mathbb{K}_n / \mathbb{E}}$ 
est surjective sur les unités.
Puisque $\varpi_{\mathbb{K}} \eta$ est 
une uniformisante de $\mathbb{E}$, on a:
$$
\upsilon (\varpi_{\mathbb{K}} \eta) = - \chi (\varpi_{\mathbb{K}}) \overline{\chi} (\eta)
$$
\end{nota}

\begin{lem}\label{LemmeHakimMurnaghan}
L'inclusion canonique:
$$
{\rm Hom}_{\langle \varpi_{\mathbb{K}} \eta \rangle L_{m/2}} (\pi_0, \upsilon) 
\subseteq {\rm Hom}_{L_{m/2}} (\overline{\gamma}_0, \mathds{1})
$$
est un isomorphisme de $\mathbb{C}$-espaces vectoriels.
\end{lem}

\begin{démo}
Pour la démonstration, on pourra se référer à \cite{HakimMurnaghan1}, proposition 6.3 page 1877. 

\end{démo}

\begin{theo}\label{CritereDistinctionTRdImpair}
La représentation $\pi$ est ${\rm GL}_m (\mathcal{D})$-distinguée si et seulement si 
$m$ est pair (i.e $n$ est pair), $\chi$ est trivial sur $\mathbb{F}^{\times}$, $\overline{\chi}$ est 
trivial sur $l_0^{\times}$ et 
$\chi (\varpi_{\mathbb{K}}) \overline{\chi} (\eta) = -1$.\\
De plus, si $\pi$ est ${\rm GL}_m (\mathcal{D})$-distinguée, on a:
$$
{\rm dim}_{\mathbb{C}} ({\rm Hom}_{{\rm GL}_m (\mathcal{D})} (\pi, \mathds{1})) = 1
$$
\end{theo}

\begin{démo} 
Supposons tout d'abord que $\pi$ est ${\rm GL}_m (\mathcal{D})$-distinguée, alors 
nécessairement $m$ est pair, $\chi$ est trivial sur $\mathbb{F}^{\times}$ 
et:
$$
{\rm Hom}_{\langle \varpi_{\mathbb{K}} w \rangle L_{m/2}} (\pi_0, \mathds{1}) \neq 0
$$
En particulier $\overline{\gamma}_0$ est $L_{m/2}$-distinguée 
donc $\overline{\chi}$ est trivial sur $l_0^{\times}$.
Rappelons que l'on a l'égalité 
$\langle \varpi_{\mathbb{K}} w \rangle L_{m/2}
= \langle \varpi_{\mathbb{K}} \eta \rangle L_{m/2}$.
Comme:
$$
{\rm Hom}_{\langle \varpi_{\mathbb{K}} \eta \rangle L_{m/2}} (\pi_0, \mathds{1})
 \subseteq {\rm Hom}_{L_{m/2}} (\overline{\gamma}_0, \mathds{1})
$$
et que l'espace ${\rm Hom}_{L_{m/2}} (\overline{\gamma}_0, \mathds{1})$ est de dimension $1$, on a:
$$
{\rm Hom}_{L_{m/2}} (\overline{\gamma}_0, \mathds{1}) 
= {\rm Hom}_{\langle \varpi_{\mathbb{K}} \eta \rangle L_{m/2}} (\pi_0, \mathds{1})
$$
donc d'après le lemme \ref{LemmeHakimMurnaghan}:
$$
{\rm Hom}_{\langle \varpi_{\mathbb{K}} \eta \rangle L_{m/2}} (\pi_0, \mathds{1}) 
= {\rm Hom}_{\langle \varpi_{\mathbb{K}} \eta \rangle L_{m/2}} (\pi_0, \upsilon)
$$
Fixons $\varphi$ forme linéaire non nulle dans 
${\rm Hom}_{\langle \varpi_{\mathbb{K}} \eta \rangle L_{m/2}} (\pi_0, \mathds{1})$ 
et notons $V$ l'espace de $\pi_0$.
Alors, pour tout $u$ dans $V$, 
$\varphi (\pi_0 (\varpi_{\mathbb{K}} \eta).u) = 
\varphi (u)$.
Or $\varphi (\pi_0 (\varpi_{\mathbb{K}} \eta).u) 
= \upsilon (\varpi_{\mathbb{K}} \eta) \varphi (u) 
= - \chi (\varpi_{\mathbb{K}}) \overline{\chi} (\eta) \varphi (u)$.
On en déduit que 
$\chi (\varpi_{\mathbb{K}}) \overline{\chi} (\eta) = -1$.
De plus, si $\pi$ est ${\rm GL}_m (\mathcal{D})$-distinguée, on a:
$$
{\rm dim} ({\rm Hom}_{{\rm GL}_m (\mathcal{D})} (\pi, \mathds{1})) 
= {\rm dim} ({\rm Hom}_{\langle \varpi_{\mathbb{K}} w \rangle L_{m/2}} (\pi_0, \mathds{1})) 
= {\rm dim} ({\rm Hom}_{L_{m/2}} (\overline{\gamma}_0, \mathds{1})) = 1
$$
Supposons à présent que $m$ est pair, 
 $\chi$ trivial sur $\mathbb{F}^{\times}$, $\overline{\chi}$ 
trivial sur $l_0^{\times}$ et 
$\chi (\varpi_{\mathbb{K}}) \overline{\chi} (\eta) = -1$. 
Comme $m$ est pair et $\overline{\chi}$  
trivial sur $l_0^{\times}$, d'après le lemme \ref{DistinctionParLeviCuspidale}:
$$
{\rm dim} ({\rm Hom}_{L_{m/2}} (\overline{\gamma}_0, \mathds{1})) = 1
$$
De plus, d'après le lemme \ref{LemmeHakimMurnaghan}, on a:
$$
{\rm Hom}_{\langle \varpi_{\mathbb{K}} \eta \rangle L_{m/2}} (\pi_0, \upsilon)
\simeq {\rm Hom}_{L_{m/2}} (\overline{\gamma}_0, \mathds{1})
$$
Puisque $\chi (\varpi_{\mathbb{K}}) \overline{\chi} (\eta) = -1$, $\upsilon$ est le caractère trivial, par suite:
$$
{\rm Hom}_{\langle \varpi_{\mathbb{K}} \eta \rangle L_{m/2}} (\pi_0, \mathds{1})
\simeq {\rm Hom}_{L_{m/2}} (\overline{\gamma}_0, \mathds{1}) \neq 0
$$
et donc $\pi$ est bien ${\rm GL}_m (\mathcal{D})$-distinguée.

\end{démo}

\section{Conditions de distinction lorsque $d$ est pair.}\label{PartieCNSDistinctionPair}

\subsection{Quelques propriétés.}

On suppose dans toute cette partie que $d$ est pair.

\begin{rmq}
Puisque $d$ est pair, on a $\delta = d/2$ et $\mu = 2m$.\\
De plus, étant donné que $2$ divise $d$, on peut supposer que 
$\mathbb{F} \subseteq \mathbb{K} \subseteq \mathcal{D}$.
\end{rmq}

Rappelons le théorème du bicommutant (on pourra se référer à \cite{DennisFarb}, Partie I paragraphe 3) :

\begin{theo}\label{bicommutant}
Soit $k$ un corps commutatif,
$S$ une $k$-algèbre centrale simple (de dimension finie) et 
$R$ une sous-algèbre simple de $S$.
On a les propriétés suivantes:
\begin{itemize}
\item[i)] $C(R)$, le commutant de $R$ dans $S$, est une algèbre simple.
\item[ii)] $[S : k] = [ R : k] [C (R) : k]$.
\item[iii)] $C ( C (R)) = R$, i.e $R$ est égal à son bicommutant.
\end{itemize}
\end{theo}

\begin{propr}
On peut identifier $\Delta$ au commutant de $\mathbb{K}$ dans $\mathcal{D}$:
$$
\Delta = \{ x \in \mathcal{D} : \forall k \in \mathbb{K}, xk = kx \}
$$
\end{propr}

\begin{démo}
Il suffit d'utiliser le théorème \ref{bicommutant} pour montrer que le commutant de $\mathbb{K}$ dans $\mathcal{D}$ 
est une $\mathbb{K}$-algèbre à division centrale d'indice $d/2$. 
De plus, si l'on note $\Delta^{'}$ le commutant de $\mathbb{K}$ dans $\mathcal{D}$, on a un 
isomorphisme de $\mathbb{K}$-algèbres:
$$
\mathcal{D} \otimes_{\mathbb{F}} \mathbb{K} \rightarrow {\rm End}_{\Delta^{'}} (\mathcal{D}), \, 
y \otimes k \mapsto (\mathcal{D} \rightarrow \mathcal{D}, x \mapsto yxk)
$$

\end{démo}

\begin{rmq}
On remarque que:
$$
[ \mathcal{D} : \mathbb{F}] = d^2 = [ \mathcal{D} : \Delta] [ \Delta : \mathbb{K}] [ \mathbb{K} : \mathbb{F}] 
= [ \mathcal{D} : \Delta] \times \left( \frac{d}{2}\right)^2 \times 2
$$
Donc $[ \mathcal{D} : \Delta] = 2$ et $\mathcal{D}$ peut être vu comme un $\Delta$-espace 
vectoriel à droite de dimension~$2$.
\end{rmq}

On montre facilement le résultat suivant :

\begin{propr}
Posons:
$$
\Phi : \mathcal{D} \otimes_{\mathbb{F}} \mathbb{K} \rightarrow {\rm End}_{\Delta} (\mathcal{D}), \, 
x \otimes k \mapsto [f_{x \otimes k} : \mathcal{D} \rightarrow \mathcal{D} , y \mapsto xyk]
$$
Alors $\Phi$ est un isomorphisme de $\mathbb{K}$-algèbres. 
\end{propr}

En appliquant le théorème de Skölem-Noether à $\sigma$, on montre la propriété suivante :

\begin{propr}\label{SkolemNoether}
Il existe $d_0$ dans $\mathcal{D}^{\times}$ tel que pour tout $x$ dans $\mathcal{D}$ et 
tout $k$ dans $\mathbb{K}$:
$$
\Phi (\sigma. (x \otimes k)) : y \mapsto f_{x \otimes k} (y d_0) d_0^{-1}
$$
De plus, $d_0$ est unique modulo la multiplication par un élément de $\Delta$ à droite.
\end{propr}

\subsection{Cas où l'extension $\mathbb{K} / \mathbb{F}$ est non ramifiée.}

On suppose dans toute cette partie que l'extension $\mathbb{K} / \mathbb{F}$ est non ramifiée.

\begin{propr}
Dans l'énoncé de la propriété \ref{SkolemNoether}, on peut choisir $d_0 = \varpi_{\mathcal{D}}$.
\end{propr}

\begin{démo}
On peut choisir $\mathbb{L}/\mathbb{F}$ une extension non ramifiée de degré $d$ contenue dans $\mathcal{D}$ 
qui soit normalisée par $\varpi_{\mathcal{D}}$. 
Puisque $2$ divise $d$ et que $\mathbb{K} / \mathbb{F}$ est non ramifiée, on peut supposer que 
$\mathbb{F} \subseteq \mathbb{K} \subseteq \mathbb{L}$. Soit:
$$
\tau : \mathbb{L} \rightarrow \mathbb{L}, x \mapsto \varpi_{\mathcal{D}} x \varpi_{\mathcal{D}}^{-1}
$$
Alors $\tau$ est un générateur du groupe de Galois ${\rm Gal} (\mathbb{L} / \mathbb{F})$.
On en déduit que $\tau_{\vert \mathbb{K}} = \sigma$ et que l'on peut choisir $d_0 = \varpi_{\mathcal{D}}$.

\end{démo}

\begin{rmq}
Comme $\mathbb{K} / \mathbb{F}$ est non ramifiée, on peut supposer que 
$\varpi_{\mathbb{F}} = \varpi_{\mathbb{K}}$. On a donc:
$$
\varpi_{\mathcal{D}}^d = \varpi_{\mathbb{F}} = \varpi_{\mathbb{K}} = \varpi_{\Delta}^{\frac{d}{2}}
$$ 
De plus, pour tout $x$ dans $\Delta^{\times}$, 
$v_{\mathcal{D}} (x) = 2 v_{\Delta} (x)$.
\end{rmq}

\begin{propr}\label{Injection}
On peut choisir $(1, \varpi_{\mathcal{D}})$ comme $\Delta$-base de $\mathcal{D}$ (à droite).
On a ainsi un isomorphisme de $\mathbb{K}$-algèbres:
$$
\Phi : \mathcal{D} \otimes_{\mathbb{F}} \mathbb{K} \rightarrow {\rm M}_2 (\Delta), \, 
x \otimes k \mapsto {\rm Mat}_{(1, \varpi_{\mathcal{D}})} (f_{x \otimes k})
$$
Via cet isomorphisme , $\varpi_{\mathcal{D}}$ s'identifie à la matrice:
$$
\Pi_{\Delta} = \left( \begin{array}{cc}
0 & \varpi_{\Delta}\\
1 & 0\\
\end{array} \right)
$$
De plus, on a $\mathcal{O}_{\mathcal{D}} = \mathcal{O}_{\Delta} \oplus \varpi_{\mathcal{D}} \mathcal{O}_{\Delta}$ 
et $\mathcal{P}_{\mathcal{D}} = \mathcal{P}_{\Delta} \oplus \varpi_{\mathcal{D}} \mathcal{O}_{\Delta}$.

\end{propr}

\begin{démo}
Il est clair que $(1, \varpi_{\mathcal{D}})$ est une $\Delta$-base de $\mathcal{D}$. 
On fixe $\mathbb{L} / \mathbb{F}$ une extension non ramifiée de degré $d$ 
contenue dans $\mathcal{D}$, alors $\mathbb{F} \subseteq \mathbb{K} \subseteq \mathbb{L}$. 
Ainsi $\mathbb{L} / \mathbb{K}$ une extension non ramifiée de degré $d/2$. On en déduit que:
$$
\mathcal{O}_{\Delta} = 1.\mathcal{O}_{\mathbb{L}} \oplus \varpi_{\Delta}.\mathcal{O}_{\mathbb{L}} \oplus 
\cdots \oplus \varpi_{\Delta}^{d/2 -1}.\mathcal{O}_{\mathbb{L}}
$$
et:
$$
\mathcal{O}_{\mathcal{D}} = 1.\mathcal{O}_{\mathbb{L}} \oplus \varpi_{\mathcal{D}}.\mathcal{O}_{\mathbb{L}} \oplus 
\cdots \oplus \varpi_{\mathcal{D}}^{d -1}.\mathcal{O}_{\mathbb{L}}
$$
Ainsi $\mathcal{O}_{\mathcal{D}} = \mathcal{O}_{\Delta} \oplus \varpi_{\mathcal{D}} \mathcal{O}_{\Delta}$ et  
$\mathcal{P}_{\mathcal{D}} = \mathcal{P}_{\Delta} \oplus \varpi_{\mathcal{D}} \mathcal{O}_{\Delta}$.
Il est immédiat que $\Phi (\varpi_{\mathcal{D}}) = \Pi_{\Delta}$.

\end{démo}

\begin{nota}\label{Immeuble}
En fixant $(1, \varpi_{\mathcal{D}}, \cdots, 1, \varpi_{\mathcal{D}})$ comme $\Delta$-base de 
$\mathcal{D}^m$, on peut identifier ${\rm End}_{\Delta} (\mathcal{D}^m)$ à ${\rm M}_{2m} (\Delta)$. 
On fixe pour la suite un isomorphisme de $\mathbb{K}$-algèbres:
$$
\widetilde{\Phi} : {\rm M}_m (\mathcal{D} \otimes_{\mathbb{F}} \mathbb{K}) \rightarrow {\rm M}_{2m} (\Delta) \simeq {\rm End}_{\Delta} (\mathcal{D}^m), \, 
[a_{i,j}] \mapsto [\Phi (a_{i,j})]
$$
(on en déduit une injection de ${\rm M}_m (\mathcal{D})$ dans ${\rm M}_{2m} (\Delta)$). On notera:

$$
\mathcal{L}_1 = (L_i^1)_{i \in \mathbb{Z}}, \, 
L_i^1 = \mathcal{P}_{\mathcal{D}}^{2i} \oplus \cdots \oplus  \mathcal{P}_{\mathcal{D}}^{2i} \subseteq \mathcal{D}^m
$$
et:
$$
\mathcal{L}_2 = (L_i^2)_{i \in \mathbb{Z}}, \, 
L_i^2 = \mathcal{P}_{\mathcal{D}}^{2i+1} \oplus \cdots \oplus  \mathcal{P}_{\mathcal{D}}^{2i+1}
$$
Pour tout $i \in \mathbb{Z}$, $\mathcal{P}_{\mathcal{D}}^{2i}$ est un $\mathcal{O}_{\Delta}$-réseau de $\mathcal{D}$.
Donc $L_i^1$ est un $\mathcal{O}_{\Delta}$-réseau de $\mathcal{D}^m$. 
De même, $L_i^2$ est un $\mathcal{O}_{\Delta}$-réseau de $\mathcal{D}^m$.
On vérifie facilement que $\mathcal{L}_1$ 
et $\mathcal{L}_2$ sont des chaînes de période $1$, ce sont donc des sommets de l'immeuble $X_{\mathbb{K}}$.\\
Pour $i$ dans $\lbrace 1, 2 \rbrace$, on note $\mathcal{A}_i$ les ordres héréditaires associés à $\mathcal{L}_i$, alors:
$$
\mathcal{A}_i = \lbrace a \in {\rm M}_{2m} (\Delta) : \forall k \in \mathbb{Z}, 
a L_k^i \subseteq L_k^i \rbrace
$$
et:
$$
\mathcal{K}_i = \mathcal{K} (\mathcal{L}_i) 
= \lbrace g \in {\rm GL}_{2m} (\Delta) : g \mathcal{A}_i g^{-1} = \mathcal{A}_i \rbrace
$$
D'après \cite{BushnellFröhlich} (théorème 1.3.2 page 217), 
$\mathcal{K}_i  
= \lbrace g \in {\rm GL}_{2m} (\Delta) : \exists n_g \in \mathbb{Z} : \forall k \in \mathbb{Z}, 
g.L_k^i = L_{k+n_g}^i \rbrace$ est un sous-groupe ouvert compact modulo le centre maximal.
\end{nota}

\begin{propr}
L'action de $\langle \sigma \rangle$ échange les deux sommets $\mathcal{L}_1$ et $\mathcal{L}_2$, i.e 
$\sigma (\mathcal{K}_1) = \mathcal{K}_2$.
\end{propr}

\begin{démo}
Soit $g \in \mathcal{K}_1$.
Soit $n_g \in \mathbb{Z}$ tel que pour tout $l \in \mathbb{Z}, g.L_l^1 = L_{l+n_g}^1$.
Identifions $ \mathcal{P}_{\mathcal{D}}^{2i} \oplus \cdots \oplus  \mathcal{P}_{\mathcal{D}}^{2i}$ 
avec $ \mathcal{P}_{\mathcal{D}}^{2i} \times \cdots \times  \mathcal{P}_{\mathcal{D}}^{2i}$.
Ainsi, pour tout $l \in \mathbb{Z}$, on a:
$$
g (\mathcal{P}_{\mathcal{D}}^{2l} \times \cdots \times  \mathcal{P}_{\mathcal{D}}^{2l}) 
= \mathcal{P}_{\mathcal{D}}^{2(l+n_g)} \times \cdots \times  \mathcal{P}_{\mathcal{D}}^{2(l+n_g)}
$$
Puisque $g \in {\rm End}_{\Delta} (\mathcal{D}^m)$, on peut supposer que:
$$
g = \widetilde{\Phi} [a_{i,j}], 
a_{i,j} = \sum_{q=1}^{r_{i,j}} d_{i,j}^{q} \otimes k_{i,j}^{q} 
$$
où $d_{i,j}^{q} \in \mathcal{D}$ et $k_{i,j}^{q} \in \mathbb{K}$.
Alors, pour tout $(x_1, \cdots, x_m) \in \mathcal{D}^m$, on a:
$$
g (x_1, \cdots, x_m) = [\Phi(a_{i,j})][x_1, \cdots, x_m] 
= (\sum_{j=1}^{m} \Phi(a_{1,j}) (x_j), \cdots, \sum_{j=1}^{m} \Phi(a_{m,j})(x_j)) 
$$
où, pour tout $i$ dans $\{ 1, \cdots, m \}$, on a 
$\Phi(a_{i,j}) (x_j) = \sum_{q=1}^{r_{i,j}} d_{i,j}^{q} x_j k_{i,j}^{q}$
et:
\begin{eqnarray*}
g^{\sigma} (x_1, \cdots, x_m)
& = & [\Phi^{\sigma} (a_{i,j})][x_1, \cdots, x_m] \\
& = & ((\sum_{j=1}^{m} \Phi(a_{1,j}) (x_j \varpi_{\mathcal{D}})) \varpi_{\mathcal{D}}^{-1}, \cdots, 
    (\sum_{j=1}^{m} \Phi(a_{m,j})(x_j \varpi_{\mathcal{D}})) \varpi_{\mathcal{D}}^{-1}) 
\end{eqnarray*}

Soit $(x_1, \cdots, x_m) \in 
\mathcal{P}_{\mathcal{D}}^{2i+1} \times \cdots \times  \mathcal{P}_{\mathcal{D}}^{2i+1}
 = L_{i}^2$.
Alors $(x_1 \varpi_{\mathcal{D}}, \cdots, x_m \varpi_{\mathcal{D}}) \in L_{i+1}^1$.
Puisque $g \in \mathcal{K}_1$, on a:
\begin{eqnarray*}
g (x_1 \varpi_{\mathcal{D}}, \cdots, x_m \varpi_{\mathcal{D}}) \in L_{i+1+n_g}^1
& \Rightarrow & \forall l \in \lbrace 1, \cdots, m \rbrace, 
                \sum_{j=1}^{m} \Phi(a_{l,j}) (x_j \varpi_{\mathcal{D}}) \in  \mathcal{P}_{\mathcal{D}}^{2(i+1+n_g)}\\
& \Rightarrow & \forall l \in \lbrace 1, \cdots, m \rbrace, 
                 (\sum_{j=1}^{m} \Phi(a_{l,j}) (x_j \varpi_{\mathcal{D}}))\varpi_{\mathcal{D}}^{-1} \in  \mathcal{P}_{\mathcal{D}}^{2(i+n_g)+1}\\
& \Rightarrow & g^{\sigma} (x_1, \cdots, x_m) \in L_{i+n_g}^2
\end{eqnarray*}
On en déduit que $g^{\sigma} \in \mathcal{K}_2$.\\
Par suite, $\sigma (\mathcal{K}_1) \subseteq \mathcal{K}_2$, et par maximalité, 
$\sigma (\mathcal{K}_1) = \mathcal{K}_2$.

\end{démo}

\subsubsection{Explicitation des injections d'immeubles.}

\begin{propr}\label{InjectionImmeubleNRPair}
Soit $A_{\mathbb{F}}$ (resp. $A_{\mathbb{K}}$) l'appartement standard de $X_{\mathbb{F}}$ (resp. $X_{\mathbb{K}}$) 
que l'on identifie à l'espace affine $\mathbb{R}^m / \mathbb{R} (1, \cdots, 1)$ 
(resp. $\mathbb{R}^{2m} / \mathbb{R} (1, \cdots, 1)$) 
et dont l'ensemble des sommets s'identifie à 
$\mathbb{Z}^m / \mathbb{Z} (1, \cdots, 1)$ (resp. $\mathbb{Z}^{2m} / \mathbb{Z} (1, \cdots, 1)$).
Soit $j : X_{\mathbb{F}} \rightarrow X_{\mathbb{K}}$ l'injection naturelle entre les immeubles.\\
Pour tout sommet $x= \overline{(x_1, \cdots, x_m)}$ de $X_{\mathbb{F}}$, où pour tout $i$ dans 
$\{ 1, \cdots, m \}$, $x_i \in \mathbb{Z}$, il existe un unique couple $(k_i, r_i)$ de $\mathbb{Z}^2$ 
tel que $x_i = 2 k_i + r_i$ et $r_i \in \{ 0, 1 \}$. Alors:
$$
j(x) = \overline{(y_1, y_2, \cdots, y_{2m-1}, y_{2m})}
$$
où, pour tout $i$ dans $\lbrace 1, \cdots, m \rbrace$, $(y_{2i-1}, y_{2i}) = (\frac{1}{2} + k_i, k_i)$
si $x_i$ est pair et, si $x_i$ est impair, $(y_{2i-1}, y_{2i}) = (1 + k_i, \frac{1}{2} + k_i)$.
\end{propr}

\begin{rmq}
Pour définir $j$, il nous faut trouver un sommet $s$ de 
$X_{\mathbb{F}}$ tel que $j(s)$ soit fixe par l'action de $\sigma$.
On remarque que, puisque $\mathcal{K}_1^{\sigma} = \mathcal{K}_2$, $\sigma$ échange 
les deux sommets $\mathcal{L}_1$ et $\mathcal{L}_2$, donc fixe leur milieu.
\end{rmq}

\begin{démo}
On fixe $V$ un $\mathcal{D}$-espace vectoriel (à droite) de dimension $m$ et 
$\mathcal{B} = (e_1, \cdots, e_m)$ une $\mathcal{D}$-base de $V$.
Alors $({\rm End}_{\mathcal{D}} (V))^{\times}$ s'identifie à ${\rm GL}_m (\mathcal{D})$.
On a vu que $(1, \varpi_{\mathcal{D}})$ est une $\Delta$-base (à droite) de $\mathcal{D}$.
Posons $\widetilde{\mathcal{B}} = (e_1.1, e_1.\varpi_{\mathcal{D}}, \cdots, e_m.1, e_m.\varpi_{\mathcal{D}})$,
ainsi $\widetilde{\mathcal{B}}$ est une $\Delta$-base de $V$ qui nous permet d'identifier 
$({\rm End}_{\Delta} (V))^{\times}$ à ${\rm GL}_{2m} (\Delta)$.
Soit $X_{\mathbb{F}}$ (resp. $X_{\mathbb{K}}$) l'immeuble de Bruhat-Tits de 
$({\rm End}_{\mathcal{D}} (V))^{\times} \simeq {\rm GL}_m (\mathcal{D})$
(resp. l'immeuble de $({\rm End}_{\Delta} (V))^{\times} \simeq {\rm GL}_{2m} (\Delta)$).
Soit $A_{\mathbb{F}}$ (resp. $A_{\mathbb{K}}$) l'appartement de $X_{\mathbb{F}}$ associé à la base 
$\mathcal{B}$ (resp. associé à la base $\widetilde{\mathcal{B}}$) que l'on identifie à l'espace affine 
$\mathbb{R}^m / \mathbb{R} (1, \cdots, 1)$ (resp. $\mathbb{R}^{2m} / \mathbb{R} (1, \cdots, 1)$) 
et dont l'ensemble des sommets s'identifie à 
$\mathbb{Z}^m / \mathbb{Z} (1, \cdots, 1)$ (resp. $\mathbb{Z}^{2m} / \mathbb{Z} (1, \cdots, 1)$).
On définit $s_0$, un sommet de $A_{\mathbb{F}}$ par:
$$
s_0 = (e_1.\mathcal{P}_{\mathcal{D}}^i \oplus \cdots \oplus e_m.\mathcal{P}_{\mathcal{D}}^i)_{i \in \mathbb{Z}} 
= \overline{(0, \cdots, 0)}
$$
On a déjà vu que l'action du groupe de Galois $\langle \sigma \rangle$ 
échange les stabilisateurs des deux sommets de $A_{\mathbb{K}}$, $\mathcal{L}_1 = 
(e_1.\mathcal{P}_{\mathcal{D}}^{2i} \oplus \cdots \oplus e_m.\mathcal{P}_{\mathcal{D}}^{2i})_{i \in \mathbb{Z}}$ et 
$\mathcal{L}_2 
= (e_1.\mathcal{P}_{\mathcal{D}}^{2i+1} \oplus \cdots \oplus e_m.\mathcal{P}_{\mathcal{D}}^{2i+1})_{i \in \mathbb{Z}}$.
On en déduit que $\sigma$ échange les deux sommets $\mathcal{L}_1$ et $\mathcal{L}_2$, et donc fixe leur milieu.
On définit $j(s_0)$ comme étant le milieu du segment $[ \mathcal{L}_1, \mathcal{L}_2]$.
Comme $\mathcal{O}_{\mathcal{D}} = \mathcal{O}_{\Delta} \oplus \varpi_{\mathcal{D}} \mathcal{O}_{\Delta}$, on a: 
$$
\mathcal{L}_1 = [e_1.\mathcal{O}_{\mathcal{D}} \oplus \cdots \oplus e_m.\mathcal{O}_{\mathcal{D}}] 
= [e_1.\mathcal{O}_{\Delta} \oplus e_1 \varpi_{\mathcal{D}}.\mathcal{O}_{\Delta}  \oplus \cdots \oplus 
e_m.\mathcal{O}_{\Delta} \oplus e_m \varpi_{\mathcal{D}}.\mathcal{O}_{\Delta}]
$$
Donc $\mathcal{L}_1 = \overline{(\underbrace{0, \cdots, 0}_{2m})}$. 
De même, puisque $\mathcal{P}_{\mathcal{D}} = \mathcal{P}_{\Delta} \oplus \varpi_{\mathcal{D}} \mathcal{O}_{\Delta}$, 
on a:
$$
\mathcal{L}_2 = [e_1.\mathcal{P}_{\mathcal{D}} \oplus \cdots \oplus e_m.\mathcal{P}_{\mathcal{D}}] 
= [e_1.\mathcal{P}_{\Delta} \oplus e_1 \varpi_{\mathcal{D}}.\mathcal{O}_{\Delta}  \oplus \cdots \oplus 
e_m.\mathcal{P}_{\Delta} \oplus e_m \varpi_{\mathcal{D}}.\mathcal{O}_{\Delta}]
$$
Et $\mathcal{L}_2 = \overline{(\underbrace{1, 0, \cdots, 1, 0}_{2m})}$.
Enfin:
\begin{eqnarray*}
j(s_0) 
& = & {\rm milieu} ([\mathcal{L}_1, \mathcal{L}_2]) 
       = \frac{1}{2} \mathcal{L}_1  
       + \frac{1}{2} \mathcal{L}_2 \\
& = & \overline{(\underbrace{\frac{1}{2}, 0, \frac{1}{2}, 0, \cdots, \frac{1}{2}, 0}_{2m})}
\end{eqnarray*}
D'après la remarque précédente, on a $\sigma (j(s_0)) = j(s_0)$, 
de plus:
$$
{\rm Stab}_{{\rm GL}_m (\mathcal{D})} (s_0) \simeq \langle \varpi_{\mathcal{D}} \rangle {\rm GL}_m (\mathcal{O}_{\mathcal{D}})
$$
Or ${\rm GL}_m (\mathcal{O}_{\mathcal{D}})$ stabilise les chaînes ($\mathcal{O}_{\mathcal{D}}$-chaînes de période $2$) 
$\mathcal{L}_1$ et $\mathcal{L}_2$ et $\varpi_{\mathcal{D}}$ échange ces deux chaînes.
On en déduit que ${\rm Stab}_{{\rm GL}_m (\mathcal{D})} (s_0)$ stabilise bien le milieu de $[ \mathcal{L}_1, \mathcal{L}_2]$, 
d'où l'inclusion ${\rm Stab}_{{\rm GL}_m (\mathcal{D})} (s_0) 
\subseteq {\rm Stab}_{{\rm GL}_{2m} (\Delta)} (j(s_0))$.\\
On va définir $j$ sur les sommets de $A_{\mathbb{F}}$.
On fixe $x = \overline{(x_1, \cdots, x_m)}$ un sommet de $A_{\mathbb{F}}$ ($x_i \in \mathbb{Z}$).
Alors $x = g_x . s_0$ où $g_x = {\rm diag} (\varpi_{\mathcal{D}}^{x_1}, \cdots, \varpi_{\mathcal{D}}^{x_m}) 
\in {\rm GL}_m (\mathcal{D})$.
On pose $j(x) = \widetilde{\Phi} (g_x). j(s_0)$. 
On a:
$$
\widetilde{\Phi} (g_x) = \left( \begin{array}{cccc}
\Pi_{\Delta}^{x_1}  & 0 & \cdots & 0\\
0 & \ddots &  & \vdots\\
\vdots & & \ddots & 0\\
0 & & 0 & \Pi_{\Delta}^{x_m}
\end{array} \right)
$$
(où $\Pi_{\Delta}$ a été définie en \ref{Injection}). 
Comme $x_i \in \mathbb{Z}$, on fixe $k_i, r_i \in \mathbb{Z}$ tels que 
$x_i = 2 k_i + r_i, \, r_i \in \lbrace 0, 1 \rbrace$.\\
Alors $\Pi_{\Delta}^{x_i} = \Pi_{\Delta}^{r_i} (\Pi_{\Delta}^2)^{k_i }
 = \Pi_{\Delta}^{r_i} \varpi_{\Delta}^{k_i } I_{2m}$.
On remarque que, pour $\alpha, \beta \in \mathbb{Z}$:
$$
\Pi_{\Delta} \left( \begin{array}{c}
\mathcal{P}_{\Delta}^{\alpha}\\
\mathcal{P}_{\Delta}^{\beta}
\end{array} \right) 
= \left( \begin{array}{c}
\mathcal{P}_{\Delta}^{\beta + 1}\\
\mathcal{P}_{\Delta}^{\alpha}
\end{array} \right)
$$
On a alors $j(x) = \overline{(y_1, y_2, \cdots, y_{2m-1}, y_{2m})}$
où, pour tout $i$ dans $\lbrace 1, \cdots, m \rbrace$:
$$
(y_{2i-1}, y_{2i}) = \left(\frac{1}{2} + k_i, 0+ k_i \right) = \left( \frac{1}{2} + k_i, k_i \right)
$$
si $x_i$ est pair ($r_i = 0$) et:
$$
(y_{2i-1}, y_{2i}) = \left( 0 + k_i + 1, \frac{1}{2} + k_i \right) 
= \left( k_i +1, \frac{1}{2} + k_i \right)
$$
si $x_i$ est impair ($r_i = 1$).
On a donc bien défini $j$ sur les sommets de $A_{\mathbb{F}}$. On définit $j$ sur l'appartement tout entier par affinité.\\
Soit $T$ le tore maximal déployé associé à l'appartement $A_{\mathbb{F}}$, 
et $N(T)$ son normalisateur.
Alors, avec le même raisonnement que dans la 
démonstration de la propriété \ref{InjectionImmeubleNRImpair}, on vérifie que 
$j$ est $N(T)$-équivariante sur les sommets en remarquant que si 
$x = \overline{(x_1, \cdots, x_m)}$ est un sommet de $A_{\mathbb{F}}$ et si $g \in N(T)$ est de la forme 
$g = t P_{\tau}$ avec $t = {\rm diag} (\varpi_{\mathcal{D}}^{\alpha_1}, \cdots, \varpi_{\mathcal{D}}^{\alpha_m})$, 
$\tau$ une permutation de $\{ 1, \cdots, m \}$ et $P_{\tau}$ la matrice de la permutation $\tau$, alors 
$j(t P_{\tau} .x) = \widetilde{\Phi} (t) \widetilde{\Phi} (P_{\tau}).j(x)$. 
Les calculs montrent que $\widetilde{\Phi} (P_{\tau}) = P_{\widetilde{\tau}}$ où 
$\widetilde{\tau}$ est une permutation de $\{ 1, \cdots, 2m-1, 2m \}$ définie par:
$$
\forall k \in \{ 1, \cdots, m \} , 
\widetilde{\tau} (2k-1) = 2 \tau (k) -1, 
\widetilde{\tau} (2k) = 2 \tau (k)
$$
On en déduit que $j$ est $N(T)$-équivariante sur les sommets, et donc sur l'appartement $A_{\mathbb{F}}$. 
Par conséquent, $j$ est bien l'injection naturelle entre les deux immeubles.

\end{démo}

\begin{propo}
Notons $C_0$ la chambre standard de $X_{\mathbb{F}}$.
Alors il n'y a aucun sommet de $X_{\mathbb{K}}$ dans l'image de la chambre standard, 
$j(C_0)$.
\end{propo}

\begin{démo}
Notons $\{s_0, s_1, \cdots, s_{m-1} \}$
les sommets de la chambre standard $C_0$, c'est-à-dire que pour tout $i$ dans $\{0, 1, \cdots, m-1 \}$, 
on a:
$$
s_i = [\mathcal{O}_{\mathcal{D}} \oplus \cdots \oplus \mathcal{O}_{\mathcal{D}} 
\oplus \underbrace{\mathcal{P}_{\mathcal{D}} \oplus \cdots \oplus \mathcal{P}_{\mathcal{D}}}_i]
$$
Or $s_i = h_i. s_0$ où 
$h_i = {\rm diag} (1, \cdots, 1, \underbrace{\varpi_{\mathcal{D}}, \cdots, \varpi_{\mathcal{D}}}_i)$.
En utilisant les mêmes notations que dans la démonstration précédente, on a $\Phi (1) = I_2$, 
et $\Phi (\varpi_{\mathcal{D}}) = \Pi_{\Delta}$.
La ${\rm GL}_m (\mathcal{D})$-équivariance de $j$ nous permet de vérifier que 
$j(s_i) = \widetilde{\Phi} (h_i).j(s_0)$
où 
$\widetilde{\Phi} (h_i)$ est la matrice diagonale par blocs 
$\widetilde{\Phi} (h_i) = {\rm diag} (I_2, \cdots, I_2, \underbrace{\Pi_{\Delta}, \cdots, \Pi_{\Delta}}_i)$.
On en déduit que:
$$
j(s_i) = \overline{(\underbrace{\frac{1}{2}, 0, \cdots, \frac{1}{2}, 0}_{2(m-i)}, 
\underbrace{1, \frac{1}{2}, \cdots, 1, \frac{1}{2}}_{2i})}
$$
Soit $t \in j(C_0)$.
Il existe $(\lambda_0, \lambda_1, \cdots, \lambda_{m-1}) \in (\mathbb{R}^{+})^{m}$ tels que:
$$
t = \sum_{i=0}^{m-1} \lambda_i j(s_i) \, \, \text{et} \, \, 
\sum_{i=0}^{m-1} \lambda_i = 1
$$
Ainsi, on peut calculer les deux premières coordonnées de $t$:
$$
t = \overline{((\sum_{i=0}^{m-1} \lambda_i) \times \frac{1}{2}, 0, \cdots)} 
= \overline{(\frac{1}{2}, 0, \cdots)} 
$$
Par conséquent, $t$ ne peut pas être à coordonnées entières dans 
$\mathbb{R}^{2m} / \mathbb{R} (1, \cdots, 1)$, donc ne peut en aucun cas être un sommet 
de $X_{\mathbb{K}}$.

\end{démo}

\subsubsection{Conditions de distinction.}

En raisonnant comme dans la démonstration de \ref{RaisonnementDistinction}, on montre le théorème suivant :

\begin{theo}\label{CritereDistinctionNRdPair}
La représentation $\pi$ n'est pas ${\rm GL}_m (\mathcal{D})$-distinguée.
\end{theo}

\begin{rmq}
Il est intéressant de noter que ce résultat de non distinction découle d'arguments purement géométriques.
\end{rmq}

\subsection{Cas où l'extension $\mathbb{K} / \mathbb{F}$ est totalement ramifiée, modérément ramifiée.}

On suppose dans toute cette partie que l'extension $\mathbb{K} / \mathbb{F}$ est totalement ramifiée, modérément ramifiée.
 On fixe les uniformisantes telles que 
 $\varpi_{\mathbb{K}}^2 = \varpi_{\mathbb{F}}$.

\begin{rmq}
 On a $\varpi_{\mathcal{D}}^d = \varpi_{\mathbb{F}} = \varpi_{\mathbb{K}}^2 = \varpi_{\Delta}^d$. 
 Ainsi, pour tout $x$ dans $\Delta^{\times}$, on a 
 \mbox{$v_{\mathcal{D}} (x) = v_{\Delta} (x)$}.
 Comme dans la partie précédente, on déduit de $\Phi$ un isomorphisme de 
 $\mathbb{K}$-algèbres, $\widetilde{\Phi}$, entre ${\rm M}_m (\mathcal{D} \otimes_{\mathbb{F}} \mathbb{K})$ 
 et ${\rm End}_{\Delta} (\mathcal{D}^m)$.
\end{rmq}

\begin{propr}
Soit $\mathcal{L} = (\mathcal{P}_{\mathcal{D}}^k \times \cdots \times \mathcal{P}_{\mathcal{D}}^k)_{k \in \mathbb{Z}}$. 
Alors $\mathcal{L}$ est une $\mathcal{O}_{\Delta}$-chaîne de réseaux de $\mathcal{D}^m$ de période $1$, 
donc s'identifie à un sommet de l'immeuble de Bruhat-Tits $X_{\mathbb{K}}$.\\
Notons, comme dans les notations \ref{NotationsImmeuble}, $\mathcal{K} = \mathcal{K} (\mathcal{L})$ 
(sous-groupe ouvert compact modulo 
le centre maximal), alors $\sigma$ stabilise $\mathcal{K}$ donc fixe le sommet $\mathcal{L}$ dans l'immeuble 
de ${\rm GL}_{2m} (\Delta)$.
\end{propr}

\begin{démo}
Puisque $\mathcal{O}_{\mathcal{D}}$ est un $\mathcal{O}_{\Delta}$-réseau 
de $\mathcal{D}$, $\mathcal{L}$ est une chaîne de période $1$. 
Vérifions que $\mathcal{K}$ est stable sous l'action du groupe de Galois $\langle \sigma \rangle$. 
Soient $g \in \mathcal{K}$ et $n_g \in \mathbb{Z}$ tels que pour tout $l$ dans $\mathbb{Z}$ :
$$ 
g(\mathcal{P}_{\mathcal{D}}^l \times \cdots \times \mathcal{P}_{\mathcal{D}}^l) 
= \mathcal{P}_{\mathcal{D}}^{l+n_g} \times \cdots \times \mathcal{P}_{\mathcal{D}}^{l+n_g}
$$
On peut supposer que $g = \widetilde{\Phi} ([a_{i,j}])$
où $a_{i,j} = \sum_{q=1}^{r_{i,j}} d_{i,j}^q \otimes k_{i,j}^{q}$ avec 
$d_{i,j}^q \in \mathcal{D}$ et $k_{i,j}^q \in \mathbb{K}$.
Ainsi, pour tout $(x_1, \cdots, x_m) \in \mathcal{D}^m$, on a:
$$
g(x_1, \cdots, x_m) 
= (\sum_{j=1}^{m} \Phi (a_{1,j})(x_j), \cdots, \sum_{j=1}^{m} \Phi (a_{m,j})(x_j))
$$
avec $\Phi (a_{i,j})(x_j) = \sum_{q=1}^{r_{i,j}} d_{i,j}^q x_j k_{i,j}^{q}$\\
et: 
$$
g^{\sigma} (x_1, \cdots, x_m) 
= ((\sum_{j=1}^{m} \Phi (a_{1,j})(x_j d_0))d_0^{-1}, \cdots, 
(\sum_{j=1}^{m} \Phi (a_{m,j})(x_j d_0))d_0^{-1})
$$
Soient $l \in \mathbb{Z}$ et 
$(x_1, \cdots, x_m) \in \mathcal{P}_{\mathcal{D}}^l \times \cdots \times \mathcal{P}_{\mathcal{D}}^l$.
Alors:
\begin{eqnarray*}
\forall j \in \{ 1, \cdots, m \}, x_j d_0 \in \mathcal{P}_{\mathcal{D}}^{l+ v_{\mathcal{D}} (d_0)}
& \Rightarrow & \forall i \in \{ 1, \cdots, m \}, \sum_{j=1}^{m} \Phi (a_{i,j})(x_j d_0) 
              \in \mathcal{P}_{\mathcal{D}}^{l+ v_{\mathcal{D}} (d_0) + n_g}\\
& \Rightarrow & \forall i \in \{ 1, \cdots, m \}, (\sum_{j=1}^{m} \Phi (a_{i,j})(x_j d_0)) d_0^{-1}
              \in \mathcal{P}_{\mathcal{D}}^{l + n_g}\\
& \Rightarrow & g^{\sigma} (x_1, \cdots, x_m) \in 
              \mathcal{P}_{\mathcal{D}}^{l + n_g} \times \cdots \times \mathcal{P}_{\mathcal{D}}^{l + n_g}
\end{eqnarray*}
On en déduit que $g^{\sigma} \in \mathcal{K}$, d'où $\sigma (\mathcal{K}) \subseteq \mathcal{K}$, et par maximalité, 
$\sigma (\mathcal{K}) = \mathcal{K}$.

\end{démo}

\begin{propr}
On peut fixer $(1, \zeta)$ une $\Delta$-base de $\mathcal{D}$ telle que 
$\mathcal{O}_{\mathcal{D}} = 1.\mathcal{O}_{\Delta} + \zeta. \mathcal{O}_{\Delta}$, ce qui nous permet d'identifier 
${\rm End}_{\Delta} (\mathcal{D})$ (resp. ${\rm End}_{\Delta} (\mathcal{D}^m)$) à ${\rm M}_2 (\Delta)$ 
(resp.  ${\rm M}_{2m} (\Delta)$). L'uniformisante $\varpi_{\mathcal{D}}$ s'identifie alors à une matrice 
$\Pi_{\Delta}$ qui vérifie, pour tout $i$ dans $\mathbb{Z}$ :
$$
\Pi_{\Delta} \left( \begin{array}{c}
\mathcal{P}_{\Delta}^i\\
\mathcal{P}_{\Delta}^i\\
\end{array} \right) 
= \left( \begin{array}{c}
\mathcal{P}_{\Delta}^{i+1}\\
\mathcal{P}_{\Delta}^{i+1}\\
\end{array} \right)
$$
\end{propr}

\begin{démo}
On définit $\mathbb{L} / \mathbb{F}$ une extension non ramifiée de degré $d$ contenue dans $\mathcal{D}$ 
normalisée par $\varpi_{\mathcal{D}}$ 
et $\tau : \mathbb{L} \rightarrow \mathbb{L}$ 
la conjugaison par $\varpi_{\mathcal{D}}$. On définit:
\begin{eqnarray*}
\mathbb{L}_0 = \mathbb{L}^{\tau^{\frac{d}{2}}} 
& = & \{ l \in \mathbb{L} : \tau^{\frac{d}{2}} (l) = \varpi_{\mathcal{D}}^{\frac{d}{2}} l \varpi_{\mathcal{D}}^{-\frac{d}{2}} = l \}\\
& = & \{ l \in \mathbb{L} : \varpi_{\mathbb{K}} l \varpi_{\mathbb{K}}^{-1} = l \}
\end{eqnarray*}
Alors $\mathbb{L} / \mathbb{L}_0$ est une extension quadratique non ramifiée. 
Comme $\Delta$ est le commutant de $\mathbb{K}$ dans $\mathcal{D}$, 
que $\mathbb{K} = \mathbb{F} [\varpi_{\mathbb{K}}]$ et que $\mathbb{F}$ 
est le centre de $\mathcal{D}$, $\Delta$ est l'ensemble des éléments de $\mathcal{D} = \langle \mathbb{L}, \varpi_{\mathcal{D}} \rangle$ 
($\mathbb{F}$-algèbre engendrée) qui commutent avec $\varpi_{\mathbb{K}}$.
Il est clair que $\varpi_{\mathcal{D}}$ 
commute avec $\varpi_{\mathbb{K}}$.
Si $l \in \mathbb{L}$, alors $l$ commute avec $\varpi_{\mathbb{K}}$ si et seulement si $l \in \mathbb{L}_0$.
Donc $\Delta = \langle \mathbb{L}_0, \varpi_{\mathcal{D}} \rangle$ ($\mathbb{K}$-algèbre engendrée).
Nous allons fixer $\zeta$ un générateur de $\mathbb{L} / \mathbb{L}_0$.
Puisque l'extension $\mathbb{L} / \mathbb{L}_0$ est quadratique non ramifiée, on peut fixer 
$a \in k_{\mathbb{L}} \backslash k_{\mathbb{L}_0}$ tel que $a^2 \in k_{\mathbb{L}_0}$ et fixer 
$\zeta \in \mathcal{O}_{\mathbb{L}}^{\times}$ tel que $\overline{\zeta} = a$.
Ainsi $\zeta^2 = u \in \mathcal{O}_{\mathbb{L}_0}^{\times}$, $k_{\mathbb{L}} = k_{\mathbb{L}_0} [a]$, 
et $\mathbb{L} = \mathbb{L}_0 [\zeta]$. 
Il est clair que $(1, \zeta)$ est une $\Delta$-base de $\mathcal{D}$ et 
$\mathcal{D} = 1.\Delta  + \zeta .\Delta$.
De plus, on vérifie facilement que $\mathcal{O}_{\mathbb{L}} = \mathcal{O}_{\mathbb{L}_0} [\zeta]$. Puisque 
$\mathcal{O}_{\mathcal{D}}$ est le $\mathcal{O}_{\mathbb{L}}$-module engendré par 
$\lbrace 1, \cdots, \varpi_{\mathcal{D}}^{d-1} \rbrace$, on a 
$\mathcal{O}_{\mathcal{D}} = 1. \mathcal{O}_{\Delta} + \zeta . \mathcal{O}_{\Delta}$.
Comme dans \ref{Injection} et \ref{Immeuble}, on définit les isomorphismes de $\mathbb{K}$-algèbre 
$\Phi$ et $\widetilde{\Phi}$ qui permettent d'identifier 
${\rm End}_{\Delta} (\mathcal{D})$ (resp. ${\rm End}_{\Delta} (\mathcal{D}^m)$) à ${\rm M}_2 (\Delta)$ 
(resp.  ${\rm M}_{2m} (\Delta)$).
On remarque que $\varpi_{\mathcal{D}} = \varpi_{\Delta} u$, 
avec $u \in \mathcal{O}_{\mathcal{D}}^{\times}$. 
On a $f_{\varpi_{\Delta} \otimes 1} (1) = \varpi_{\Delta} \in \Delta$
et:
$$
f_{\varpi_{\Delta} \otimes 1} (\zeta) = \varpi_{\Delta} \zeta 
= (\varpi_{\Delta} \zeta \varpi_{\Delta}^{-1}) \varpi_{\Delta} 
$$
Puisque $\zeta \in \mathcal{O}_{\mathbb{L}}^{\times}$ et que $\varpi_{\Delta}$ 
normalise $\mathbb{L}$, on a:
$$
\varpi_{\Delta} \zeta \varpi_{\Delta}^{-1} \in \mathcal{O}_{\mathbb{L}} = \mathcal{O}_{\mathbb{L}_0} [\zeta]
$$
Fixons $l_0, l_1 \in \mathcal{O}_{\mathbb{L}_0}$ tels que 
$\varpi_{\Delta} \zeta \varpi_{\Delta}^{-1} = l_0 + \zeta l_1$.
Puisque $\varpi_{\mathbb{F}}$ est une uniformisante de $\mathbb{L}_0$ et que 
$\varpi_{\Delta}^d = \varpi_{\mathbb{F}}$, il existe $m_0, m_1 \in \mathbb{N}$ tels que:
$$
v_{\Delta} (l_i) = d \times m_i, \, i \in \lbrace 1, 2 \rbrace
$$
On a donc 
$f_{\varpi_{\Delta} \otimes 1} (\zeta) = l_0 \varpi_{\Delta} + \zeta (l_1 \varpi_{\Delta})$ 
et:
$$
\varpi_{\Delta} \simeq \widetilde{\Pi}_{\Delta} = {\rm Mat}_{(1, \zeta)} (f_{\varpi_{\Delta} \otimes 1}) 
= \left( \begin{array}{cc}
         \varpi_{\Delta} & l_0 \varpi_{\Delta} \\
         0  &  l_1 \varpi_{\Delta} \\
         \end{array} \right)
$$
Il est clair que $\widetilde{\Pi}_{\Delta}$ vérifie la propriété annoncée. 
De plus, puisque $\varpi_{\mathcal{D}} = \varpi_{\Delta} u$, avec 
$u$ dans $\mathcal{O}_{\mathcal{D}}^{\times}$, l'uniformisante 
$\varpi_{\mathcal{D}}$ s'identifie à une matrice $\Pi_{\Delta}$, 
qui vérifie la même propriété que $\widetilde{\Pi}_{\Delta}$.

\end{démo}

\subsubsection{Explicitation des injections d'immeubles.}

\begin{propr}\label{InjectionImmeubleTRPair}
Soit $A_{\mathbb{F}}$ (resp. $A_{\mathbb{K}}$) l'appartement standard de $X_{\mathbb{F}}$ (resp. $X_{\mathbb{K}}$) 
que l'on identifie à l'espace affine 
$\mathbb{R}^m / \mathbb{R} (1, \cdots, 1)$ (resp. $\mathbb{R}^{2m} / \mathbb{R} (1, \cdots, 1)$) 
et dont l'ensemble des sommets s'identifie à 
$\mathbb{Z}^m / \mathbb{Z} (1, \cdots, 1)$ (resp. $\mathbb{Z}^{2m} / \mathbb{Z} (1, \cdots, 1)$).
L'injection naturelle entre les immeubles est donnée par:
$$
j : A_{\mathbb{F}} \rightarrow A_{\mathbb{K}} , \, 
\overline{(x_1, \cdots, x_m)} \mapsto \overline{(x_1,x_1, x_2, x_2, \cdots, x_m, x_m)}
$$
\end{propr}

\begin{rmq}
Plus précisément, on fixe $V$ un $\mathcal{D}$-espace vectoriel de dimension $m$ et 
une $\mathcal{D}$-base de $V$, $\mathcal{B} = (e_1, \cdots, e_m)$.
Alors $({\rm End}_{\mathcal{D}} (V))^{\times}$ s'identifie à ${\rm GL}_m (\mathcal{D})$.
On a vu que $(1, \zeta)$ est une $\Delta$-base (à droite) de $\mathcal{D}$.
Posons alors:
$$
\widetilde{\mathcal{B}} = (e_1.1, e_1.\zeta, \cdots, e_m.1, e_m.\zeta)
$$
Alors $\widetilde{\mathcal{B}}$ est une $\Delta$-base de $V$ qui nous permet d'identifier 
$({\rm End}_{\Delta} (V))^{\times}$ à ${\rm GL}_{2m} (\Delta)$. 
Les ensembles $X_{\mathbb{F}}$ (resp. $X_{\mathbb{K}}$) sont en fait les immeubles de Bruhat-Tits de 
$({\rm End}_{\mathcal{D}} (V))^{\times} \simeq {\rm GL}_m (\mathcal{D})$
(resp. $({\rm End}_{\Delta} (V))^{\times} \simeq {\rm GL}_{2m} (\Delta)$) 
et $A_{\mathbb{F}}$ (resp. $A_{\mathbb{K}}$) l'appartement de $X_{\mathbb{F}}$ associé à la base 
$\mathcal{B}$ (resp. associé à la base $\widetilde{\mathcal{B}}$).
\end{rmq}

\begin{démo}
On définit $j$ sur $A_{\mathbb{F}}$ par:
$$ 
j(\overline{(x_1, \cdots, x_m)}) = \overline{(x_1,x_1, x_2, x_2, \cdots, x_m, x_m)}
$$
pour tout $\overline{(x_1, \cdots, x_m)}$ dans $A_{\mathbb{F}}$. 
Soit $T$ le tore maximal déployé associé à l'appartement $A_{\mathbb{F}}$. 
Alors, $T$ s'identifie au tore diagonal dans 
${\rm GL}_m (\mathcal{D})$, 
$T = \{ {\rm diag} (t_1, \cdots, t_m) : t_i \in \mathbb{F}^{\times} \}$.
Comme dans la démonstration de la propriété \ref{InjectionImmeubleNRPair},
on note $N(T)$ son normalisateur, alors $N(T) = T_0 \mathcal{S}_m$.
Rappelons qu'en fixant $(1, \zeta)$ comme $\Delta$-base de $\mathcal{D}$, on a une injection
$\mathcal{D} \subseteq {\rm M}_2 (\Delta)$
qui nous permet d'identifier $\varpi_{\mathcal{D}}$ à une matrice $\Pi_{\Delta} \in {\rm M}_2 (\Delta)$.
De plus, pour tout $i \in \mathbb{Z}$:
$$
\mathcal{P}_{\mathcal{D}}^i \simeq \mathcal{P}_{\Delta}^i \oplus \mathcal{P}_{\Delta}^i \, \, \text{et} \, \, 
\Pi_{\Delta} \left( \begin{array}{c}
\mathcal{P}_{\Delta}^i\\
\mathcal{P}_{\Delta}^i\\
\end{array} \right) 
= \left( \begin{array}{c}
\mathcal{P}_{\Delta}^{i+1}\\
\mathcal{P}_{\Delta}^{i+1}\\
\end{array} \right)
$$
Le sous-groupe ouvert compact maximal ${\rm GL}_{2m} (\mathcal{O}_{\Delta})$ 
est stable sous l'action du groupe de Galois 
$\langle \sigma \rangle = {\rm Gal} (\mathbb{K} / \mathbb{F})$.
Fixons $s_0 = [\mathcal{O}_{\mathcal{D}} \oplus \cdots \oplus \mathcal{O}_{\mathcal{D}} ] 
= \overline{(\underbrace{0, \cdots, 0}_m)}$ un sommet de $A_{\mathbb{F}}$.
Alors \mbox{$j(s_0) = [\mathcal{O}_{\Delta} \oplus \cdots \oplus \mathcal{O}_{\Delta} ]
 = \overline{(\underbrace{0, \cdots, 0}_{2m})}$} et ${\rm Stab}_{{\rm GL}_{2m} (\Delta)} (j(s_0)) 
= \langle \varpi_{\Delta} \rangle {\rm GL}_{2m} (\mathcal{O}_{\Delta})$
est stable sous l'action de $\langle \sigma \rangle$.
On en déduit que $j(s_0)$ est un point de l'immeuble fixé par l'action du groupe de Galois.
De plus:
$$
{\rm Stab}_{{\rm GL}_{m} (\mathcal{D})} (s_0) 
= \langle \varpi_{\mathcal{D}} \rangle {\rm GL}_{m} (\mathcal{O}_{\mathcal{D}}) 
\subseteq {\rm Stab}_{{\rm GL}_{2m} (\Delta)} (j(s_0)) 
= \langle \varpi_{\Delta} \rangle {\rm GL}_{2m} (\mathcal{O}_{\Delta})
$$
Il nous reste à vérifier que $j$ est bien $N(T)$-équivariante. 
Si $x = \overline{(x_1, \cdots, x_m)} \in A_{\mathbb{F}}$ 
et $g \in N(T)$, alors 
on peut supposer que $g$ s'écrit sous la forme $g = t P_{\tau}$ où $\tau$ est 
une permutation de $\{ 1, \cdots, m \}$, $P_{\tau}$ la matrice de la permutation $\tau$ et 
$t$ est une matrice diagonale de la forme 
${\rm diag} (\varpi_{\mathcal{D}}^{\delta_1}, \cdots, \varpi_{\mathcal{D}}^{\delta_m})$ où 
chaque $\delta_i$ appartient à $\mathbb{Z}$.
On vérifie alors que $j(t P_{\tau}.x) = \widetilde{\Phi} (t) \widetilde{\Phi} (P_{\tau}).j(x)$  
en remarquant que $\widetilde{\Phi} (P_{\tau}) = P_{\widetilde{\tau}}$ où 
$\widetilde{\tau}$ est une permutation de $\{ 1, \cdots, 2m-1, 2m \}$ définie par:
$$
\widetilde{\tau} (2k-1) = 2 \tau (k) -1, 
\widetilde{\tau} (2k) = 2 \tau (k)
$$
pour tout $k$ dans $\{ 1, \cdots, m \}$.

\end{démo}

\begin{nota}
On notera $C_0$ et $\widetilde{C}_0$ les chambres standards de $X_{\mathbb{F}}$ et $X_{\mathbb{K}}$ 
respectivement.
Alors $C_0$ est un simplexe de dimension maximale dont les sommets sont 
$\{ s_0, s_1, \cdots, s_{m-1} \}$ tels que pour tout $i$ dans $\lbrace 0, 1, \cdots, m-1 \rbrace$, on a:
$$
s_i = [\mathcal{O}_{\mathcal{D}} \oplus \cdots \oplus \mathcal{O}_{\mathcal{D}} 
\oplus \underbrace{\mathcal{P}_{\mathcal{D}} \oplus \cdots \oplus \mathcal{P}_{\mathcal{D}}}_i] 
= \overline{(\underbrace{0, \cdots, 0}_{m-i}, \underbrace{1, \cdots, 1}_i)}
$$
De même, l'ensemble des sommets de $\widetilde{C}_0$ est $\lbrace S_0, S_1, \cdots, S_{2m-1} \rbrace$
tel que pour tout $i$ dans $\lbrace 0, 1, \cdots, 2m-1 \rbrace$, on a:
$$
S_i = [\mathcal{O}_{\Delta} \oplus \cdots \oplus \mathcal{O}_{\Delta} 
\oplus \underbrace{\mathcal{P}_{\Delta} \oplus \cdots \oplus \mathcal{P}_{\Delta}}_i]
= \overline{(\underbrace{0, \cdots, 0}_{2m-i}, \underbrace{1, \cdots, 1}_i)}
$$
\end{nota}

\begin{propo}\label{SommetsImageTRPair}
Il y a exactement $m$ sommets de $X_{\mathbb{K}}$ dans l'image de la chambre standard, 
$j(C_0)$, qui sont 
$S_0, S_2, \cdots, S_{2m-2}$.
\end{propo}

\begin{démo}
Pour tout $i$ dans $\lbrace 0, 1, \cdots, m-1 \rbrace$, on a 
$s_i = h_i.s_0$
où:
$$
h_i = {\rm diag} (1, \cdots, 1, \underbrace{\varpi_{\mathcal{D}}, \cdots, \varpi_{\mathcal{D}}}_i)
$$
Par ${\rm GL}_m (\mathcal{D})$-équivariance de $j$, on a $j(s_i) 
= \widetilde{\Phi} (h_i). j(s_0) = \widetilde{\Phi} (h_i).S_0$
avec $\widetilde{\Phi} (h_i)$ la matrice diagonale par blocs de ${\rm GL}_{2m} (\Delta)$:
$$
\widetilde{\Phi} (h_i) = {\rm diag} (I_2, \cdots, I_2, \underbrace{\Pi_{\Delta}, \cdots, \Pi_{\Delta}}_i)
$$
Il est alors immédiat que $j(s_i) = S_{2i}$
car $\Pi_{\Delta} . \mathcal{O}_{\mathcal{D}} = \mathcal{P}_{\mathcal{D}}$. 
Réciproquement, si $t$ dans $j(C_0)$ est aussi un sommet de $X_{\mathbb{K}}$, alors $t$ est à coordonnées entières et 
il existe $(\lambda_0, \cdots, \lambda_{m-1})$ dans $(\mathbb{R}^{+})^m$ tels que:
$$
t = \sum_{i=0}^{m-1} \lambda_i S_{2i} \, \, \text{et} \, \, 
\sum_{i=0}^{m-1} \lambda_i = 1
$$
Soit $i_0$ dans $\lbrace 0, \cdots, m-1 \rbrace$ le plus grand indice $k$ tel que $\lambda_k \neq 0$. Alors, la 
première coordonnée non nulle de $t$ est:
$$
\sum_{i=i_0}^{m-1} \lambda_i = \lambda_{i_0}
$$
Comme $t$ est à coordonnées entières, on a forcément $\lambda_{i_0} = 1$ 
et $t = S_{2i_0}$. 
Par conséquent, les seuls sommets de $X_{\mathbb{K}}$ qui sont aussi dans $j(C_0)$ sont:
$$
\lbrace j(s_0), j(s_1), \cdots, j(s_{m-1}) \rbrace 
= \lbrace S_0, S_2, S_4, \cdots, S_{2m-2} \rbrace
$$

\end{démo}

\subsubsection{Conditions de distinction.}

Avec des calculs analogues à la démonstration de \ref{PropositionIntersection}, 
on montre le résultat suivant :

\begin{lem}
 On a:
$$
{\rm GL}_m (\mathcal{D}) \cap \langle \varpi_{\mathbb{K}} \rangle 
{\rm GL}_{2m} (\mathcal{O}_{\Delta}) 
= \langle \varpi_{\mathcal{D}}^{d/2} \rangle 
{\rm GL}_{m} (\mathcal{O}_{\mathcal{D}})
$$
\end{lem}

\begin{propo}
On un isomorphisme de $\mathbb{C}$-espaces vectoriels:
$$
{\rm Hom}_{{\rm GL}_m (\mathcal{D})} (\pi, \mathds{1}) 
\simeq {\rm Hom}_{\langle \varpi_{\mathcal{D}}^{d/2} \rangle {\rm GL}_m (\mathcal{O}_{\mathcal{D}})} 
(\pi_0, \mathds{1})
$$
Ainsi, la représentation $\pi$ est ${\rm GL}_m (\mathcal{D})$-distinguée si et seulement si la représentation 
$\pi_0$ est $\langle \varpi_{\mathcal{D}}^{d/2} \rangle {\rm GL}_m (\mathcal{O}_{\mathcal{D}})$-distinguée.
\end{propo}

\begin{démo}
On pose $\mathcal{K} = \langle \varpi_{\mathbb{K}} \rangle {\rm GL}_{2m} (\mathcal{O}_{\Delta})$ et
 $K = {\rm GL}_{2m} (\mathcal{O}_{\Delta})$.
On remarque que 
$K$ est le sous-groupe parahorique $\mathcal{A} (j(s_0))^{\times}$,
donc pour tout $g$ dans ${\rm GL}_{2m} (\Delta)$, 
$gKg^{-1} = \mathcal{A} (g.j(s_0))^{\times}$.
En raisonnant à nouveau comme en \ref{RaisonnementDistinction} et en utilisant la proposition 
\ref{SommetsImageTRPair}, on vérifie que:
$$
{\rm Hom}_{{\rm GL}_m (\mathcal{D})} (\pi, \mathds{1}) 
\simeq {\rm Hom}_{{\rm GL}_m (\mathcal{D}) \cap \langle \varpi_{\mathbb{K}} \rangle 
{\rm GL}_{2m} (\mathcal{O}_{\Delta}) } 
(\pi_0, \mathds{1})
$$
Le lemme précédent nous permet de conclure directement.

\end{démo}

\begin{rmq}
En remarquant que $(\varpi_{\mathcal{D}}^{d/2})^2 = \varpi_{\mathbb{F}}$ et que 
$\mathcal{O}_{\mathbb{F}}^{\times} \subseteq \mathcal{O}_{\mathcal{D}}^{\times} 
\subseteq {\rm GL}_m (\mathcal{O}_{\mathcal{D}})$, on vérifie facilement que 
si $\pi$ est ${\rm Gl}_m (\mathcal{D})$-distinguée alors son caractère central $\chi$ 
est trivial sur $\mathbb{F}^{\times}$.
\end{rmq}

\begin{nota}
On note $l = k_{\mathbb{K},n}$ et $l_0$ une extension de degré $m$ de $k_{\Delta}$ telle que $l_0 \subseteq l$. 
On a alors le diagramme d'extensions de corps finis suivant:
$$
\xymatrix{ & l = k_{\mathbb{K},n} = k_{\mathcal{D},m} = k_{\Delta,2m} & \\ 
         k_{\mathcal{D}} \ar@{-}[ru]^m \ar@{-}[rd]_2  
         &  & l_0 \ar@{-}[lu]^2 \ar@{-}[ld]_m\\
          & k_{\Delta} & }
$$
On a $k_{\mathcal{D}} = k_{\Delta} [\alpha]$, avec $\alpha^2 \in k_{\Delta}$. 
On peut voir $k_{\mathcal{D}}$ comme une sous-$k_{\Delta}$-algèbre de ${\rm M}_2 (k_{\Delta})$:
$$
k_{\mathcal{D}} \simeq \left\lbrace \left( \begin{array}{cc}
x & y \alpha^2 \\
y & x \\
\end{array} \right): x, y \in k_{\Delta} \right\rbrace
$$
On peut donc injecter ${\rm M}_m (k_{\mathcal{D}})$ dans ${\rm M}_{2m} (k_{\Delta})$ par blocs.
On notera:
$$
\beta = \left( \begin{array}{cc}
0 & \alpha^2 \\
1 & 0 \\
\end{array} \right) \, \, \text{et} \, \, 
w = {\rm diag} (\beta, \cdots, \beta) \in {\rm GL}_{2m} (k_{\Delta})
$$
On vérifie facilement que le commutant de $\beta$ dans ${\rm M}_{2} (k_{\Delta})$ est $k_{\mathcal{D}}$.\\ 
On note $\tau : {\rm GL}_{2m} (k_{\Delta}) \rightarrow {\rm GL}_{2m} (k_{\Delta}), \, x \mapsto wxw^{-1}$. 
Alors $\tau$, la restriction de ${\rm Ad} (w)$ (la conjugaison par $w$) à ${\rm GL}_{2m} (k_{\Delta})$, 
 est une involution dont l'ensemble des points fixes est ${\rm GL}_m (k_{\mathcal{D}})$. \\
Enfin, on notera $T$ l'image de $l^{\times}$ dans ${\rm GL}_m (k_{\mathcal{D}}) \subseteq {\rm GL}_{2m} (k_{\Delta})$ 
et ${\rm N}_{{\rm GL}_{2m} (k_{\Delta})} (T)$ le normalisateur de $T$ dans ${\rm GL}_{2m} (k_{\Delta})$.
\end{nota}

\begin{lem}\label{LemmeGLmkdDistinction}
La représentation $\overline{\gamma}_0$ de ${\rm GL}_{2m} (k_{\Delta})$ est 
${\rm GL}_m (k_{\mathcal{D}})$-distinguée si et seulement si $\overline{\chi}$ est trivial sur 
$l_0^{\times}$. Dans ce cas, ${\rm Hom}_{{\rm GL}_m (k_{\mathcal{D}})} (\overline{\gamma}_0, \mathds{1})$ 
est de dimension $1$.
\end{lem}

\begin{démo}
On cherche une condition nécessaire et suffisante pour que $\overline{\gamma}_0$ soit 
${\rm GL}_m (k_{\mathcal{D}})$-distinguée. Pour cela, d'après l'article \cite{Lusztig}, 
il faut déterminer l'ensemble suivant:
$$
\Xi = \{ g \in {\rm GL}_{2m} (k_{\Delta}) : \tau (g^{-1} T g) = g^{-1} T g \} 
$$
lorsque $\Xi$ est vu comme partie du double quotient 
$T \backslash {\rm GL}_{2m} (k_{\Delta}) / {\rm GL}_m (k_{\mathcal{D}})$. 
Puisque $w \in T$, on a, pour tout $x$ dans ${\rm GL}_{2m} (k_{\Delta})$, 
$x \in \Xi$ si et seulement si $xwx^{-1} \in {\rm N}_{{\rm GL}_{2m} (k_{\Delta})} (T)$.
En utilisant Skölem-Noether, on vérifie qu'il existe $g_0 \in {\rm GL}_{2m} (k_{\Delta})$ 
tel que la restriction de ${\rm Ad} (g_0)$ à $T$ soit d'ordre 
$2m$ et égale au Frobenius de $l$ sur $k_{\Delta}$ (${\rm Ad} (g_0)_{\vert T} = {\rm Frob}_{l/k_{\Delta}}$) 
de sorte que ${\rm N}_{{\rm GL}_{2m} (k_{\Delta})} (T)$ soit de la forme suivante:
$$
{\rm N}_{{\rm GL}_{2m} (k_{\Delta})} (T) = T. \langle g_0 \rangle
$$
Soit $x \in \Xi$, alors il existe $t$ dans $T$ et $k$ dans $\mathbb{Z}$ tels que 
$x w x^{-1} = t g_0^k$. Ainsi:
$$
(x w x^{-1})^2 = w^2 = \alpha^2 I_{2m} = t g_0^k t g_0^{-k} g_0^{2k} 
= \underbrace{(t {\rm Frob}_{l/k_{\Delta}}^k (t))}_{t^{'} \in T} g_0^{2k} 
= t^{'} g_0^{2k}
$$
On en déduit que $g_0^{2k} = (t^{'})^{-1} w^2 \in T$. 
Puisque $T$ est commutatif, pour tout $t$ dans $T$:
$$
{\rm Ad} (g_0^{2k}) (t) = t = ({\rm Ad} (g_0)_{\vert T})^{2k} (t)
$$
Or ${\rm Ad} (g_0)_{\vert T}$ est d'ordre $2m$, donc nécessairement $m$ divise $k$. 
 On peut ainsi supposer que $k=0$ ou $k=m$.
\begin{itemize}
\item[$\ast$] Supposons que $k=0$. 
Dans ce cas, $x w x^{-1} = t \in T$ donc commute avec $w$. 
Or, $t^2 = w^2$ donc $(w^{-1} t)^2 = I_{2m}$. 
Ainsi, $w^{-1} t$ est un élément d'ordre $1$ ou $2$ de $T = l^{\times}$ qui est un groupe cyclique. 
Il y a seulement deux possibilités: ou bien $t = w$ 
ou bien $t = -w$.\\
Si $t=w$, alors $xwx^{-1} = w$, donc $wxw^{-1} = \tau (x) = x$ et 
$x \in {\rm GL}_m (k_{\mathcal{D}})$. Comme on cherche $x$ modulo un élément de 
${\rm GL}_m (k_{\mathcal{D}})$ à droite, on peut se restreindre au cas où $x = I_{2m}$.\\
Sinon, $t = -w$, alors $xwx^{-1} = -w$ et $\tau (x^{-1}) = - x^{-1}$. \\
On définit par blocs $x^{-1} = [X_{i,j}]_{1 \leq i, j \leq m}$, $X_{i,j} \in {\rm M}_2 (k_{\Delta})$.
Après calcul, on constate que $xwx^{-1} = -w$ si et seulement si $x^{-1} = [X_{i,j}]$ avec 
$X_{i,j} = (s_{i,j} + z_{i,j} \beta) \alpha_0$ (où $s_{i,j}, z_{i,j}$ appartiennent à $k_{\Delta}$) si et seulement s'il 
existe $t$ dans ${\rm GL}_{m} (k_{\mathcal{D}})$ tel que:
$$
x^{-1} = t. \Lambda_0, \, \, \text{où} \, \, 
\Lambda_{0} = {\rm diag} (\alpha_0, \cdots, \alpha_0)
$$
(car $(1, \beta)$ est une $k_{\Delta}$-base de $k_{\mathcal{D}}$). Ainsi, 
en remarquant que $\Lambda_0^{-1} = \Lambda_0$:
$$
\{ x \in {\rm GL}_{2m} (k_{\Delta}) : xwx^{-1} = -w \} 
= \{ \Lambda_0 t : t \in {\rm GL}_m (k_{\mathcal{D}}) \}
$$
Comme on cherche $x$ modulo la multiplication par un élément de 
${\rm GL}_m (k_{\mathcal{D}})$ à droite, on peut se restreindre au cas 
$x = \Lambda_0$.
\item[$\ast$] Supposons maintenant que $k=m$.
On a $xwx^{-1} = t g_0^m$ donc:
$$
(xwx^{-1})^2 = w^2 = t {\rm Frob}_{l/k_{\Delta}}^m (t) g_0^{2m} 
= t {\rm Frob}_{l/k_{\Delta}}^m (t) = t {\rm Frob}_{l/l_0} (t)
= {\rm N}_{l/l_0} (t)
$$
Fixons $t_0 \in T$ tel que ${\rm N}_{l/l_0} (t_0) = w^2$, alors 
${\rm N}_{l/l_0} (tt_0^{-1}) = 1$, par théorème $90$ de Hilbert, il existe $y$ dans $l^{\times}$ tel que:
$$
tt_0^{-1} = \frac{y}{{\rm Frob}_{l/l_0} (y)} = y g_0^m y^{-1} g_0^{-m}
$$
Par conséquent, comme $T$ est commutatif, $y$ et $t_0$ commutent:
$$
xwx^{-1} = t_0 y g_0^{m} y^{-1} = y t_0 g_0^{m} y^{-1} 
\Rightarrow (y^{-1} x)w(y^{-1} x)^{-1} = t_0 g_0^{m}
$$
Comme $y^{-1} \in T$ et que l'on cherche $x$ modulo la multiplication par un élément de $T$ à gauche, il y a 
un seul cas à considérer: le cas où $xwx^{-1} = t_0 g_0^m$.\\
Déterminons combien il y a d'éléments $x$ dans le double quotient 
$T \backslash {\rm GL}_{2m} (k_{\Delta}) / {\rm GL}_{m} (k_{\mathcal{D}})$ tels que 
$xwx^{-1} = t_0 g_0^m$. 
Soient $x_1$ et $x_2$ dans ${\rm GL}_{2m} (k_{\Delta})$ tels que 
$x_1 w x_1^{-1} = t_0 g_0^m = x_2 w x_2^{-1}$, alors, puisque $k_{\mathcal{D}} = k_{\Delta} (w)$:
$$
x_1 k_{\mathcal{D}} x_1^{-1} = k_{\Delta} (x_1 w x_1^{-1}) = k_{\Delta} (x_2 w x_2^{-1}) 
= x_2 k_{\mathcal{D}} x_2^{-1}
$$ 
On en déduit que 
${\rm Ad} (x_2^{-1} x_1)$, la conjugaison par $x_2^{-1} x_1$, fixe point par point les éléments de 
$k_{\Delta}$. 
Par théorème de Skölem-Noether, on sait qu'il existe $\gamma$ dans ${\rm GL}_{2m} (k_{\Delta})$ tel que 
${\rm Frob}_{k_{\mathcal{D}} / k_{\Delta}}$ soit la restriction de ${\rm Ad} (\gamma)$ à $k_{\mathcal{D}}$.
 On peut également supposer que ${\rm Frob}_{k_{\mathcal{D}} / k_{\Delta}} (w) = -w$.\\
Si ${\rm Ad} (x_2^{-1} x_1)$ fixe tous les éléments de 
$k_{\mathcal{D}}$ alors $x_2^{-1} x_1 \in {\rm GL}_{m} (k_{\mathcal{D}})$   
et dans ce cas $x_1$ et $x_2$ sont dans la même double classe dans 
$T \backslash {\rm GL}_{2m} (k_{\Delta}) / {\rm GL}_{m} (k_{\mathcal{D}})$.\\
Sinon, ${\rm Ad} (x_2^{-1} x_1)_{\vert k_{\mathcal{D}}} = {\rm Ad} (\gamma)_{\vert k_{\mathcal{D}}}$ et donc 
$x_1 \in x_2 \gamma {\rm GL}_{m} (k_{\mathcal{D}})$.
Or: 
$$
(x_2 \gamma) w (x_2 \gamma)^{-1} = - x_2 w x_2^{-1} \neq t_0 g_0^m
$$
On déduit de tout ceci qu'on peut se restreindre au cas $x = x_0$ avec $x_0 w x_0^{-1} = t_0 g_0^m$.\\
Ensuite, on vérifie facilement qu'un tel élément $x$ tel que $xwx^{-1} = t_0 g_0^m$ existe bien. 
En effet, posons $\delta = t_0 g_0^m$, alors $\delta^2 = w^2 \in k_{\Delta}$ et n'est pas un carré dans 
$k_{\Delta}$, donc engendre une extension (de corps) de degré $2$ de $k_{\Delta}$ et:
$$
k_{\Delta} [\delta] = k_{\Delta} [w] = k_{\mathcal{D}}
$$
Les éléments $\delta$ et $w$ ont même polynôme minimal sur $k_{\Delta}$, à savoir 
$X^2 - w^2$. On en déduit un morphisme de $k_{\Delta}$-algèbres (bien défini):
$$
\psi : k_{\Delta} [w] \rightarrow k_{\Delta} [\delta], P(w) \mapsto P(\delta)
$$
et donc un morphisme injectif $\psi : k_{\Delta} [w] \hookrightarrow {\rm GL}_{2m} (k_{\Delta})$ 
où ${\rm GL}_{2m} (k_{\Delta})$ est une $k_{\Delta}$-algèbre centrale simple. 
Par théorème de Skölem-Noether, il existe $x$ dans ${\rm GL}_{2m} (k_{\Delta})$ tel que:
$$
\forall z \in k_{\Delta} [w], \psi (z) = x z x^{-1}
$$
En particulier, $\psi (w) = \delta = x w x^{-1} = t_0 g_0^m$.\\

\end{itemize}

Nous allons en déduire ${\rm dim} ({\rm Hom}_{{\rm GL}_m (k_{\mathcal{D}})} (\overline{\gamma}_0, \mathds{1}))$. 
D'après l'article de Lusztig, on sait que:
$$
{\rm dim} ({\rm Hom}_{{\rm GL}_m (k_{\mathcal{D}})} (\overline{\gamma}_0, \mathds{1})) 
= \sum_{x \in \Xi} r(x)
$$
où $r(x) \in \{ -1, 0, 1 \}$. D'après ce qui précède, il y a seulement $3$ élements de $\Xi$ à considérer.\\
Si $x = {\rm Id}$, on sait que $r(x) \neq 0$ si et seulement si la restriction de 
$\overline{\chi}$ à $T \cap {\rm GL}_m (k_{\mathcal{D}})$ et la fonction $\varepsilon_{T}$ de 
Lusztig coïncident.
Comme $T$ est connexe, $\varepsilon_{T}$ est la fonction 
constante égale à $1$ sur $T \cap {\rm GL}_m (k_{\mathcal{D}}) = T$, or le caractère $\overline{\chi}$ ne 
peut pas être trivial sur $T$.
Finalement, $r({\rm Id}) = 0$.\\
De même, on vérifie que si $x = \Lambda_0$, on a 
$\Lambda_0^{-1} T \Lambda_0 \cap {\rm GL}_m (k_{\mathcal{D}}) = \Lambda_0^{-1} T \Lambda_0$ donc est connexe. 
Par suite, $r(\Lambda_0) = 0$.\\
Finalement, il y a un seul élément qui intervient dans le calcul de 
${\rm dim} ({\rm Hom}_{{\rm GL}_m (k_{\mathcal{D}})} (\overline{\gamma}_0, \mathbb{C}))$. 
Fixons $x_0$ dans ${\rm GL}_{2m} (k_{\Delta})$ tel que $x_0 w x_0^{-1} = t_0 g_0^m$. 
Soit $t \in T$, alors, en utilisant la commutativité de $T$, on montre que:
\begin{eqnarray*}
x_0^{-1} t x_0 \in x_0^{-1} T x_0 \cap {\rm GL}_m (k_{\mathcal{D}}) 
& \Leftrightarrow & (x_0 w x_0^{-1})t = t (x_0 w x_0^{-1}) 
  \Leftrightarrow t_0 g_0^m t = t t_0 g_0^m\\
& \Leftrightarrow & t_0 {\rm Frob}_{l/l_0} (t) = t t_0 = t_0 t
  \Leftrightarrow t \in l_0^{\times}\\
\end{eqnarray*}
On en déduit que $x_0^{-1} T x_0 \cap {\rm GL}_m (k_{\mathcal{D}}) = x_0^{-1} l_0^{\times} x_0$. 
Comme $l_0^{\times}$ est connexe, la fonction $\varepsilon_{x_0^{-1} T x_0}$ est triviale sur 
$x_0^{-1} T x_0 \cap {\rm GL}_m (k_{\mathcal{D}}) = x_0^{-1} l_0^{\times} x_0$. 
On en déduit que $r (x_0) \neq 0$ si et seulement si $\overline{\chi}$ est trivial sur $l_0^{\times}$, 
et dans ce cas:
$$
{\rm dim} ({\rm Hom}_{{\rm GL}_m (k_{\mathcal{D}})} (\overline{\gamma}_0, \mathds{1})) 
= r({\rm Id}) + r(\Lambda_0) + r(x_0) = r(x_0) = 1
$$

\end{démo}

\begin{nota}
On notera $l$ le corps résiduel de $\mathbb{K}_n$ (extension de degré $2m$ de $k_{\Delta}$) 
et $l_0$ celui de $\mathbb{K}_{n/2}$. 
On fixe $\eta$ dans $l \backslash l_0$ tel que $\eta^2 \in l_0^{\times}$. Alors 
$\varpi_{\mathbb{K}} \eta$ est une uniformisante de $\mathbb{E}$ où $\mathbb{E}$ est défini par le diagramme d'extensions 
de corps suivant:
$$
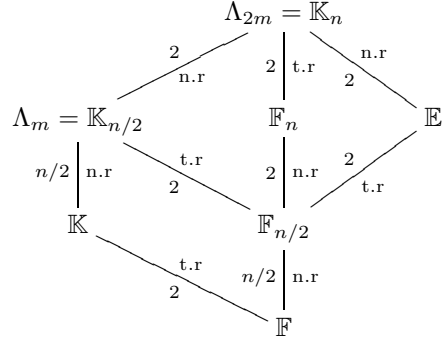
\xymatrix{ & \Lambda_{2m} = \mathbb{K}_n & \\ 
         \Lambda_{m} = \mathbb{K}_{n/2} \ar@{-}[ru]_{{\rm n.r}}^2 \ar@{-}[rd]_2^{{\rm t.r}} \ar@{-}[d]_{n/2}^{{\rm n.r}} 
         & \mathbb{F}_n  \ar@{-}[u]_{{\rm t.r}}^2 \ar@{-}[d]_2^{{\rm n.r}} 
         & \mathbb{E} \ar@{-}[lu]_{{\rm n.r}}^2 \ar@{-}[ld]_2^{{\rm t.r}}\\
         \mathbb{K} \ar@{-}[rd]_2^{{\rm t.r}} & \mathbb{F}_{n/2} \ar@{-}[d]_{n/2}^{{\rm n.r}} & \\
          & \mathbb{F} & }
$$
Par théorème de Skölem-Noether, on sait que l'on peut fixer $\gamma$ dans ${\rm GL}_{2m} (k_{\Delta})$ tel que 
la restriction de ${\rm Ad} (\gamma)$, la conjugaison par $\gamma$, à $k_{\mathcal{D}}$ soit 
le Frobenius ${\rm Frob}_{k_{\mathcal{D}} / k_{\Delta}}$. En particulier, la conjugaison par $\gamma^2$ 
fixe point par point les éléments de $k_{\mathcal{D}}$, donc $\gamma^2 \in {\rm GL}_m (k_{\mathcal{D}})$ 
et, puisque $k_{\mathcal{D}} = k_{\Delta} [w]$, on peut supposer que $\gamma w \gamma^{-1} = -w$.\\
Soit $\mathcal{L} = (L_i)_{i \in \mathbb{Z}}$ la chaîne d'$\mathcal{O}_{\Delta}$-réseaux de 
$\mathcal{D}^m$ définie par $L_i = \mathcal{P}_{\mathcal{D}}^i \oplus \cdots \oplus \mathcal{P}_{\mathcal{D}}^i$. 
Avec l'identification $\mathcal{P}_{\mathcal{D}} \simeq \mathcal{O}_{\Delta} \oplus \mathcal{O}_{\Delta}$, 
on remarque que:
$$
\forall i \in \mathbb{Z}, \varpi_{\mathbb{K}}^{-1} \varpi_{\mathcal{D}}^{d/2} L_i = L_i 
\Rightarrow \varpi_{\mathbb{K}}^{-1} \varpi_{\mathcal{D}}^{d/2} \in \mathcal{A} (\mathcal{L})^{\times} 
= {\rm GL}_{2m} (\mathcal{O}_{\Delta})
$$
Soit $u$ dans ${\rm GL}_{2m} (\mathcal{O}_{\Delta})$, 
tel que $\varpi_{\mathbb{K}} u = \varpi_{\mathcal{D}}^{d/2}$. On sait que la conjugaison par $\varpi_{\mathcal{D}}^{d/2}$ 
et donc par $\varpi_{\mathbb{K}} u$ engendre ${\rm Frob}_{k_{\mathcal{D}} / k_{\Delta}}$. 
Ainsi, pour tout $x$ dans $k_{\mathcal{D}}$, en utilisant que $\varpi_{\mathbb{K}}$ est central, on a:
$$
\gamma x \gamma^{-1} = \varpi_{\mathbb{K}} u x u^{-1} \varpi_{\mathbb{K}}^{-1} = u x u^{-1} 
\Rightarrow (\gamma^{-1} u) x (\gamma^{-1} u)^{-1} = x
$$
On en déduit que $u \in \gamma {\rm GL}_m (\mathcal{O}_{\mathcal{D}})$ et donc 
$\langle \varpi_{\mathcal{D}}^{d/2} \rangle {\rm GL}_m (\mathcal{O}_{\mathcal{D}}) 
= \langle \varpi_{\mathbb{K}} u \rangle {\rm GL}_m (\mathcal{O}_{\mathcal{D}}) 
= \langle \varpi_{\mathbb{K}} \gamma \rangle {\rm GL}_m (\mathcal{O}_{\mathcal{D}})$. 
On définit $\upsilon$ un caractère de $\langle \varpi_{\mathbb{K}} \gamma \rangle {\rm GL}_m (\mathcal{O}_{\mathcal{D}})$ par:
$$
\upsilon ((\varpi_{\mathbb{K}} \gamma)^s x) = \chi ((\varpi_{\mathbb{K}} \eta)^s) 
\chi_{\Lambda_{2m} / \mathbb{E}} ((\varpi_{\mathbb{K}} \eta)^s)
$$
pour tout $x$ dans ${\rm GL}_m (\mathcal{O}_{\mathcal{D}})$ et tout $s$ dans $\mathbb{Z}$, 
où $\chi_{\mathbb{K}_n / \mathbb{E}} = \chi_{\Lambda_{2m} / \mathbb{E}}$ est le caractère quadratique de 
$\mathbb{K}_n / \mathbb{E}$.
Comme $\mathbb{K}_n / \mathbb{E}$ est non ramifiée, la norme ${\rm N}_{\mathbb{K}_n / \mathbb{E}}$ 
est surjective sur les unités.
Puisque $\varpi_{\mathbb{K}} \eta$ est 
une uniformisante de $\mathbb{E}$, on a:
$$
\upsilon (\varpi_{\mathbb{K}} \gamma) = - \chi (\varpi_{\mathbb{K}}) \overline{\chi} (\eta)
$$
\end{nota}

En utilisant la même idée que dans \cite{HakimMurnaghan1} lors de la démonstration de la proposition 6.3, on montre :

\begin{lem}
L'inclusion canonique:
$$
{\rm Hom}_{\langle \varpi_{\mathbb{K}} \gamma \rangle {\rm GL}_m (\mathcal{O}_{\mathcal{D}})} (\pi_0, \upsilon) 
\subseteq {\rm Hom}_{{\rm GL}_m (k_{\mathcal{D}})} (\overline{\gamma}_0, \mathds{1})
$$
est en fait un isomorphisme de $\mathbb{C}$-espaces vectoriels.
\end{lem}

\begin{démo}
On a clairement l'inclusion canonique:
$$
{\rm Hom}_{\langle \varpi_{\mathbb{K}} \gamma \rangle {\rm GL}_m (\mathcal{O}_{\mathcal{D}})} (\pi_0, \upsilon) 
\subseteq {\rm Hom}_{{\rm GL}_m (k_{\mathcal{D}})} (\overline{\gamma}_0, \mathds{1})
$$
Si $\overline{\gamma}_0$ n'est pas ${\rm GL}_m (k_{\mathcal{D}})$-
distinguée, alors:
$$
{\rm Hom}_{\langle \varpi_{\mathbb{K}} \gamma \rangle {\rm GL}_m (\mathcal{O}_{\mathcal{D}})} (\pi_0, \upsilon) 
= {\rm Hom}_{{\rm GL}_m (k_{\mathcal{D}})} (\overline{\gamma}_0, \mathds{1}) 
= \{ 0 \}
$$
Sinon, supposons que $\overline{\gamma}_0$ est ${\rm GL}_m (k_{\mathcal{D}})$-distinguée, alors 
$\overline{\chi}$ est trivial sur $l_0^{\times}$.
Soit $\varphi : V \rightarrow \mathbb{C}$ une forme linéaire non nulle dans 
${\rm Hom}_{{\rm GL}_m (k_{\mathcal{D}})} (\overline{\gamma}_0, \mathds{1})$ 
(où $V$ est l'espace de $\pi_0$). 
Posons $\psi = \varphi \circ \pi_0 (\gamma)$. L'application $\psi$ est aussi une forme linéaire non nulle. 
Alors pour tout $v$ dans $V$, pour tout $x$ dans ${\rm GL}_m (\mathcal{O}_{\mathcal{D}})$, on pose 
$u = \pi_0 (\gamma).v \in V$ et on a:
$$
\psi (\overline{\gamma}_0 (\overline{x}).v) = \psi (\pi_0 (x).v) = \varphi (\pi_0 (\gamma x \gamma^{-1}).u) 
= \varphi (u) = \psi (v)
$$
On en déduit que $\psi \in {\rm Hom}_{{\rm GL}_m (k_{\mathcal{D}})} (\overline{\gamma}_0, \mathds{1})$, 
qui est un espace de dimension $1$, donc $\psi$ et $\varphi$ sont colinéaires. Soit $c \in \mathbb{C}^{\times}$ 
tel que $\psi = c \varphi$. Comme $\gamma^2 \in {\rm GL}_m (k_{\mathcal{D}})$, on a:
$$
\forall v \in V, \varphi (v) = \varphi (\pi_0 (\gamma^2).v) = c^2 \varphi (v)
$$
On en déduit que $c^2 = 1$ et, par suite, $\psi = \varphi$ ou $\psi = - \varphi$. 
Pour montrer que les deux espaces 
${\rm Hom}_{\langle \varpi_{\mathbb{K}} \gamma \rangle {\rm GL}_m (k_{\mathcal{D}})} (\pi_0, \upsilon)$ 
et ${\rm Hom}_{{\rm GL}_m (k_{\mathcal{D}})} (\overline{\gamma}_0, \mathds{1})$ sont égaux, il suffit de montrer 
que $\varphi$ appartient à  
${\rm Hom}_{\langle \varpi_{\mathbb{K}} \gamma \rangle {\rm GL}_m (k_{\mathcal{D}})} (\pi_0, \upsilon)$. 
Pour cela, il faut vérifier que pour tout $x$ dans ${\rm GL}_m (k_{\mathcal{D}})$ et tout $u$ dans $V$:
$$
\varphi (\pi_0 (\varpi_{\mathbb{K}} \gamma) \pi_0 (x).u) = \upsilon (\varpi_{\mathbb{K}} \gamma x) \varphi (u) 
\, \, \text{i.e} \, \, 
\chi (\varpi_{\mathbb{K}}) \psi (u) = - \chi (\varpi_{\mathbb{K}}) \overline{\chi} (\eta) \varphi (u)
$$
Donc cela revient à montrer que $\psi = \varphi$ si et seulement si $\overline{\chi} (\eta) = -1$.\\
Posons $\mathcal{K} = \langle \gamma \rangle {\rm GL}_m (k_{\mathcal{D}})$. 
On a $\psi = \varphi$ si et seulement si pour tout $v$ dans $V$, 
$\varphi (\pi_0 (\gamma).v) = \varphi (v)$ i.e si 
$\varphi$ appartient à ${\rm Hom}_{\mathcal{K}} (\pi_0, \mathds{1})$. 
Or, il est clair que 
${\rm Hom}_{\mathcal{K}} (\pi_0, \mathds{1}) \subseteq 
{\rm Hom}_{{\rm GL}_m (k_{\mathcal{D}})} (\overline{\gamma}_0, \mathds{1})$ 
et comme le dernier espace est de dimension $1$, on en déduit immédiatement que 
$\psi = \varphi$ si et seulement si 
${\rm Hom}_{\mathcal{K}} (\pi_0, \mathds{1}) \neq 0$.
On remarque que $\overline{\chi}$ est trivial sur $l_0^{\times}$, donc sur $k_{\Delta}^{\times}$, 
qui est le centre de ${\rm GL}_{2m} (k_{\Delta})$. 
On peut donc voir $\pi_0$ comme une représentation de 
$\overline{{\rm GL}_{2m} (k_{\Delta})} = {\rm PGL}_{2m} (k_{\Delta})$. 
On pose $\overline{\mathcal{K}} = \mathcal{K} / k_{\Delta}^{\times}$ l'image de 
$\mathcal{K}$ dans ${\rm PGL}_{2m} (k_{\Delta})$. 
Notre objectif est donc de montrer que 
${\rm Hom}_{\overline{\mathcal{K}}} (\pi_0, \mathds{1}) \neq 0$ si et seulement si 
$\overline{\chi} (\eta) = -1$. 
Notons $\overline{\tau}$ la conjugaison par $w$ dans ${\rm PGL}_{2m} (k_{\Delta})$. 
Il est clair que $\overline{\tau}$ fixe point par point les éléments de 
$\overline{{\rm GL}_m (k_{\mathcal{D}})}$. De plus:
$$
\gamma w \gamma^{-1} = -w 
\Rightarrow w \gamma w^{-1} = - \gamma 
\Rightarrow \overline{\tau} (\overline{\gamma}) = \overline{\gamma}
$$
Par conséquent, $\overline{\mathcal{K}} \subseteq {\rm PGL}_{2m} (k_{\Delta})^{\overline{\tau}}$ 
(les points fixes de $\overline{\tau}$ dans ${\rm PGL}_{2m} (k_{\Delta})$) et 
$\overline{\mathcal{K}}$ contient la composante neutre de ${\rm PGL}_{2m} (k_{\Delta})^{\overline{\tau}}$. 
On va donc pouvoir appliquer les résultats de l'article de Lusztig \cite{Lusztig} à l'espace symétrique 
$({\rm PGL}_{2m} (k_{\Delta}), \overline{\mathcal{K}})$. 
Notons $\overline{T}$ l'image de $l^{\times}$ dans ${\rm PGL}_{2m} (k_{\Delta})$ et:
$$
\overline{\Xi} = \{ \overline{g} \in {\rm PGL}_{2m} (k_{\Delta}) : 
\overline{\tau} (\overline{g}^{-1} \overline{T} \overline{g}) = \overline{g}^{-1} \overline{T} \overline{g} \}
$$
vue comme partie du double quotient 
$\overline{T} \backslash {\rm PGL}_{2m} (k_{\Delta}) / \overline{{\rm GL}_{m} (k_{\mathcal{D}})}$. 
En reprenant les calculs de la démonstration du lemme \ref{LemmeGLmkdDistinction}, on vérifie 
qu'il n'y a qu'une seule double classe à considérer : on fixe $\overline{x_0} \in 
 {\rm PGL}_{2m} (k_{\Delta})$ tel que (en reprenant les notations de la démonstration du lemme \ref{LemmeGLmkdDistinction}) 
$\overline{x_0 w x_0^{-1}} = \overline{t_0 g_0^m}$. Alors, il est clair que pour tout 
$t$ dans $T$, si $\overline{x_0^{-1} t x_0} \in 
\overline{x_0}^{-1} \overline{T} \overline{x_0} \cap \overline{\mathcal{K}}$ alors 
$\overline{{\rm Frob}_{l/l_0} (t)} = \overline{t}$ donc il existe $c$ dans 
$k_{\Delta}^{\times}$ tel que ${\rm Frob}_{l/l_0} (t) = c t$. En utilisant le fait que 
$(1, \eta)$ est une $l_0$-base de $l$, on montre que les seuls éléments de $l^{\times}$ vérifiant 
${\rm Frob}_{l/l_0} (t) = c t$ sont les éléments de $l_0^{\times} \cup \eta l_0^{\times}$. 
Réciproquement, lors de la démonstration de \ref{LemmeGLmkdDistinction}, nous avons vu que 
$x_0^{-1} l_0^{\times} x_0 \subseteq {\rm GL}_{m} (k_{\mathcal{D}})$. 
Il nous reste donc à vérifier que $x_0^{-1} \eta x_0 \in \langle \gamma \rangle {\rm GL}_{m} (k_{\mathcal{D}})$. 
On remarque que:
$$
(x_0 w x_0^{-1}) \eta (x_0 w x_0^{-1})^{-1} = t_0 (g_0^{m} \eta g_0^{-m}) t_0^{-1} 
= t_0 {\rm Frob}_{l/l_0} (\eta) t_0^{-1} = t_0 (-\eta) t_0^{-1} = -\eta
$$
car $T$ est commutatif. On en déduit que 
$\tau (x_0^{-1} \eta x_0) = - x_0^{-1} \eta x_0$. Puisque $\tau (\gamma) = - \gamma$, on a:
$$
\tau (\gamma (x_0^{-1} \eta x_0)) = \tau (\gamma) \tau (x_0^{-1} \eta x_0) = 
(-\gamma) (-x_0^{-1} \eta x_0) = \gamma (x_0^{-1} \eta x_0)
$$
Par conséquent $\gamma (x_0^{-1} \eta x_0) \in {\rm GL}_m (k_{\mathcal{D}})$ et 
$x_0^{-1} \eta x_0 \in \langle \gamma \rangle {\rm GL}_{m} (k_{\mathcal{D}})$.
 Finalement:
$$
\overline{x_0}^{-1} \overline{T} \overline{x_0} \cap \overline{\mathcal{K}} = 
(\overline{x_0}^{-1} l_0^{\times} \overline{x_0}) / k_{\Delta}^{\times} \cup 
(\overline{x_0}^{-1} \eta l_0^{\times} \overline{x_0}) / k_{\Delta}^{\times}
$$
D'après la formule de Lusztig, 
${\rm dim} ({\rm Hom}_{\overline{\mathcal{K}}} (\pi_0, \mathds{1})) \neq 0$ si
 et seulement si $\overline{\chi}^{\overline{x_0}}$ et la fonction 
$\varepsilon_{\overline{x_0}^{-1} \overline{T} \overline{x_0}}$ coïncident sur 
$\overline{x_0}^{-1} \overline{T} \overline{x_0} \cap \overline{\mathcal{K}}$. 
Par hypothèse, on sait que $\overline{\chi}$ est trivial sur $l_0^{\times}$ 
donc $\overline{\chi}^{\overline{x_0}}$ est trivial sur 
$(\overline{x_0}^{-1} l_0^{\times} \overline{x_0}) / k_{\Delta}^{\times}$ et, comme 
$(\overline{x_0}^{-1} l_0^{\times} \overline{x_0}) / k_{\Delta}^{\times}$ est connexe, la fonction 
$\varepsilon_{\overline{x_0}^{-1} \overline{T} \overline{x_0}}$ est aussi triviale sur 
$(\overline{x_0}^{-1} l_0^{\times} \overline{x_0} / k_{\Delta}^{\times})$. 
En utilisant la définition de la fonction $\varepsilon$, on montre facilement que  
$\varepsilon_{\overline{x_0}^{-1} \overline{T} \overline{x_0}} (x_0^{-1} \eta x_0) 
= \varepsilon_{\overline{T}} (\eta)$. Un calcul rapide (cf. \cite{HakimMurnaghan1} démonstration de 
la proposition 6.3) nous permet de vérifier que 
$\varepsilon_{\overline{T}} (\eta) = -1$. On en déduit immédiatement que 
 ${\rm dim} ({\rm Hom}_{\overline{\mathcal{K}}} (\pi_0, \mathds{1})) \neq 0$ si et seulement si 
$\overline{\chi} (\eta) = -1$.\\
On a bien le résultat annoncé.

\end{démo}

Avec un raisonnement analogue à la démonstration du théorème \ref{CritereDistinctionTRdImpair}, le lemme 
précédent nous permet de montrer le théorème suivant :

\begin{theo}\label{CritereDistinctionTRPair}
La représentation $\pi$ est ${\rm GL}_m (\mathcal{D})$-distinguée si et seulement si 
$\chi$ est trivial sur $\mathbb{F}^{\times}$, $\overline{\chi}$ est trivial sur 
$l_0^{\times}$ et $\chi (\varpi_{\mathbb{K}}) \overline{\chi} (\eta) = -1$.\\
De plus, si $\pi$ est ${\rm GL}_m (\mathcal{D})$-distinguée, on a:
$$
{\rm dim}_{\mathbb{C}} ({\rm Hom}_{{\rm GL}_m (\mathcal{D})} (\pi, \mathds{1})) = 1
$$
\end{theo}

\section{Conclusion.}\label{PartieConclusion}

Grâce aux résultats des théorèmes \ref{CritereDistinctionNRdImpair}, 
\ref{CritereDistinctionNRdPair}, \ref{CritereDistinctionTRPair} et 
\ref{CritereDistinctionTRdImpair}, on a les conditions nécessaires et suffisantes de 
${\rm GL}_m (\mathcal{D})$-distinction suivantes pour les représentations 
cuspidales de niveau $0$ de ${\rm GL}_{\mu} (\Delta)$, images d'une cuspidale par Jacquet-Langlands :

\begin{theo}\label{SyntheseCNSDistinctionCuspidalesNiv0}
Soit $\pi \in \mathcal{R}_0^2 ({\rm GL}_{\mu} (\Delta))$ une cuspidale de niveau $0$, image d'une 
cuspidale de niveau~$0$ de ${\rm GL}_n (\mathbb{K})$ par la correspondance de Jacquet-Langlands, et 
$(\mathbb{K}_n / \mathbb{K}_{\delta}, \chi)$ la paire admissible modérée associée à $\pi$ 
(en particulier $\chi$ est un caractère modéré de $\mathbb{K}_n^{\times}$, et on note $\overline{\chi}$ 
la restriction de $\chi$ à $\mathcal{O}_{\mathbb{K}_n}^{\times}$ vue comme caractère de $k_{\mathbb{K},n}^{\times}$).
\begin{itemize}
\item[$\ast$] Si $\mathbb{K} / \mathbb{F}$ est non ramifiée et $n$ est pair, la représentation $\pi$ n'est pas 
${\rm GL}_m (\mathcal{D})$-distinguée.
\item[$\ast$] Si $\mathbb{K} / \mathbb{F}$ est totalement ramifiée et $n$ est impair, la représentation $\pi$ n'est pas 
${\rm GL}_m (\mathcal{D})$-distinguée.
\item[$\ast$] Si $\mathbb{K} / \mathbb{F}$ est non ramifiée et $n$ est impair. 
Soit $\tau$ un générateur du groupe de Galois ${\rm Gal} (k_{\mathbb{K},n} / k_{\mathcal{D}})$. 
Alors, la représentation $\pi$ est 
${\rm GL}_m (\mathcal{D})$-distinguée si et seulement si $\chi$ est trivial sur $\mathbb{F}^{\times}$ et s'il existe 
$\alpha$ dans le groupe de Galois ${\rm Gal} (k_{\mathbb{K},n} / k_{\Delta})$ tel que 
$\overline{\chi}^{-1} \circ \alpha = \overline{\chi} \circ \tau$.
\item[$\ast$] Si $\mathbb{K} / \mathbb{F}$ est totalement ramifiée et $n$ est pair. 
Soit $l_0$ le corps résiduel de $\mathbb{K}_{n/2}$. On fixe $\varpi_{\mathbb{K}}$ telle que 
$\varpi_{\mathbb{K}}^2 = \varpi_{\mathbb{F}}$ et 
$\eta$ dans $k_{\mathbb{K},n}^{\times} \backslash l_0^{\times}$ tel que $\eta^2 \in l_0^{\times}$. 
Alors, la représentation $\pi$ est 
${\rm GL}_m (\mathcal{D})$-distinguée si et seulement si $\chi$ est trivial sur $\mathbb{F}^{\times}$, 
$\overline{\chi}$ est trivial sur $l_0^{\times}$ et 
$\chi (\varpi_{\mathbb{K}}) \overline{\chi} (\eta) = -1$.
\end{itemize}

\end{theo}

En utilisant ces conditions, on montre que la correspondance de Jacquet-Langlands préserve la distinction 
pour les cuspidales de niveau $0$ au sens suivant :

\begin{theo}
Si $\rho \in \mathcal{R}_0^2 ({\rm GL}_n (\mathbb{K}))$ est une cuspidale (de niveau $0$), alors 
$\rho$ est ${\rm GL}_n (\mathbb{F})$-distinguée si et seulement si 
$JL(\rho)$ est ${\rm GL}_m (\mathcal{D})$-distinguée.
\end{theo}

\begin{démo}
Soit $\rho \in \mathcal{R}_0^2 ({\rm GL}_n (\mathbb{K}))$ une cuspidale de niveau $0$ 
de paire admissible modérée associée $(\mathbb{K}_n, \chi)$.
Soit $\pi = JL(\rho)$. On a vu en \ref{JacquetLanglandsCuspidales} que $\pi$ est 
une représentation cuspidale (de niveau $0$) de ${\rm GL}_{\mu} (\Delta)$ et 
que la paire admissible modérée associée à $\pi$ est aussi $\chi$.
\begin{itemize}
\item[$\ast$] \textbf{Supposons que l'extension $\mathbb{K} / \mathbb{F}$ est non ramifiée.}
D'après le théorème \ref{CritereDistinctionNRdImpair} 
(cas $d=1$), on sait que $\rho$ est ${\rm GL}_n (\mathbb{F})$-distinguée si et seulement si 
$n$ est impair, $\chi$ est trivial sur $\mathbb{F}^{\times}$ et 
$\overline{\chi}^{-1} \simeq_{k_{\mathbb{K}}} \overline{\chi} \circ \tau$ où 
$\langle \tau \rangle = {\rm Gal} (k_{\mathbb{K},n}/ k_{\mathbb{F}})$.
De même, en utilisant les théorèmes \ref{CritereDistinctionNRdImpair} et \ref{CritereDistinctionNRdPair}, 
on sait que $\pi$ est ${\rm GL}_m (\mathcal{D})$-distinguée si et seulement si $n$ est impair 
(donc $d$ est impair), $\chi$ est trivial sur $\mathbb{F}^{\times}$ et 
$\overline{\chi}^{-1} \simeq_{k_{\Delta}} \overline{\chi} \circ \widetilde{\tau}$ où 
$\langle \widetilde{\tau} \rangle = {\rm Gal} (k_{\mathbb{K},n}/ k_{\mathcal{D}})$.
Il nous suffit donc de regarder le cas où $d$ et $n$ sont impairs. 
Dans ce cas, on a le diagramme d'extensions de corps fini suivant:
$$
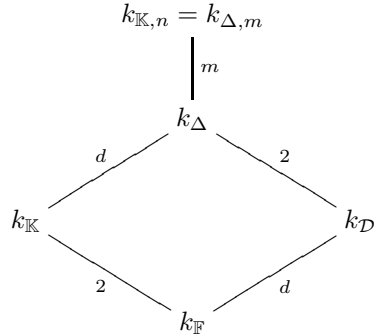
\xymatrix{ & k_{\mathbb{K},n} = k_{\Delta,m} & \\
             & k_{\Delta} \ar@{-}[u]_m & \\ 
         k_{\mathbb{K}} \ar@{-}[ru]^d \ar@{-}[rd]_2 
                & & k_{\mathcal{D}} \ar@{-}[lu]_2 \ar@{-}[ld]^d\\
          & k_{\mathbb{F}} & }
$$
Rappelons que $\overline{\chi}$ est un caractère $k_{\mathbb{K}}$-régulier.\\
Supposons tout d'abord que $\rho$ est ${\rm GL}_n (\mathbb{F})$-distinguée, alors 
$n$ est impair, $\chi$ est trivial sur $\mathbb{F}^{\times}$ et 
$\overline{\chi}^{-1} \simeq_{k_{\mathbb{K}}} \overline{\chi} \circ \tau$. 
Montrons que $\overline{\chi}^{-1} \simeq_{k_{\Delta}} \overline{\chi} \circ \widetilde{\tau}$.
On a $\langle \tau \rangle = {\rm Gal} (k_{\mathbb{K},n}/ k_{\mathbb{F}})$ donc 
$\langle \tau^2 \rangle = {\rm Gal} (k_{\mathbb{K},n}/ k_{\mathbb{K}})$. Il existe $k$ dans $\mathbb{Z}$ 
tel que \mbox{$\overline{\chi}^{-1} = \overline{\chi} \circ \tau \circ \tau^{2k} 
= \overline{\chi} \circ \tau^d \circ \alpha$} avec 
$\alpha = \tau^{2k+1-d}$ et $\tau^d = \widetilde{\tau}$. 
La restriction de $\alpha$ à $k_{\Delta}$, $\alpha_{\vert k_{\Delta}}$, appartient au groupe de Galois 
${\rm Gal} (k_{\Delta} / k_{\mathbb{K}})$ (car $2k+1-d$ est pair). Fixons $\varphi$ un générateur de 
${\rm Gal} (k_{\Delta} / k_{\mathbb{K}})$ et $\widetilde{\varphi}$ un élément de 
${\rm Gal} (k_{\mathbb{K},n} / k_{\mathbb{K}})$ qui prolonge $\varphi$. Il existe alors $r \in \mathbb{Z}$ 
et $\gamma$ dans ${\rm Gal} (k_{\mathbb{K},n} / k_{\Delta})$ tel que 
$\alpha = \widetilde{\varphi}^r \circ \gamma$. Ainsi:
$$
\overline{\chi}^{-1} = \overline{\chi} \circ \delta, \, \, 
\delta = \widetilde{\tau} \circ \widetilde{\varphi}^r \circ \gamma
$$
On en déduit que $\overline{\chi} \circ \delta^2 = \overline{\chi}$. Comme 
$\delta^2 = \widetilde{\tau}^2 \circ \widetilde{\varphi}^{2r} \circ \gamma^2 \in 
{\rm Gal} (k_{\mathbb{K},n} / k_{\mathbb{K}})$, la $k_{\mathbb{K}}$-régularité de $\overline{\chi}$ 
impose que $\delta^2 = Id$. 
Par conséquent $\widetilde{\varphi}^{2r} = \widetilde{\tau}^{-2} \circ \gamma^{-2}$. 
Puisque $\langle \widetilde{\tau}^{2} \rangle = {\rm Gal} (k_{\mathbb{K},n} / k_{\Delta})$, la 
restriction de $\widetilde{\tau}^{-2} \circ \gamma^{-2}$ à $k_{\Delta}$ est l'identité. On 
en déduit que $\widetilde{\varphi}^{2r}_{\vert k_{\Delta}} = \varphi^{2r} = Id$. 
Or $\varphi$ est d'ordre $d$ donc $d$ divise $2r$, et comme $d$ est impair, $d$ divise $r$. 
Ainsi $\varphi^r = Id$ et on peut choisir $\widetilde{\varphi} = Id$.
On déduit de tout ceci que:
$$
\overline{\chi}^{-1} = \overline{\chi} \circ \widetilde{\tau} \circ \gamma 
\simeq_{k_{\Delta}} \overline{\chi} \circ \widetilde{\tau}
$$
Finalement, si $\rho$ est ${\rm GL}_n (\mathbb{F})$-distinguée alors 
$\pi$ est ${\rm GL}_m (\mathcal{D})$-distinguée.\\
Réciproquement, supposons que $\pi$ est ${\rm GL}_m (\mathcal{D})$-distinguée. Alors 
$n$ est impair 
(donc $d$ est impair), $\chi$ est trivial sur $\mathbb{F}^{\times}$ et 
$\overline{\chi}^{-1} \simeq_{k_{\Delta}} \overline{\chi} \circ \widetilde{\tau}$. 
Montrons que $\overline{\chi}^{-1} \simeq_{k_{\mathbb{K}}} \overline{\chi} \circ \tau$.\\
Comme précédemment, on remarque que $\widetilde{\tau} = \tau^d$ et que 
$\langle \tau^{2d} \rangle = {\rm Gal} (k_{\mathbb{K},n} / k_{\Delta})$. 
Il existe $k \in \mathbb{Z}$ tel que:
$$
\overline{\chi}^{-1} = \overline{\chi} \circ \tau^d \circ \tau^{2dk} 
= \overline{\chi} \circ \tau \circ \tau^{d-1 + 2dk}
$$
On a $\tau^{2dk} \in {\rm Gal} (k_{\mathbb{K},n} / k_{\mathbb{K}})$ et, 
puisque $d-1$ est pair, il n'est pas premier avec $2n$, l'ordre de 
${\rm Gal} (k_{\mathbb{K},n} / k_{\mathbb{F}})$, donc engendre un sous-groupe de 
${\rm Gal} (k_{\mathbb{K},n} / k_{\mathbb{F}})$ contenu dans ${\rm Gal} (k_{\mathbb{K},n} / k_{\mathbb{K}})$. 
Par suite, $\tau^{d-1 + 2dk} \in {\rm Gal} (k_{\mathbb{K},n} / k_{\mathbb{K}})$ et:
$$
\overline{\chi}^{-1} \simeq_{k_{\mathbb{K}}} \overline{\chi} \circ \tau
$$
\item[$\ast$] \textbf{Supposons que l'extension $\mathbb{K} / \mathbb{F}$ est totalement ramifiée, 
modérément ramifiée.}
On a directement le résultat en utilisant les théorèmes \ref{CritereDistinctionTRPair} et \ref{CritereDistinctionTRdImpair}.

\end{itemize}

\end{démo}


\vspace{5mm}

\noindent
Charlène Coniglio-Guilloton\\
Charlene.Coniglio@math.univ-poitiers.fr\\
Département de Mathématiques\\
Téléport 2 - BP 30179\\
Boulevard Marie et Pierre Curie\\
86962 Futuroscope Chasseneuil Cedex\\
France

\end{document}